\documentclass[11pt]{article}

\usepackage[tbtags]{amsmath}
\usepackage{amsfonts,amssymb,mathrsfs,amscd,comment}
\usepackage{amsthm}
\usepackage[matrix,arrow,curve]{xy}
\usepackage{longtable}
\usepackage{graphicx}
\def\numberlineI#1{{#1. }}
\def\numberlineII#1{{#1 }}

\def\aaaaa{\let\numberline\numberlineII}
\def\bbbbb{\let\numberline\numberlineI}

\newtheorem{theorem}{Theorem}[section]

\newtheorem{lemma}{Lemma}[section]
\newtheorem{proposition}{Proposition}[section]
\newtheorem{corollary}{Corollary}[section]

\newtheorem{condition}{Condition}
\newtheorem{ncondition}{Condition}

\newtheorem{definition}{Definition}[section]
\newtheorem{remark}{Remark}[section]
\newtheorem{example}{Example}[section]

\newcommand{\llangle}{\langle\kern-.2ex\langle}
\newcommand{\rrangle}{\rangle\kern-.2ex\rangle}
\renewcommand{\mod}{\operatorname{mod}}
\renewcommand{\pmod}[1]{\ (\operatorname{mod}#1)}

\numberwithin{equation}{section}

\newcounter{aa}
\newenvironment{numa}[1]{%
\begin{list}{{\rm(\alph{aa})}}{%
\usecounter{aa}\setlength{\itemindent}{17pt}\setlength{\listparindent}{12pt}
            \setlength{\topsep}{-1pt}
        \setlength{\leftmargin}{0pt}\setlength{\labelsep}{1ex}}}{\end{list}}

\begin{document}
 

\title
{Hill's formula}
\author{Sergey Bolotin\thanks{
Steklov Mathematical Institute, Moscow and  UW-Madison, USA}
and
 Dmitry Treschev\thanks{
Steklov Mathematical Institute, Moscow and
Lomonosov Moscow State University}}

\date{February 25, 2010}

\maketitle

\footnotetext{This research was carried out with the support of the Programme
``Mathematical Control Theory'' of the Presidium of the RAS and
the RFBR (grant no.   08-01-00681-a).}

\begin{abstract}
In his study of periodic orbits of the 3 body problem, Hill obtained
a formula relating the characteristic polynomial of the monodromy
matrix of a periodic orbit and an infinite determinant of the
Hessian of the action functional. A mathematically correct
definition of the Hill determinant and a proof of Hill's formula were
obtained later by Poincar\'e. We give two multidimensional
generalizations of Hill's formula: to discrete Lagrangian systems
(symplectic twist maps) and continuous Lagrangian systems.
We discuss additional aspects which appear in the presence
of symmetries or reversibility. We also study the change of the Morse index of
a periodic trajectory after the reduction of order in a system with symmetries.
Applications are given to the problem of stability of periodic orbits.
\end{abstract}


\setcounter{tocdepth}{2} \tableofcontents

\section{Introduction}
\label{sec1}

In 1886, in his study of lunar orbits, Hill \cite{1}    discovered a formula
which expresses the characteristic polynomial of the monodromy matrix for a
second order time periodic differential equation in terms of the determinant
 of a certain infinite matrix.
Here is a slightly modified version of this result. 
Consider  Hill's equation
\begin{equation}
\label{eq1.1}
\ddot x=a(t)x,
\end{equation}
where
$$
a(t)=\sum_{k=-\infty}^{+\infty} a_k e^{ikt}
$$
is a real       $2\pi$-periodic function. Let       $\rho$ and~$\rho^{-1}$
be  eigenfunctions of the monodromy matrix. Hill showed that
\begin{equation}
\label{eq1.2}
\frac{\rho+\rho^{-1}-2}{e^{2\pi}+e^{-2\pi}-2}=\det H,
\end{equation}                                                 where
   $H$ is the infinite matrix\,\footnote{Hill's matrix was slightly different.}
\begin{equation}
\label{eq1.3}
H=\biggl(\frac{k^2\delta_{jk}+a_{k-j}}{k^2 +1}\biggr)_{j,k\in\mathbb{Z}},
\end{equation} and
  $\delta_{jk}$ is the Kronecker symbol.

Hill computed $\det H$ approximately  replacing $H$ by a $3\times 3$
matrix, which gave  quite a good approximation.  He used equation
\eqref{eq1.2}        to find the multipliers approximately.
Astronomical tables obtained by this method are well-known.

Hill's argument was not rigorous because he did not prove convergence for
the infinite determinant  $\det H$.
Several years later Poincar\'e \cite{2}        explained an exact meaning of
the Hill determinant and presented a rigorous proof of Hill's formula. Hill's  result
 entered textbooks on differential equations, but was almost forgotten by
dynamical systems community until the end of the
XXth century when an analogue
 of equation \eqref{eq1.2}        appeared for discrete Lagrangian systems
in \cite{3} and independently in \cite{4}. Here $H$  turned out to be the
finite Hessian matrix associated with the action functional at the
critical point generated by the periodic solution. In \cite{5} (see also
\cite{6} a general form of     Hill's formula was obtained for a periodic
solution
of an arbitrary Lagrangian system on a manifold. In this case $H$  is a
properly regularized Hessian operator of the action functional at the
critical point determined by a periodic solution.

Both discrete and continuous versions of Hill's formula give non-trivial
information
on the dynamical stability of the periodic orbit in terms of its Morse index.
Recently this connection was
investigated by means of symplectic geometry (see, for example,
 \cite{7} and~\cite{8}). However, the approach
based on the Hill determinant is sometimes simpler and provides additional
insight to the problem.

As mentioned, there are two similar but formally different cases:
\begin{itemize}
\item[--]
Continuous Lagrangian system with configuration manifold $M$ and
$\tau$-periodic Lagrangian
$\mathscr{L}(x,\dot x,t)$ on $TM\times\mathbb{R}$
 which is strictly convex in the velocity.
Solutions of the Lagrangian system will be called trajectories.
 Then $\tau$-periodic trajectories $\gamma$ are critical points of the action
functional
$$
\mathscr{A}(\gamma)=\int_0^\tau \mathscr{L}\bigl(\gamma(t),\dot \gamma(t),t\bigr)\,dt
$$
on the set of $\tau$-periodic curves      $\gamma\colon\mathbb{R}\to M$.

\item[--]
Discrete Lagrangian system with Lagrangian $L(x,y)$ on $M\times M$
satisfying certain non-degeneracy condition. Then periodic
trajectories are $n$-periodic sequences $\mathbf{x}=(x_i)_{i\in\mathbb{Z}}$ which
are critical points of the action functional on $M^n$:
$$
\mathscr{A}(\mathbf{x})=\sum_{i=1}^n L(x_i,x_{i+1}),\qquad
x_{i+n}=x_i.
$$
\end{itemize}

Usually one case can be reduced to the other, but this reduction may be cumbersome.
Hence it makes sense to consider both cases separately.

Both versions of Hill's formula look similar. Let $P$ be the monodromy matrix
 of the periodic trajectory, $h$ the second variation of the action
functional at the periodic trajectory, and $H$ the corresponding  Hessian operator.
Then
\begin{equation}
\label{eq1.4}
\det(P-I)=\sigma(-1)^m\beta\det H,
\end{equation}
where  $m=\dim M$ and~$\sigma=\pm 1$
takes care of orientation.
The coefficient $\beta$ is a positive scaling factor.

The operator $H$ is self-adjoint in a proper Hilbert space.
For continuous systems, $H$ is an unbounded operator, so it needs to
be regularized. For example, for Hill's equation
\eqref{eq1.1}, $H$ is a Sturm--Liouville operator.

Another version of Hill's formula, a generalization of
\eqref{eq1.4}, has the form
\begin{equation}
\label{eq1.5}
\rho^{-m}\det(P-\rho I)=\sigma(-1)^m\beta\det H_\rho,\qquad
\rho\in\mathbb{C},
\end{equation}
where   $H_\rho$ is   the $\rho$-Hessian which coincides with the ordinary
Hessian for $\rho=1$. It is self-adjoint if $|\rho|=1$.  Since $P$
is symplectic, both sides of  \eqref{eq1.5} are polynomials of degree $m$ in
$\rho+\rho^{-1}$.

Hill's formula  \eqref{eq1.4} has many dynamical
applications. The first one is the well known statement that the
Poincar\'e degeneracy of a periodic trajectory (that is, the condition that
 1 is an eigenvalue of $P$) is equivalent to
the variational degeneracy (the condition $\det H=0$).

Another application concerns dynamical instability of a periodic trajectory.
It is based on the observation that
the inequality $\det(P-I)<0$ implies the existence of a real multiplier
(that is,  an eigenvalue of $P$) $\rho>1$.
Indeed, $F(\rho)=\det(P-\rho I)=\det(\rho I-P)\to\infty$ as $\rho\to+\infty$, and
so $F(1)<0$ implies the existence
of a root $\rho>1$.  Thus $\gamma$ has a positive Lyapunov exponent and is
exponentially unstable.

If $\det H\ne 0$, we have \
 $\operatorname{sign}\det H=(-1)^{\operatorname{ind} H}$,
where  $\operatorname{ind} H$
is the Morse index of the periodic trajectory. Hence if the periodic
trajectory is nondegenerate, then by \eqref{eq1.4}
the inequality $\sigma(-1)^{m+\operatorname{ind} H}<0$ implies
exponential instability in a `physical' (with $\beta>0$) system.

In some cases it is possible to prove that for any $|\rho|=1$ the Hessian
$H_\rho$ is positive definite and therefore the equation $\det H_\rho = 0$
has no solutions on the unit circle. Then we obtain exponential
instability, in fact, total hyperbolicity for the corresponding periodic
trajectory (Propositions \ref{pr4.5} and \ref{pr7.2}).

Below we present other dynamical consequences of Hill's formula.

Note that  the connections between dynamical and geometrical
properties of periodic orbits are not restricted to     Hill's
formula. We mention here interesting relations  between stability
properties and the structure (index, signature, and so on) of a quadratic
first integral of the linearized system (\cite{9} and  \cite{10}). Many
interesting results follow from the index formula in symplectic
geometry (\cite{7} and  \cite{8}). Some of our results may be
regarded as Lagrangian versions of the results of
\cite{11} and  \cite{12}.

Hill's formula is potentially most useful for the study of periodic
orbits obtained by variational methods. Many such orbits were
obtained recently in celestial mechanics by minimization of the
action functional on appropriate classes of curves, see, for example,
\cite{13}--~\cite{15}. The most famous example is
the figure eight orbit, see \cite{13}.        However, due to
rotational and other symmetries, none of these periodic orbits are
nondegenerate minimum points of the action.

In applications periodic trajectories are usually degenerate.  For
example, any periodic orbit of an autonomous  continuous Lagrangian
system is degenerate.  In this case the variational equation has a
$\tau$-periodic solution, $\dot\gamma$, and a linear first integral,
the linearization of the energy integral. Another reason for such a
degeneracy (now in both discrete and continuous cases) is the
presence of a symmetry group, preserving the Lagrangian. This
degeneracy also gives $\tau$-periodic solutions and linear integrals
for the variational equation. For degenerate periodic trajectory
equation \eqref{eq1.4}     is useless because both sides vanish. A
nondegenerate version of Hill's formula can be obtained with the help of
the reduction procedure. We consider the case when the Lie algebra $V$ of
symmetry vector fields for the variational equation is commutative, the
dimension of the generalized unit eigenspace $N$ of $P$ is $2k$, where $k=
\dim V$ (no further degeneracy) and a condition, called
the non-degeneracy of the trajectory $\mod V$, holds. The latter condition has a
Lagrangian nature rather than Hamiltonian.

The reduced Hill's formula looks similarly, but the corresponding monodromy
and Hessian operators $\tilde P$ and $H^\perp$ act on smaller (reduced) spaces, and
are nondegenerate if all the symmetries are taken into account:
$$
\det(\widetilde P-I)=\sigma^\perp(-1)^{m-k}\beta^\perp\det H^\perp.
$$
Here
      $\sigma^\perp\in\{1,-1\}$ is the `reduced orientation' and
  $\beta^\perp>0$.

Now an interesting question appears on the relation between $\sigma$ and
$\sigma^\perp$ as well as between $\operatorname{ind} H$ and
  $\operatorname{ind} H^\perp$, the Morse indices of the Hessians in the
original and reduced systems. Indeed, $\sigma$ and $\operatorname{ind} H$
are often known for solutions obtained by variational methods, while
$\sigma^\perp$ and $\operatorname{ind} H^\perp$ appear in stability problems.
The following construction explains our answer to this question.

Let $h$ and $h^\perp$ be bilinear second variation forms corresponding to the
operators $H$ and $H^\perp$ respectively. The forms $h$ and $h^\perp$ are defined
on the vector spaces $X$ and $X^\perp$ of variations along the periodic orbit, for
the original and reduced Lagrangian system, respectively. The procedure of
the order reduction gives a canonical projection $\Pi\colon  X\to X^\perp$.

For any $\zeta\in V$ and $\eta\in X$ we have $h(\zeta,\eta) = 0$.
Therefore $h$ defines a bilinear form $\hat h$ on $\widehat X = X / V$ and
$\operatorname{ind} h=\operatorname{ind} \hat h$. The spaces $\widehat X$ and
$X^\perp$ admit the expansions
$$
\widehat X=\widehat\Omega\oplus\widehat Y^0\oplus\widehat Z,\qquad
\widehat X^\perp=\Omega^\perp\oplus Y^\perp
$$
with the following properties:
\begin{itemize}
\item[1)] $\dim\widehat\Omega=\dim\Omega^\perp=k$;

\item[2)] the spaces  $\widehat\Omega$, $\widehat Y^0$, and
$\widehat Z$ are $\hat h$-orthogonal, while~$\Omega^\perp$ and
  $Y^\perp$ are $h^\perp$-orthogonal;

\goodbreak
\item[3)]
the
restriction $\hat h\big|_{\widehat Z}$ is nondegenerate in the discrete
case and positive definite in the continuous case;

\item[4)] $\Pi(\widehat Z\,)=0$, while the restrictions
$\Pi\big|_{\widehat Y^0}\colon \widehat Y^0\to Y^\perp$ and
$\Pi\big|_{\widehat\Omega}\colon \widehat\Omega\to \Omega^\perp$
are linear isomorphisms;

\item[5)] the forms $\hat h\big|_{\widehat Y^0}$ and~$h^\perp\big|_{Y^\perp}$
coincide in the sense that  $h^\top\big|_{Y^\perp}=h^\perp\big|_{Y^\perp}$,
where  $\hat h=h^\top\circ\Pi$;

\item[6)] $h^\top\big|_{\Omega^\perp}-h^\perp\big|_{\Omega^\perp}=\chi$,
where~$\chi$  is positive definite in the continuous case.
\end{itemize}

In a convenient basis we give an explicit expression for the matrices
$h^\top\big|_{\Omega^\perp}$ and~$h^\perp\big|_{\Omega^\perp}$.
Using these expressions we show that
$$
\sigma (-1)^{\operatorname{ind} H}=\sigma^\perp (-1)^{\operatorname{ind} H^\perp+\operatorname{ind} b},
$$
where the
quadratic form  $b$ on the generalized eigenspace
$N=\operatorname{Ker}(P-\nobreak I)^2$ is defined by $b(v)=\omega((P-I)v,v)$,
where $\omega$ is the symplectic structure.      

In some cases $(-1)^{\operatorname{ind} b}$ has a clear dynamical meaning.
For example, suppose that the degeneracy appears solely because the continuous
Lagrangian
is autonomous. Then $\dim V = 1$. The periodic trajectory $\gamma$ belongs to
a smooth family of periodic trajectories. Let $E$ and $\tau$ be
the energy and the period along this family. Then (see Lemma \ref{lem6.4})
$$
(-1)^{\operatorname{ind} b}=-\operatorname{sign} \frac{dE}{d\tau}\,.
$$
Suppose that the periodic trajectory $\gamma$ of an autonomous Lagrangian
system has no other degeneracy. Then $k=1$ and by the reduced Hill's formula
it has a real multiplier $\rho > 1$ provided that
\begin{equation}
\label{eq1.6}
\sigma (-1)^{m+\operatorname{ind} H} \frac{dE}{d\tau}<0.
\end{equation}

The sign of the quantity $dE / d\tau$ is easily computed, for example, in
the problem of the motion of a point in $\mathbb{R}^m$  in a homogeneous
potential force field.

It turns out that closed geodesics do not satisfy the condition of
non-degeneracy $\mod V$. However we show that inequality \eqref{eq1.6}
 still implies the existence of a multiplier $\rho > 1$ provided no extra
degeneracy takes place (Corollary \ref{cor5.2}). (Note that in this case
$dE / d\tau < 0$).

As mentioned above, the subject of this paper is closely related to the
theory of Maslov--Morse index for periodic orbits of Hamiltonian systems, see,
for instance,  \cite{7},~\cite{8},~\cite{11},~\cite{12}).
Some of our results can be obtained by these purely symplectic methods. Others
are Lagrangian, and so do not have direct symplectic formulation. The
situation is similar to the relation between Hamiltonian and Lagrangian
systems: Hamiltonian theory is simpler, more general, and more powerful. Nevertheless
for many problems the Lagrangian approach is essential.

This paper splits in two parts: discrete and continuous. Although the majority
of constructions
and statements in the discrete and continuous parts are analogous, there are
many technical differences
which forced us to deal with these two cases separately.

The plan of the paper is as follows. In \S\,\ref{sec2} we first recall the
definition and basic properties of discrete Lagrangian systems (DLS). This
material is well known to specialists, but these objects are not as standard as
their continuous analogues.

Then we present several versions of Hill's formula for a periodic trajectory of
a DLS. As an application, we give some sufficient conditions for the
instability of  periodic trajectories. Several statements concern stability problem
 for billiard systems in arbitrary dimension. For example, any $n$-periodic trajectory
$\mathbf{x}$ of a billiard  system
inside a hypersurface in  $\mathbb{R}^{m+1}$ such that
    $(-1)^{m+n+\operatorname{ind}\mathbf{x}}<0$
is exponentially unstable (Corollary \ref{cor2.6}). As far as we know,
 there are very few publications about stability of periodic trajectories in multidimensional
DLS. Here we mention \cite{16}     and \cite{17}, where trajectories of
period 2 are studied.

In \S\,\ref{sec3} we consider DLS with symmetry. We present a discrete version
 of Routh's procedure of order reduction and a reduced version of Hill's
formula where the degeneracy which appears due to symmetry is removed. We
also give a formula for the difference  between the Morse index of a periodic trajectory
of the original system and the Morse index of the
corresponding periodic orbit of the reduced system.

In \S\,\ref{sec4} we study reversible DLS,
that is, discrete Lagrangian systems with the Lagrangian $L$ invariant under
time reversal combined with an involution $S\colon M\to M$,  $S^2=\mathrm{id}$.
 Thus $L(S(x),S(y))=L(y,x)$ for any $x,y\in M$. For any trajectory
 $\mathbf{x}=(x_i)$ of the DLS, the sequence
$\widetilde{\mathbf{x}}=(Sx_{-i})$ is also a trajectory. If
$\mathbf{x}=\widetilde{\mathbf{x}}$
modulo  a translation, the trajectory $\widetilde{\mathbf{x}}$ is called
$S$-reversible. Then the corresponding space of variation splits into a
direct sum of spaces of odd and even variations with respect to $S$. Hill's determinant
also admits  splitting into a product of two determinants.  
Reversible periodic trajectories are also critical
points of another action functional $\mathscr{A}_+$ which is obtained from
the original one, $\mathscr{A}$, by restriction to the space of even
variations. Morse index of a trajectory with respect to $\mathscr{A}_+$ is in
general different from that computed with respect to $\mathscr{A}$.

Any $S$-reversible trajectory $\mathbf{x}$ of a DLS has 0,1, or 2
fixed points of the involution $S$. According to this we say that
$\operatorname{type}\mathbf{x}= 0,1$, or 2. One of application, presented in
\S\,\ref{sec4} is as follows (Corollary \ref{cor4.4}).
Suppose that the billiard surface $M\subset\mathbb{R}^{m+1}$ is symmetric relative to
 a hyperplane and $S$ denotes this symmetry. Let $\mathbf{x}$ be an
$S$-reversible periodic billiard trajectory of type $\tau\in\{0,1,2\}$ which
is a nondegenerate minimum of the `half-length' $\mathscr{A}_+$. If
$m+\tau$ is odd, then $\mathbf{x}$ is exponentially unstable.

In \S\,\ref{sec5} the continuous part of the paper starts. The main
technical difference of the continuous case is the infinite dimension of the
space of variations. Because of this the definition of
Hill's determinant needs more care. We give  a construction defining
the Hill determinant and present several versions of Hill's formula analogous
to the ones in the discrete case. Then we give applications to
instability of periodic orbits of Lagrangian systems including the case of
closed geodesics. A typical statement from this part
(in fact, going back to Poincar\'e) is as follows. Let $\gamma$ be a
nondegenerate closed geodesic on an $m$-dimensional manifold and
$\sigma (-1)^{m+\operatorname{ind}\gamma} > 0$. Then $\gamma$ is
exponentially unstable (Corollary \ref{cor5.3}).

In \S\,\ref{sec6} we discuss the role of  symmetries and give  a version
of Hill's formula
which eliminates the corresponding degeneracy. Then we study the relation
between the Morse index of the periodic trajectory of the original system and
the corresponding periodic solution of the reduced system. We present some
applications of this formula to the problem of stability for
Lagrangian systems with symmetry.

Finally, in \S\,\ref{sec7} we consider a reversible CLS.    The
Lagrangian~$\mathscr{L}$ of an $S$-reversible NLS is compatible with the
involution $S$ in the following sense:
$$
\mathscr{L}\bigl(S(x),dS(x)\dot x,t\bigr)=\mathscr{L}(x,-\dot x,-t).
$$
As in the discrete case, the functional $\mathscr{A}_+$ corresponding to even
variations is defined. The main questions are the relation between
the indices of an $S$-reversible periodic trajectory with respect to
$\mathscr{A}$ and $\mathscr{A}_+$ and the relation between the index
with respect to $\mathscr{A}_+$ and stability properties. We show that in
many cases the computation of $\operatorname{ind}\gamma\;\bmod 2$ may be
performed on variations from a  $2m$-dimensional space.

The authors are grateful to V.\,V.~Kozlov for very useful discussions.

\section{Discrete case}
\label{sec2}

\subsection{Discrete Lagrangian systems (DLS)}
\label{ssec2.1}
Let  $M$ be an $m$-dimensional manifold and  $L$ a
smooth\,\footnote{Actually, $C^2$ is enough.} function on
  $M^2=M\times M$.
Denote
\begin{equation}
\label{eq2.1}
\partial_1 L(x,y)=\frac{\partial L(x,y)}{\partial x}\,,\qquad
\partial_2 L(x,y)=\frac{\partial L(x,y)}{\partial y}
\end{equation}
and let
$$
B(x,y)=-\partial_1\partial_2L(x,y).
$$
In local coordinates,
\begin{equation}
\label{eq2.2}
B(x,y)=-\biggl(\frac{\partial^2 L}{\partial y_j\,\partial x_i}\biggr).
\end{equation}

In invariant terms,     $B(x,y)$ is a linear operator $T_xM\to T_y^*M$,
or a bilinear form on  $T_xM\times T_yM$. We say that~$L$
is a discrete Lagrangian if it satisfies the following condition.

\medskip\noindent{\bf
 Twist condition.}
$B(x,y)$ is nondegenerate for all  $x,y\in M$.
 \medskip

Any discrete Lagrangian  $L$ locally defines a map
$$
T\colon M^2\to M^2,\qquad
T(x,y)=(y,z),
$$
where     $z=z(x,y)$ is determined by the equation
\begin{equation}
\label{eq2.3}
\frac{\partial}{\partial y}\bigl(L(x,y)+L(y,z)\bigr)=
\partial_2 L(x,y)+\partial_1 L(y,z)=0.
\end{equation}

In general,  $T$ is a multivalued map (relation) with the graph
$$
\Gamma=\bigl\{(x,y,y,z)\in M^2\times M^2:
\partial_2 L(x,y)+\partial_1 L(y,z)=0\bigr\}.
$$
The dynamical system determined by $T$ is called the discrete Lagrangian
system (DLS) with configuration space $M$ and Lagrangian $L$.

\begin{remark}
\label{rem2.1}{\rm
In this paper we deal with a small neighbourhood of a periodic orbit.
Hence it is sufficient to assume that the non-degeneracy condition holds locally.}
\end{remark}

It is easy to check (see, for example, \cite{18}) that $T$ is
 symplectic with respect to the symplectic 2-form $\omega=B(x,y)\,dx\wedge dy$,
\begin{equation}
\label{eq2.4}
\omega(\mathbf{u},\mathbf{v})=\langle B(x,y)u_1,v_2\rangle-
\langle B(x,y)v_1,u_2\rangle,\qquad
\mathbf{u}=(u_1,u_2),\quad
\mathbf{v}=(v_1,v_2)
\end{equation}
(${\langle\,\cdot\,{,}\,\cdot\,\rangle}$ is the canonical
pairing of a covector on a vector). 

\begin{remark}
\label{rem2.2}{\rm
Let us pass to Hamiltonian variables by the map
$S\colon M^2\to T^*M$, \,$(x,y)\mapsto (x,p_x)$, \,$p_x=-\partial_1 L(x,y)$.
It is locally
invertible and replaces $T$ by a locally defined map
$F=STS^{-1}\colon T^*M\to T^*M$. The map $F$  is symplectic with respect
to the standard symplectic form $dp_x\wedge dx$ on $T^*M$, and $L$
is the generating function of $F$:
$$
F(x,p_x)=(y,p_y),\qquad
p_x=-\partial_1 L(x,y),\quad p_y=\partial_2
L(x,y).
$$
Such a symplectic map $F$ is usually called a twist map.}
\end{remark}

The map $T$ remains the same after multiplication of the Lagrangian by a constant,
after addition of a constant to $L$, and after the so-called gauge
transformation
$$
L(x,y) \mapsto L(x,y)+f(x)-f(y)
$$
with an arbitrary smooth function $f$ on~$M$.

A typical example of DLS is the multidimensional standard map:
\begin{equation}
\label{eq2.5}
L(x,y)=\frac{1}{2}\,\bigl\langle B(x-y),x-y\bigr\rangle
-\frac{1}{2}\,\bigl(V(x)+V(y)\bigr),\qquad
x,y\in \mathbb{R}^m,
\end{equation}
where $B$ is a symmetric constant nondegenerate  matrix.\footnote{One can
replace the potential $(V(x)+V(y))/2$ by $V(x)$ or $V(y)$ because they
are all gauge-equivalent.}

Consider a  domain  in  $\mathbb{R}^{m+1}$  bounded by a smooth convex
hypersurface $M$.  The billiard system  is a DLS with the Lagrangian
$L(x,y) = |x-y|$ on $M\times M$. Let $\langle B(x,y)v,w\rangle$ be the
bilinear form on $T_xM\times T_yM$ corresponding to the operator
$B(x,y):T_xM\to T_y^*M$. A computation gives
\begin{equation}
\label{eq2.6}
\langle B(x,y)v,w\rangle=\frac{\langle v,w\rangle-
\langle v,e\rangle \langle w,e\rangle}{|x-y|}\,,\qquad
e=\frac{x-y}{|x-y|}\,.
\end{equation}
We may identify  $T_xM$ and~$T_yM$ by an isomorphism
$\Pi(x,y)\colon T_x M \to T_y M$, which is the parallel projection in
  $\mathbb{R}^{m+1}$ along the segment $[x,y]$: $\Pi v
   =    
v\allowbreak
\pmod e$. Then
$$
\langle B(x,y)v,\Pi(x,y)v\rangle
=\frac{|v|^2-\langle v,e\rangle^2}{|x-y|}>0,\qquad
v\in T_xM\setminus\{0\}.
$$
We orient  $M$ as the boundary. Since       $\Pi(x,y)$ changes orientation,
we obtain

\begin{proposition}
\label{pr2.1}
$\det B(x,y)<0$.
\end{proposition}

Since the image and the range of $B(x,y)$ are different, $\det B(x,y)$ is not
invariantly defined, but its sign is.
The fact that the map $B(x,y)$ is  nondegenerate, provided the
hyperplanes $T_x M$ and $T_y M$ are not parallel to each other in
$\mathbb{R}^m$, is well-known; for a recent reference see \cite{19}.

In \cite{18} the reader can find many examples of  (mostly integrable)  DLS,
including multidimensional ones. 

For a continuous Lagrangian system (CLS) with Lagrangian
$\mathscr{L}(x,\dot x)$, an analogue of the operator $B(x,y)$ is the matrix
$\mathscr{L}_{\dot x\dot x}(x,\dot x)$ of second partial derivatives.
Indeed, consider a DLS on $\mathbb{R}^m$ with the Lagrangian
$L(x,y)=\mathscr{L}(x,(y-x)/\varepsilon)$. In the limit as $\varepsilon\to 0$,
orbits of DLS converge to orbits of the CLS with the Lagrangian $\mathscr{L}$.
 A computation shows that
$$
\varepsilon^2 B(x,y)=\mathscr{L}_{\dot x\dot x}(x,(y-x)/\varepsilon)+
O(\varepsilon).
$$
In particular, for an analogue of a positive definite Lagrangian system,
$\det B(x,y)>0$.

For $m\ge 2$ there is no universally accepted discrete analogue of
positive definite continuous Lagrangian systems. Indeed, in general $B$ is
not symmetric and, moreover, its symmetry does not have an invariant
meaning since $B$ and $B^*$ are defined on different spaces.
Note that the 1-dimensional Aubry--Mather theory was developed
for twist maps, while multidimensional theory is well developed
 for continuous positive definite Lagrangian systems.

The most common definition of a positive definite DLS is as follows.
Let $M=\mathbb{R}^m$ and suppose $L$ satisfies the following conditions
(see, for example, \cite{20}): 

\begin{itemize}
\item[--] the function  $L(x,x+v)$ is periodic in  $x\in\mathbb{T}^m$
and superlinear
in $v\in\mathbb{R}^m$;
\item[--] for any    $x\in\mathbb{R}^m$ the map
$y\mapsto \partial_1 L(x,y)$ is a diffeomorphism of $\mathbb{R}^m$.
\end{itemize}

Then $L$ is a generating function of a globally defined symplectic twist
map                of
$\mathbb{T}^m\times\mathbb{R}^m$. Evidently, such
  $L$ satisfies $\det B(x,y)>0$.

\subsection{Discrete Hill determinant}
\label{ssec2.2}
Let   $(x_i)_{i\in\mathbb{Z}}$, \,$x_{i+n}=x_i$, be an
$n$-periodic trajectory of a DLS, that is,
        $T(x_{i-1},x_i)=(x_i,x_{i+1})$ for all  $i$. The periodic orbit is
determined by       $\mathbf{x}=(x_1,\dots,x_n)\in M^n$, and a cyclic permutation
of  $\mathbf{x}$ gives the same orbit.
By       \eqref{eq2.3},
\begin{equation}
\label{eq2.7}
\partial_2 L(x_{i-1},x_i)+\partial_1 L(x_i,x_{i+1})=0, \qquad
i=1,\dots,n,
\end{equation}
where  $x_0=x_n$ and~$x_1=x_{n+1}$. Thus,          $\mathbf{x}$ is a
critical point of the action functional
$$
\mathscr{A}({\mathbf{x}})=L(x_1,x_2)+L(x_2,x_3)+\dots+L(x_n,x_1),\qquad
\mathbf{x}\in M^n.
$$

The point $p=(x_1,x_2)$ is a fixed point of the map $T^n\colon M^2\to
M^2$.   The linear approximation to dynamics of $T$ near the periodic
trajectory is
determined by the linear Poincar\'e map $P=DT^n(p)\colon W\to W$,
$W=T_{p}M^2$. In local coordinates, $P$ becomes the monodromy
matrix defined uniquely up to a similarity $P\mapsto S^{-1}PS$.
Eigenvalues of $P$ are called multipliers of the periodic orbit.
They determine dynamical properties of the periodic
trajectory in the linear approximation.

Let
$$
\mathbf{H}=\frac{\partial^2 \mathscr{A}(\mathbf{x})}{\partial \mathbf{x}^2}
$$
be the Hessian matrix of  $\mathscr{A}$ at the critical point~$\mathbf{x}$.
Denote
$$
B_{i}=B(x_i,x_{i+1}),\qquad
x_{n+1}=x_1.
$$

\begin{theorem}[{\rm discrete Hill formula}]
\label{th2.1}
\begin{gather}
\label{eq2.8}
\det(P-I)=\frac{(-1)^{m}\det\mathbf{H}}{\prod_{i=1}^n \det B_{i}}
=\sigma (-1)^m \beta \det\mathbf{H},
\\
\label{eq2.9}
\sigma(\mathbf{x})=\operatorname{sign}\prod_{i=1}^n \det B_{i},\qquad
\beta=\biggl|\,\prod_{i=1}^n \det B_{i}\biggr|^{-1}.
\end{gather}
\end{theorem}

For `physical' discrete Lagrangians  the geometrical meaning of
$\sigma$ is the orientability: the trajectory $\mathbf{x}$ is, in a certain
sense, orientable if $\sigma(\mathbf{x})>0$ and non-orientable
otherwise. For example, this is true if DLS is obtained by discretization
of  a positive definite CLS. By Proposition \ref{pr2.1},
 for a billiard $n$-periodic trajectory $\mathbf{x}$\; $\sigma(\mathbf{x})
= (-1)^n$. Therefore in this sense billiard periodic trajectories with
odd period are non-orientable. Note that
$\sigma$ is replaced by $(-1)^{mn}\sigma$ if we replace $L$ by $-L$.

\subsection{Invariant meaning of Hill's formula}
\label{ssec2.3}
The left-hand side of \eqref{eq2.8}      obviously does not depend
on the choice of local coordinates in $M$. However an invariant
meaning of the right hand side is a priori not clear. Let us explain
why   it is coordinate independent.  Let $E_i=T_{x_i}M$. Then
$B_i=B(x_i,x_{i+1})$  is a linear operator $E_i\to E_{i+1}^*$, and
$$
  A_i = \partial_{22}L(x_{i-1},x_i) + \partial_{11}L(x_i,x_{i+1})
$$
is  a symmetric operator $A_i\colon E_i\to E_i^*$.

The Hessian   of $\mathscr A$ at the critical
point $\mathbf{x}\in M^n$ is a
symmetric bilinear form $h$ on $X=T_\mathbf{x} M^n=E_1\times\dots\times
E_n$ given by
\begin{equation}
\label{eq2.10}
h(\mathbf{u},\mathbf{v})=\sum_{i=1}^n \bigl(\langle A_iu_i,v_i\rangle
-\langle B_{i-1}u_{i-1},v_i\rangle-\langle B_{i}^*u_{i+1},v_i\rangle\bigr),
\end{equation}
where
$$
\mathbf{u}=(u_1,\dots,u_n),\quad
\mathbf{v}=(v_1,\dots,v_n),\quad
u_0=u_n,\quad
u_{n+1}=u_1.
$$

The form~$h$ is represented by a symmetric operator
$\mathbf{H}\colon X\to X^*$:
$$
h(\mathbf{u},\mathbf{v})=\langle \mathbf{H}\mathbf{u},\mathbf{v}\rangle,
\qquad
\mathbf{u},\mathbf{v}\in X,
$$
where
$$
(\mathbf{H}\mathbf{u})_i=A_iu_i-B_{i-1}u_{i-1}-B_i^*u_{i+1},\qquad
i=1,\dots,n.
$$

Define linear operators      $\mathbf{A},\mathbf{B}\colon X\to X^*$
by
\begin{equation}
\label{eq2.11}
(\mathbf{A}\mathbf{u})_i=A_iu_i,\quad
(\mathbf{B}\mathbf{u})_i=-B_{i-1}u_{i-1},\quad
(\mathbf{B}^*\mathbf{u})_i=-B_i^*u_{i+1}.
\end{equation}
Then
$$
\mathbf{H}=\mathbf{A}+\mathbf{B}+\mathbf{B}^*.
$$

Since the~$B_i$ are nondegenerate, $\mathbf{B}$ is invertible. If we introduce
local coordinates, then
$\mathbf{B}$ becomes an $(mn\times mn)$-matrix, and
\begin{equation}
\label{eq2.12}
\det\mathbf{B}=(-1)^m \prod_{i=1}^n \det B_{i}.
\end{equation}

Hence Hill's formula takes the invariant form\,\footnote{Since
$\det\mathbf{H},\det\mathbf{B}$ are linear operators of 1-dimensional spaces
 $\wedge^{mn}X\to \wedge^{mn} X^*$, their quotient is a well defined scalar.}
\begin{equation}
\label{eq2.13}
\det(P-I)=\frac{\det\mathbf{H}}{\det\mathbf{B}}=
\det(\mathbf{B}^{-1}\mathbf{H}).
\end{equation}

The equation  $\mathbf{H}\mathbf{u}=0$ gives the variational system of the
periodic trajectory  $\mathbf{x}$:
\begin{equation}
\label{eq2.14}
A_iu_i-B_{i-1}u_{i-1}-B_i^*u_{i+1}=0,\qquad
u_i\in E_i,\quad
i\in\mathbb{Z}.
\end{equation}
This is the linear approximation to the system~\eqref{eq2.7}
near the periodic trajectory   $\mathbf{x}$. More precisely, if $u_i$
is any solution of the variational system, then the linearized map
                               $P_i=dT(x_{i-1},x_i)$ acts as
$$
P_i(u_{i-1},u_{i})=(u_{i},u_{i+1}),\qquad
u_{i+1}=(B_i^*)^{-1}(A_iu_i-B_{i-1}u_{i-1}).
$$
The kernel of  $\mathbf{H}\colon X\to X^*$ is the set of
              $n$-periodic solutions    $\mathbf{v}=(v_i)$, \,$v_{i+n}=v_i$,
of
  \eqref{eq2.14}.

The variational system is a linear Lagrangian system.

\begin{definition}
\label{def2.1}
A linear periodic discrete Lagrangian system $(E,\Lambda)$
is defined by $n$-periodic sequences $E=(E_i)_{i\in Z}$  of vector spaces
and linear operators $A_i\colon E_i\to E_i^*$,
$B_i \colon E_i\to E_{i+1}^*$, where $A_i$ is symmetric and $B_i$ is
nondegenerate. The Lagrangian is
\begin{equation}
\label{eq2.15}
\Lambda_i(u_i,u_{i+1})=\frac{1}{2}\,\langle A_iu_i,u_i\rangle-
\langle B_{i}u_{i},u_{i+1}\rangle.
\end{equation}
Trajectories of    $(E,\Lambda)$ are sequences $\mathbf{u}=(u_i)$ such that
$$
\partial_{u_i}\bigl(\Lambda_{i-1}(u_{i-1},u_i)+\Lambda_i(u_i,u_{i+1})\bigr)=0.
$$
\end{definition}

Thus trajectories of $(E,\Lambda)$  satisfy the variational system
\eqref{eq2.14}
and are extremals of the quadratic action functional
\begin{equation}
\label{eq2.16}
\frac{1}{2}\,h(\mathbf{u},\mathbf{u})=\sum_{i=1}^n \Lambda_i(u_i,u_{i+1}).
\end{equation} The
system  $(E,\Lambda)$ is the linearization of $(M,L)$ at~$\mathbf{x}$.

\subsection{Generalized Hill determinant}
\label{ssec2.4}
Let us define a generalization of the Hessian $\mathbf{H}$.
 Let $S^1=\{\rho\in\mathbb{C}:|\rho|=1\}$.
For any $\rho\in S^1$, let $X_\rho$ be the space of all
quasiperiodic complex vector sequences $\mathbf{u}=(u_j)_{j\in\mathbb{Z}}$
 such that $u_{j+n}=\rho u_j$.
Here $u_j$ lies in the complexification of $E_j$ which we will denote
$E_j$ for simplicity. The Hessian of the
action defines a Hermitian form  on $X_\rho$:
$$
h(\mathbf{u},\overline{\mathbf{v}})=\sum_{j=1}^n
\bigl(\langle A_j u_j,\bar{v}_j\rangle-
\langle B_{j-1}u_{j-1},\bar{v}_j\rangle
-\langle B_{j}^*u_{j+1},\bar{v}_j\rangle\bigr).
$$

Since a quasiperiodic sequence $(u_i)_{i\in\mathbb{Z}}$ is determined by
$\mathbf{u}=(u_1,\dots,u_n)\in
E_1\times\dots\times E_n=X$, we identify $X_\rho$ with
$X$ (more precisely, with the complexification of $X$). Then we obtain
the  Hermitian form
$$
h_\rho(\mathbf{u},\overline{\mathbf{v}})= \langle
\mathbf{H}_\rho\mathbf{u},\overline{\mathbf{v}}\rangle,\qquad
\mathbf{u},\mathbf{v}\in X,
$$
where  $\mathbf{H}_\rho\colon X\to X^*$ is given by
\begin{gather*}
(\mathbf{H}_\rho\mathbf{u})_j=A_ju_j-B_{j-1}u_{j-1}-B_j^*u_{j+1},
\\
j=1,\dots,n,\quad
u_0=\rho^{-1}u_n,\quad
u_{n+1}=\rho u_1.
\end{gather*}
Similarly we define an operator  $\mathbf{B}_\rho\colon X\to X^*$:
$$
(\mathbf{B}_\rho\mathbf{u})_{j}=-B_{j-1}u_{j-1},\qquad
j=1,\dots,n,\quad
u_0=\rho^{-1}u_n.
$$
Then
$$
\mathbf{H}_\rho=\mathbf{A}+\mathbf{B}_\rho+\mathbf{B}_\rho^*.
$$

The operators  $\mathbf{H}_\rho$, $\mathbf{B}_\rho$ make sense
for any non-zero    $\rho\in\mathbb{C}$.

\begin{theorem}[{\rm generalized Hill formula}]
\label{th2.2}
For any    $\rho\in\mathbb{C}$
\begin{equation}
\label{eq2.17}
\det(P-\rho I)=\frac{\det \mathbf{H}_\rho}{\det \mathbf{B}_\rho}\,.
\end{equation}
\end{theorem}

Since   $\det\mathbf{B}_\rho=\rho^{-m}\det\mathbf{B}$, we obtain
\begin{equation}
\label{eq2.18}
\rho^{-m}\det(P-\rho I)=\frac{\det\mathbf{H}_\rho}{\det\mathbf{B}}\,.
\end{equation}

This is an invariant version of the result of \cite{3}.
For $\rho=1$, \eqref{eq2.18}  gives
\eqref{eq2.8}.

Both sides in \eqref{eq2.18}  are polynomials of degree $m$ in the
Hill   {\it discriminant} $\rho+\rho^{-1}$ with senior coefficient 1.
Indeed, the characteristic polynomial $F(\rho)=\det(P-\rho I)$ of the
symplectic operator $P$  satisfies
$F(\rho)=\rho^{2m}F(\rho^{-1})$.  Hence $\rho^{-m}F(\rho)$
is a symmetric polynomial in $\rho$ and $\rho^{-1}$.
Thus it is a function of $\rho+\rho^{-1}$.

\goodbreak
In coordinates,           $\mathbf{H}_\rho$ is an $(mn\times mn)$-matrix
which coincides with $\mathbf{H}$ with two exceptions:
in the upper right $(m\times m)$-block, $-B_n$ is replaced by $-\rho^{-1} B_n$
and in the lower
left           $(m\times m)$-block,  $-B_n^*$ is replaced by  $-\rho B_n^*$.

\medskip\noindent{\it
Proof Theorem~\ref{th2.2}.}
Let us show that  $G(\rho)=\det(\mathbf{B}_\rho^{-1}\mathbf{H}_\rho)$ is
a polynomial of degree~$2m$ with
senior coefficient  equal to  1.

Let     $\nu=\rho^{1/n}$ and make a change of variables
$\mathbf{u}\mapsto \mathbf{w}$, \,$u_j=\nu^j w_j$. Then the
operator~$\mathbf{H}_\rho$ is replaced by
$\widehat{\mathbf{H}}_\rho\colon X\to X^*$, where
$$
(\widehat{\mathbf{H}}_\rho\mathbf{w})_j=A_j w_j-\nu^{-1}B_{j-1}w_{j-1}-
\nu B_{j}^*w_{j+1},\qquad
w_n=w_0,\quad
w_{n+1}=w_1.
$$
Hence
$$
\widehat{\mathbf{H}}_\rho=\mathbf{A}+\nu^{-1}\mathbf{B}+\nu\mathbf{B}^*.
$$
Similarly,  $\mathbf{B}_\rho$ is replaced by
   $\widehat{\mathbf{B}}_\rho=\nu^{-1}\mathbf{B}$.
By the invariance of a determinant,
$$
G(\rho)=\det\bigl(\widehat{\mathbf{B}}_\rho^{-1}\widehat{\mathbf{H}}_\rho\bigr)
=\det(\nu\mathbf{B}^{-1}\mathbf{A}+
I+\nu^2\mathbf{B}^{-1}\mathbf{B}^*)
$$
is a polynomial of order $2mn$ in~$\nu$. Since   $\mathbf{B}_\rho$ and
$\mathbf{H}_\rho$ are linear in  $\rho$ and $\rho^{-1}$,    $G(\rho)$ is
a polynomial in
  $\rho$, $\rho^{-1}$. Thus,          $G(\rho)$ is a polynomial of order~$2m$
in~$\rho$. The senior coefficient is $\det(\mathbf{B}^{-1}\mathbf{B}^*)=1$.

We have $\mathbf{H}_\rho\mathbf{u}=0$ if and only if  $\mathbf{u}$
satisfies the variational system and $u_n=\rho u_0$, $u_{n+1}=\rho u_1$, or
 $Pw=\rho w$, where $w=(u_0,u_1)$. Hence $F(\rho)=\det (P-\rho I)=0$
is equivalent to $\det \mathbf{H}_\rho=0$.
Thus the polynomials $F(\rho)$ and $G(\rho)$ have the
 same roots and so they coincide.\qed

\subsection{Some applications}
\label{ssec2.5}
Identity  \eqref{eq2.8}
implies  that dynamical non-degeneracy of a periodic trajectory is
equivalent to the geometric non-degeneracy:
$\det \mathbf{H}=0\,\Leftrightarrow\,\det(P-I)=0$. Actually, the proof of
Theorem  \ref{th2.2} was based on this fact.

Equation  \eqref{eq2.8} gives
$$
\sigma(\mathbf{x})(-1)^m\det\mathbf{H}\det(P-I)>0.
$$

\begin{corollary}
\label{cor2.1}
Suppose that     $\sigma(\mathbf{x})(-1)^{m
+\operatorname{ind}\mathbf{H}}<0$. Then  $\mathbf{x}$ is dynamically
unstable: there is a real multiplier                $\rho>1$.
\end{corollary}

For example, the hypothesis holds when        $\sigma(\mathbf{x}) (-1)^m<0$
and
   $\mathbf{x}$ is a nondegenerate local minimum of the action     $\mathscr{A}$.

\begin{corollary}
\label{cor2.2}
Suppose that     $\sigma(\mathbf{x})(-1)^{m+n}<0$ and~$\mathbf{x}$ is a
nondegenerate local maximum of the action  $\mathscr{A}$. Then  $\mathbf{x}$
has a real multiplier             $\rho>1$.
\end{corollary}

Indeed, it is sufficient to use the following

\begin{proposition}
\label{pr2.2}
If   $\det(P-I)<0$, there is a real positive multiplier  $\rho>1$.
\end{proposition}

\noindent{\it Proof.}
Consider the characteristic polynomial
$F(\rho) = \det (P - \rho I)$.
Its roots are the multipliers of the periodic solution $\mathbf{x}$. We
have  $F(+\infty) = +\infty$ and $F(1) = \det (P - I) < 0$. Then
there exists a real root $\rho > 1$.
\qed

\begin{corollary}
\label{cor2.3}
If   $\sigma(\mathbf{x})(-1)^{\operatorname{ind}\mathbf{H}_{-1}}<0$,
there is a real positive multiplier       $\rho<-1$.
\end{corollary}

Indeed, for        $\rho=-1$ Hill's formula~\eqref{eq2.18} gives
$$
\sigma(\mathbf{x})\det\mathbf{H}_{-1}\det(I+P)>0.
$$

Let   $\mathbf{x}^2$ be the iterate of a periodic trajectory $\mathbf{x}$,
that is, the corresponding $2n$-periodic trajectory.

\begin{corollary}
\label{cor2.4}
Suppose $\sigma(\mathbf{x})(-1)^{\operatorname{ind}\mathbf{H}(\mathbf{x}^2)-
\operatorname{ind}\mathbf{H}(\mathbf{x})}<0$. Then  $\mathbf{x}$
is exponentially unstable. 
\end{corollary}

\noindent{\it Proof.}
Since $2n$-periodic vector fields along $\mathbf{x}^2$
are split into $n$-periodic and $n$-antiperiodic ones,
   $\operatorname{ind}\mathbf{H}(\mathbf{x}^2)=\operatorname{ind}\mathbf{H}(\mathbf{x})+
\operatorname{ind}\mathbf{H}_{-1}(\mathbf{x})$. It remains to use Corollary
\ref{cor2.3}.
\qed

\medskip

In the case $m=1$ there is a possibility to identify hyperbolicity or ellipticity of a periodic trajectory in terms of the index.

\begin{corollary}
\label{cor2.5}
Suppose that $m=1$. Then a nondegenerate periodic trajectory $\mathbf{x}$
 is hyperbolic if and only if
$\operatorname{ind}\mathbf{x}^2$ if even and elliptic if and only
if    $\operatorname{ind}\mathbf{x}^2$ is odd.
\end{corollary}

\noindent{\it Proof.} The
hyperbolicity of~$\mathbf{x}$ is equivalent to the hyperbolicity of
  $\mathbf{x}^2$.
For $m=1$,  $\mathbf{x}^2$ is hyperbolic if and only if it has
a multiplier                   $\rho>1$. This is equivalent to the inequality
$\sigma(\mathbf{x}^2)(-1)^{1+\operatorname{ind}\mathbf{H}(\mathbf{x}^2)}<0$.
It remains to note that $\sigma(\mathbf{x}^2)=\sigma^2(\mathbf{x})>0$.
\qed \medskip

Consider the convex billiard bounded by a hypersurface  $M$ in
  $\mathbb{R}^{n+1}$. Then the corresponding action is length and, by
Proposition \ref{pr2.1}, $\det B(x,y)>0$. Therefore,
$\sigma=(-1)^n$ and we obtain

\begin{corollary}
\label{cor2.6}
Suppose          $(-1)^{m+n+\operatorname{ind}\mathbf{H}(\mathbf{x})}<0$.
Then  $\mathbf{x}$ is exponentially unstable by Corollary
\ref{cor2.1}.
\end{corollary}

In particular,  $\mathbf{x}$ is exponentially unstable in
each of the following two cases 
\begin{itemize}
\item[$\bullet$] if  $m$ is odd and~$\mathbf{x}$ is a
nondegenerate local maximum of the billiard length functional;
\item[$\bullet$] if   $m+n$ is odd and $\mathbf{x}$ is a
nondegenerate local minimum of the billiard length functional.
\end{itemize}

For $m=1$ by the Birkhoff theorem \cite{21} (see also \cite{6}), any convex
billiard system has (at least) two periodic trajectories of period $n$ with rotation
number $k<n$, where one of them is a maximum of length, and hence generically
hyperbolic. The other has index 1, and so $\det(P-I)>0$. This implies that
the trajectory has no real multipliers $>1$. Indeed, if such a
multiplier exists, then the other one is also real     and greater than 1.
This contradicts $\det P=1$.

The problem of stability for billiard trajectories of period 2 is
systematically studied in the recent paper \cite{17}.  In this case,
the characteristic polynomial, as a function of $\rho+\rho^{-1}$, can be
presented as a determinant of some $m\times m$ matrix. This matrix is
explicitly determined by the matrices of second fundamental forms of the
surface $M$ at the end points of the trajectory.

The requirement for $n$ to be even in the hyperbolicity condition for a
periodic trajectory of minimal length ($m=1$) at the first glance looks
somewhat strange because the billiard trajectory minimizing $\mathscr{A}$
is naturally associated with a locally shortest closed geodesic on a
two-dimensional Riemannian manifold. Such geodesics due to
Poincar\'e \cite{2} are known to be hyperbolic. However one should
keep in mind that this Poincar\'e's result is valid only for
orientable geodesics (see details in \S\,\ref{sec5}) while a
periodic billiard trajectory with an odd period should be associated with a
non-orientable geodesic.

\begin{figure}[ht]
\begin{center}
\includegraphics{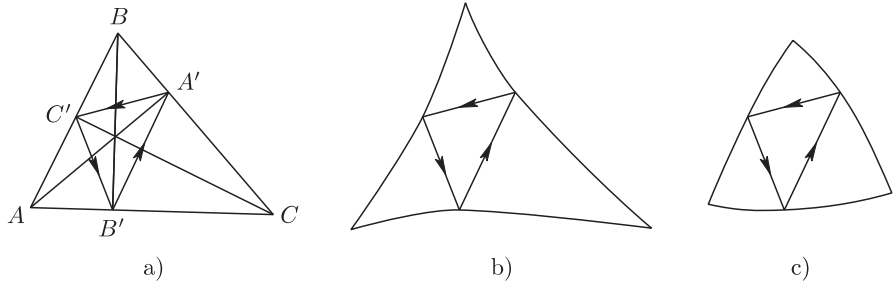}
\caption{Deformation of a parabolic periodic orbit}
\label{fig1}
\end{center}
\end{figure}

A simple example of an elliptic action minimizing billiard
trajectory with odd period can be constructed as follows. Let the
billiard curve be an acute-angled triangle $ABC$. Then by a
well-known theorem from planimetry the  projections
$A'$, $B'$, and $C'$ of the vertices to the opposite sides form a triangle \
(the orthotriangle)
which presents a local nondegenerate
minimum of the billiard action (Fig. \ref{fig1},~a). The
corresponding periodic trajectory is parabolic: its multipliers
  $\rho_{1,2}$ are equal to~$-1$.

A small deformation of the billiard curve does not destroy the
periodic trajectory $A'B'C'$ and just slightly deforms it. If the
boundary curve becomes concave, we obtain a Sinai billiard
\cite{22}. In this case the trajectory is hyperbolic (Fig.
\ref{fig1},~b)). If the boundary curve becomes strictly convex
(the curvature gets positive (Fig. \ref{fig1},~c)), then the
trajectory becomes elliptic still having a locally minimal action
provided the deformation is small.

\section{Continuous symmetry in a DLS}
\label{sec3}

\subsection{Discrete symmetry}
\label{ssec3.1} A
diffeomorphism $\psi\colon M\to M$ is a discrete {\it symmetry\/} of the
Lagrangian  $L$ if the map
$\widetilde\psi=\psi\times \psi$:
$M^2\to\nobreak M^2$,
\,$\widetilde\psi(x,y)=\bigl(\psi(x),\psi(y)\bigr)$ preserves  $L$:
$$
L\bigl(\psi(x),\psi(y)\bigr)=L(x,y).
$$
A more general definition is that  $\psi$ preserves~$L$
up to a cocycle:
$$
L\bigl(\psi(x),\psi(y)\bigr)=L(x,y)+f(x)-f(y).
$$
If there exists a function~$g$ such that  $g\circ \psi-g=f$, then  $\psi$
preserves the
gauge-equivalent Lagrangian          $\widehat L(x,y)=L(x,y)+g(x)-g(y)$.

\begin{proposition}
\label{pr3.1}
A symmetry takes a trajectory of a DLS into a trajectory. Thus,
$\widetilde\psi\circ T=T\circ \widetilde\psi$.
\end{proposition}

\noindent{\it Proof.}
Since   $\psi$ preserves the action functional~$\mathscr{A}$,
it takes critical points to critical points.
\qed \medskip

If a DLS $(M,L)$   admits a discrete symmetry group $\Gamma$,  then, in
principle,  symmetry can be removed by a factorization $\widetilde 
M=M/\Gamma$ of the configuration space. However it is useful to
keep in mind the following two aspects.

1.\kern1ex  Since in general $L(x,y) \ne L(g_1x,g_2y)$ for $g_1\ne g_2$, the
Lagrangian becomes multivalued after the factorization $M/\Gamma$.
This phenomenon is effectively used in the construction of a symbolic
dynamics by the method of anti-integrable limit, see \cite{23} and  \cite{24}
(a more general setup is discussed in \cite{25}, a continuous
analogue is presented in \cite{26}).

2.\kern1ex
 A periodic trajectory of the original system can turn into a trajectory
of the factorized system with a smaller period. Therefore the trajectory
can lose orientability. Moreover, the configuration space itself can lose
orientability. This happens, for example, in the case of a billiard system
inside a convex hypersurface $M\subset\mathbb{R}^{m+1}$ symmetric with respect to
the origin. Then $G = \{\mathrm{id},S\}$, where $S(x)=-x$.
Then $\widehat M=M/G$ is homeomorphic to the $m$-dimensional projective
space which is non-orientable for $m$ even.

\subsection{Noether symmetry}
\label{ssec3.2}
Let   $\mathbf{w}$ be a smooth vector field on  the configuration space
$M$ and~$\psi_s\colon M\to M$
its phase flow. We say that    $\mathbf{w}$ is a {\it symmetry field\/}
for the DLS if $\psi_s$ is a symmetry for~$L$ for all  $s$.

Define the vector fields $\mathbf{w}_1=(\mathbf{w},0)$ and
$\mathbf{w}_2=(0,\mathbf{w})$ on~$M^2$. Let
$\widetilde{\mathbf{w}}=\mathbf{w}_1+\mathbf{w}_2=(\mathbf{w},\mathbf{w})$
be the vector field corresponding to the group action $\widetilde
\psi_s=\psi_s\times \psi_s$. We have an equivalent version of the definition:
$\mathbf{w}$ is a symmetry  field for~$L$ if and only if
\begin{equation}
\label{eq3.1} D_{\widetilde{\mathbf{w}}} L(x,y)=D_{\mathbf{w}}
f(x)-D_{\mathbf{w}} f(y)
\end{equation}
for some function  $f$ on~$M$. If we replace~$L$
by its proper calibration, equation  \eqref{eq3.1} can be replaced by
\begin{equation}
\label{eq3.2}
D_{\widetilde{\mathbf{w}}}L=D_{\mathbf{w}_1}L+D_{\mathbf{w}_2}L=0.
\end{equation}

\begin{proposition}
\label{pr3.2}
Let   $\mathbf{w}$ be a symmetry field for $L$. Then
$$
\mathscr{J}=D_{\mathbf{w}_1}L=-D_{\mathbf{w}_2}L
$$
is a first integral of the corresponding DLS, that is,
        $\mathscr{J}\circ T=\mathscr{J}$.
\end{proposition}

\noindent{\it Proof.}
Suppose that     $(y,z)=T(x,y)$. Then
\begin{align*}
\mathscr{J}(x,y)-\mathscr{J}(y,z)&=D_{\mathbf{w}_2} L(y,z)-
D_{\mathbf{w}_2} L(x,y)
\\
&
=D_{\mathbf{w}_2} L(y,z)+D_{\mathbf{w}_1} L(y,z)
=D_{\widetilde{\mathbf{w}}} L(y,z)=0.
\end{align*}
The last expression vanishes by   \eqref{eq3.2}.
\qed \medskip

We call     $\mathscr{J}$ the Noether integral.

\begin{proposition}
\label{pr3.3} Let   $\mathbf{w}$ be a symmetry field. Then the group
action $\widetilde\psi_s\colon M^2\to M^2$ preserves~$\mathscr{J}$.
Equivalently, $\widetilde{\mathbf{w}}$ is tangent to the level surfaces
$$
N_c =\mathscr{J}^{-1}(c)\subset M^2.
$$
\end{proposition}

\noindent{\it Proof.} The
derivative of  $\mathscr{J}$ along  $\widetilde{\mathbf{w}}$ is
$$
D_{\widetilde{\mathbf{w}}}\mathscr{J}=(D_{\mathbf{w}_1}+
D_{\mathbf{w}_2})D_{\mathbf{w}_1} L=0
$$
by  \eqref{eq3.2} because the differential operators  $D_{\mathbf{w}_1}$
and  $D_{\mathbf{w}_1}+D_{\mathbf{w}_2}$ commute.
\qed \medskip

Note that $\mathscr{J}$ is the Hamiltonian generating the group $\widetilde
\psi_s$ of
symplectic transformations with respect to the symplectic form
\eqref{eq2.4} on~$M^2$.

\subsection{Routh reduction of order}
\label{ssec3.3}
Suppose system       $(M,L)$ admits
commuting independent symmetry fields     $\mathbf{w}^1,\dots,\mathbf{w}^k$:
$$
[\mathbf{w}^\alpha,\mathbf{w}^\beta]=0, \qquad
\alpha,\beta=1,\dots,k.
$$
Then the flows $\psi_{s}^\alpha\colon M\to M$ of symmetry fields commute.
Let  $G$ be the corresponding commutative group acting on        $M$ by
$$
x\mapsto\psi_s(x)=\psi_{s_1}^1\circ\dots\circ\psi_{s_k}^k(x),\qquad
s\in\mathbb{R}^k.
$$
 In general the flows
$\psi_{s_\alpha}^\alpha$ may be incomplete, and then $G$ is a local group acting on $M$.
 Since we are interested in a neighbourhood of a periodic orbit, these
non-local questions are irrelevant for us.
 Suppose that $\widetilde M=M/G$ is a smooth manifold and $\pi\colon
M\to\widetilde M$ a smooth fibration (at least locally this is
 always true).

Let   $\mathscr{J}^\alpha$ be the Noether integral of~$(M,L)$
corresponding to $\mathbf{w}^\alpha$ and let
  $\mathscr{J}=(\mathscr{J}^1,\dots,\mathscr{J}^k)$~
be the corresponding vector integral. We fix the value   $c\in\mathbb{R}^k$
and restrict $T$ to the level set    $N_c=\mathscr{J}^{-1}(c)\subset M^2$.
By Proposition  \ref{pr3.3}, the group~$G$ acts on  $N_c$.
If   $\widetilde M$ is a smooth manifold, then  $\widetilde N=N_c/G$
is also a smooth manifold, and $T$ defines a map
$\widetilde T\colon\widetilde N\to \widetilde N$. It is symplectic with respect
to the quotient symplectic structure $\widetilde\omega$
on~$\widetilde N$. We would like to represent~$\widetilde T$
as a discrete Lagrangian system with
quotient configuration space  $\widetilde M$ and Lagrangian $\widetilde L$
on $\widetilde M\times \widetilde M$. For this reduction we need

\medskip\noindent
{\bf Non-degeneracy assumption.} The
matrix  $G=(g^{\alpha\beta})$,
\begin{equation}
\label{eq3.3}
g^{\alpha\beta}(x,y)=
\langle B(x,y)\mathbf{w}^\alpha(x),\mathbf{w}^\beta(y)\rangle,\qquad
x,y\in M,
\end{equation}
is nondegenerate.
\medskip

First let                $c=0$. Let   $f(x,y)\in\mathbb{R}^k$
be the critical point of the function
$s\mapsto L(x,\psi_s(y))$, provided it exists and is unique.  Note that
the Hessian of this function equals $G(x,\psi_s(y))$, and so is nondegenerate.
The reduced Lagrangian is defined by
\begin{equation}
\label{eq3.4}
\widetilde L(x,y)=L\bigl(x,\psi_{f(x,y)}(y)\bigr).
\end{equation}
Since   $\widetilde L(x,y)=\widetilde L(\tilde x,\tilde y)$
depends only on~$\tilde x=\pi(x)$ and~$\tilde y=\pi(y)$,
it is a function on        $\widetilde M\times \widetilde M$.

Suppose now that  $c=(c^1,\dots,c^k)\ne 0$.  Locally there exist smooth
functions $\phi_\beta$,
$\beta=1,\dots, k$, on $M$ such that $D_{\mathbf{w}^\alpha}
\phi_\beta=\delta_{\beta}^\alpha$.
In general there are topological obstructions  to the existence of single
valued  globally defined $\phi_\alpha$. However, if
$\pi\colon M\to\widetilde M$ has fibre $\mathbb{R}^k$, then $\phi_\alpha$
exist globally.
Since we work in a neighbourhood of a periodic orbit, this is irrelevant for us.

Replace $L$ by gauge-equivalent Lagrangian\,\footnote{We use
 Einstein's sum rule with respect to repeated Greek indices,
but not Latin indices.}
$$
\widehat L(x,y)=L(x,y)+c^\alpha\bigl(\phi_\alpha(y)-\phi_\alpha(x)\bigr).
$$
Then  $\mathscr{J}^\alpha$ is replaced by
$$
\widehat{\!\!\mathscr J}^\alpha=\mathscr{J}^\alpha- c^\beta
D_{\mathbf{w}^\alpha}\phi_\beta=\mathscr{J}^\alpha-c^\alpha.
$$
Now    $\hat c=0$ and so  $\widetilde L$ can be defined by
          \eqref{eq3.4} with~$L$ replaced by $\widehat L$.

Here is Routh' Theorem for discrete Lagrangian systems.

\begin{proposition}
\label{pr3.4}
The projection
$\pi\colon M\to\widetilde M$ takes trajectories of the system $(M,L)$ with
  $\mathscr{J}=c$ to trajectories of the reduced Lagrangian system
$(\widetilde M,\widetilde L\,)$.
\end{proposition}

Next we give a coordinate version of the  Routh reduction.
Since the result is local, it is sufficient to perform the reduction
near a given trajectory $\mathbf{x}^0=(x_i^0)$.  Since the vector fields
$\mathbf{w}^\alpha$ are independent and commute, in a neighbourhood $U_i$ of the
point $x_i^0$ there are local coordinates $y_i\in\mathbb{R}^{m-k}$,
$z_i\in\mathbb{R}^k$ such that $\mathbf{w}^\alpha|_{U_i} =\partial / \partial z_i^\alpha$.
Similarly to the continuous case coordinates $z_i^\alpha$ are called
cyclic. The variables $y_i\in\mathbb{R}^{m-k}$ are local coordinates on
$\widetilde M$. Equation \eqref{eq3.2} means that
$$
L(x_i,x_{i+1})=\mathscr{L}(y_i,y_{i+1},u_i),\qquad
u_i=z_{i+1}-z_i.
$$
By \eqref{eq3.3}, the matrix
\begin{equation}
\label{eq3.5}
G_i=-\biggl(\frac{\partial^2 L}
{\partial z_i^\alpha\,\partial z_{i+1}^\beta}\biggr)=
\biggl(\frac{\partial^2 \mathscr{L}}
{\partial u_i^\alpha\,\partial u_i^\beta}\biggr)=\bigl(g_i^{\alpha\beta}\bigr)
\end{equation}
is nondegenerate.

Without loss of generality we assume that $c=0$. Then
\begin{equation}
\label{eq3.6}
\mathscr{J}=\partial_{u_i}\mathscr{L}(y_i,y_{i+1},u_i)=0.
\end{equation}
Equation  \eqref{eq3.6} can be locally solved with respect to
$u_i=f_i(y_i,y_{i+1})$. Then the Routh function is defined by
\begin{equation}
\label{eq3.7}
\widetilde L(y_i,y_{i+1})=\mathscr{L}\bigl(y_i,y_{i+1},f_i(y_i,y_{i+1})\bigr).
\end{equation}
Hence
\begin{equation}
\label{eq3.8}
\partial_{y_i}\widetilde L(y_i,y_{i+1})=
\partial_{y_i}\mathscr{L}(y_i,y_{i+1},u_i)\big|_{u_i=f(y_i,y_{i+1})},
\end{equation}
and similarly for the derivative with respect to $y_{i+1}$.

Suppose          $x_i=(y_i,z_i)$ is a trajectory of the system $(M,L)$
with~$\mathscr{J}=0$. Then
$$
\partial_{y_i}\bigl(\mathscr{L}(y_{i-1},y_i,u_{i-1})+
\mathscr{L}(y_i,y_{i+1},u_i)\bigr)=0.
$$
By  \eqref{eq3.8},
$$
\partial_{y_i}\bigl(\widetilde L(y_{i-1},y_i)+
\widetilde L(y_i,y_{i+1})\bigr)=0,
$$
so      $(y_i)$ is a trajectory of the reduced system
$(\widetilde M,\widetilde L\,)$.

To finish the proof of Routh's Theorem, it remains to show that
$\widetilde L$ is a discrete Lagrangian, that is,   it satisfies the twist
condition. This follows from

\begin{lemma}
\label{lem3.1}
Let   $\widetilde B_i=\widetilde B(y_i,y_{i+1})=
-\partial_{12}\widetilde L(y_i,y_{i+1})$ and~$B_i=B(x_i,x_{i+1})$. Then
\begin{equation}
\label{eq3.9}
\det \widetilde B_i=\frac{\det B_i}{\det G_i} \ne 0.
\end{equation}
\end{lemma}

\noindent{\it Proof.}
Putting   $(g_{\alpha\beta i})=G_i^{-1}$ we have
\begin{equation}
\label{eq3.10}
\frac{\partial}{\partial y_i}f_\alpha(y_i,y_{i+1})=-g_{\alpha\beta i}\,
\frac{\partial^2}{\partial y_i\,\partial u_i^\beta}
\mathscr{L}(y_i,y_{i+1},u_i)\big|_{u_i=f_i(y_i,y_{i+1})}.
\end{equation}
We differentiate  \eqref{eq3.8} with respect to  $y_{i+1}$
using  \eqref{eq3.10}:
$$
-\widetilde B_i=\frac{\partial^2\mathscr{L}}{\partial y_i\,\partial y_{i+1}}
+\frac{\partial f_\alpha}{\partial y_{i+1}}\,
\frac{\partial^2 \mathscr{L}}{\partial u_i^\alpha\,\partial y_i}
=\frac{\partial^2 L}{\partial y_i\,\partial y_{i+1}}+g_{\alpha\beta i}\,
\frac{\partial^2 L}{\partial y_{i+1}\,\partial z_{i+1}^\alpha}\,
\frac{\partial^2 L}{\partial y_i\,\partial z_i^\beta}\,,
$$
or in the matrix notation
$\widetilde B_i^{jk}=B_i^{jk}-B_i^{\alpha j}g_{\alpha\beta i}B_i^{k\beta}$.
It remains to make an exercise in linear algebra.
\qed \medskip

\subsection{Symplectic reduction for the Poincar\'e map}
\label{ssec3.4}
Suppose that the periodic trajectory is degenerate. Then the linear Poincar\'e
map $P\colon  W\to W$ has
a unit eigenvalue: there exists $w\ne 0$ such that $Pw=w$.
Since $P$ is symplectic, $\omega(w,Pu)=\omega(w,u)$, and so
$J_w(u)=\omega(w,u)$ is a linear first integral of $P$.
Then it is possible to reduce $P$ to a linear symplectic map
$\widetilde P\colon \widetilde W\to \widetilde W$ of lower dimension.

This section deals with  symplectic linear algebra, and
the origin of the symplectic map $P$ is  irrelevant.
In particular, the notations below will be used both for discrete and continuous
Lagrangian systems.

Suppose there are several eigenvectors corresponding to unit eigenvalue.
Let $V\subset\{w\in W: Pw=w\}$. Then $P$ has a first integral $J\colon
W\to V^*$: for $w\in V$, $\langle J(u),w\rangle=J_w(u)$.
We assume that $V$ is isotropic: $\omega\big|_V=0$.
Then $V\subset J^{-1}(0)$. We put $\widetilde W=J^{-1}(0)/V$.

\begin{proposition}[{\rm Poincar\'e}]
\label{pr3.5}
$P$ generates a reduced symplectic operator
$\widetilde P\colon \widetilde W\to \widetilde W$ such that the diagram
$$
\xymatrix
@C=30pt
{V \ar@{->}[d]\ar@{->}[r]^P& V\ar@{->}[d]
\\
\widetilde W \ar@{->}[r]^{\widetilde P}&\widetilde W}
$$
is commutative. Furthermore,
$$
\det(P-\rho I_W)=(1-\rho)^{2k}\det(\widetilde P-\rho I_{\widetilde W}),\qquad
k=\dim V.
$$
\end{proposition}

\subsection{Routh reduction for linear discrete Lagrangain systems}
\label{ssec3.5}
Next we translate Proposition \ref{pr3.5}     to the language of
the variational system, that is, the linear Lagrangian system $(E,\Lambda)$.
To any eigenvector $w$ of the Poincar\'e map $P$
there corresponds a non-zero $n$-periodic solution
$\mathbf{w} = (w_i)$ of the variational system. To the periodic solution
$\mathbf{w}$ there corresponds a
linear periodic first integral
$$
I_j(u_{j},u_{j+1})=\langle B_{j}w_{j},u_{j+1}\rangle-
\langle B_{j}u_{j},w_{j+1}\rangle.
$$
Indeed, if         $\mathbf{u}=(u_j)$ is a solution of   \eqref{eq2.14}, then
\begin{align*}
0&=\langle A_j w_j-B_j^* w_{j+1}-B_{j-1} w_{j-1}, u_j\rangle
-\langle A_j u_j - B_j^* u_{j+1} - B_{j-1} u_{j-1}, w_j\rangle
\\
&=I_{j}(u_j,u_{j+1})-I_{j-1}(u_{j-1},u_{j}).
\end{align*}
Hence
\begin{equation}
\label{eq3.11}
I_{j-1}(u_{j-1},u_{j})=I_j(u_j,u_{j+1}).
\end{equation}
In fact,               $I_j(u_j,u_{j+1})=J_w(u)$, where
    $J_w$
is the  integral of the Poincar\'e map
and $u\in W$ corresponds to the trajectory $(u_j)$.

Suppose now there are several eigenvectors and let $V\subset\operatorname{Ker}
(P-I)$ be an isotropic subspace. Denote by $\Gamma\subset X$ the set of periodic
trajectories corresponding to $V$. Let $ w^\alpha\in V$, $\alpha=1,\dots,k$,
be a basis in $V$. Then the corresponding periodic trajectories
$\mathbf{w}^\alpha=(w_i^\alpha)$ form a basis in $\gamma$. Let
$$
I_j^\alpha(u_j,u_{j+1})=\langle B_{j}w_{j}^\alpha,u_{j+1}\rangle-
\langle B_{j}u_{j},w_{j+1}^\alpha\rangle
$$
be the corresponding integrals of the variational system.
Since $V$ is isotropic, the integrals commute:
\begin{equation}
\label{eq3.12}
I_j^\alpha\bigl(w_j^\beta,w_{j+1}^\beta\bigr)=
\langle B_{j}w_{j}^\alpha,w_{j+1}^\beta\rangle
-\langle B_{j}w_{j}^\beta,w_{j+1}^\alpha\rangle
=\omega(w^\alpha,w^\beta)=0.
\end{equation}
We sometimes write $I_j=(I_j^1,\dots,I_j^k)$.

Below we need several non-degeneracy conditions.

\begin{condition}
\label{conA}
The symmetric matrix
$$
G_i=\bigl(g_i^{\alpha\beta}\bigr),\qquad
g_i^{\alpha\beta}=\langle B_iw_i^\alpha,w_{i+1}^\beta\rangle,
$$ is nondegenerate for
all  $i$. Denote
$(g_{\alpha\beta i})=\bigl(g_i^{\alpha\beta}\bigr)^{-1}=G_i^{-1}$.
\end{condition}

\begin{condition}
\label{conB}
The matrix
\begin{equation}
\label{eq3.13}
\overline{G}=\sum_{i=1}^n G_i^{-1}=(\bar g_{\alpha\beta}),\qquad
\bar g_{\alpha\beta}=\sum_{i=1}^n g_{\alpha\beta i},
\end{equation}
is nondegenerate.
\end{condition}

 Many of our results hold without condition~\ref{conB}, so we impose it later.
 Condition~\ref{conA} is used  almost everywhere, so we impose it now. In the
case of CLS, an analogue of condition  \ref{conA}
 is also introduced, but finally it turns out to be unessential, see
\S\,\ref{ssecA.3}.
An analogue of condition \ref{conB} is always satisfied for CLS.

  Denote
$$
F_i=\{w_i:\mathbf{w}\in\Gamma\}=\operatorname{span}(w_i^1,\dots,w_i^k)\subset
E_i.
$$
The reduced Poincar\'e map $\widetilde P$
corresponds to the reduced linear Lagrangian system $(\widetilde E,\widetilde
 \Lambda)$
with $\widetilde E_i=E_i/F_i$ which is obtained by the Routh reduction of
the system $(E,\Lambda)$. Under the non-degeneracy condition  \ref{conA},
$\dim F_i=k$  and the reduced configuration space $\widetilde E_i=E_i/F_i$ can
 be identified with
$$
E_i^\perp=\{u\in E_i:\langle B_{i-1}w_{i-1}^\alpha,u\rangle=0, \
\alpha=1,\dots,k\}
$$
via the projection $\Pi_i\colon E_i\to E_i^\perp$:
\begin{equation}
\label{eq3.14}
\Pi_i u=u_i-g_{\alpha\beta i}\langle B_{i-1}w_{i-1}^\alpha,u\rangle w_i^\beta.
\end{equation}
We represent any vector $u_i\in E_i$ as
\begin{equation}
\label{eq3.15}
u_i=v_i+\lambda_{\beta i} w_i^\beta,\qquad
v_i=\Pi_i u_i\in E_i^\perp, \quad
\lambda_{\beta i}=g_{\alpha\beta\,i-1}\langle B_{i-1}w_{i-1}^\alpha,u\rangle.
\end{equation}
The Routh reduction for DLS is described by the following

\begin{theorem}
\label{th3.1}
Let $\mathbf{u}=(u_i)$,\, $u_i\in E_i$, be a trajectory of
the system $(E,\Lambda)$ such that $I_i(u_i,u_{i+1})=0$. Then $\mathbf{v}=
(v_i)$,\, $v_i=\Pi_i u_i\in E_i^\perp$, is a  trajectory of the
 linear Lagrangian system
$(E^\perp,\Lambda^\perp)$ with the Lagrangian
$$
\Lambda_i^\perp(v_i,v_{i+1})=\frac{1}{2}\,\langle A_i v_i,v_i\rangle
-\langle B_iv_i,v_{i+1}\rangle-\frac{1}{2}\,\langle C_i v_i,v_i\rangle,
$$
where
$$
\langle C_iv_i,v_i\rangle=
g_{\alpha\beta i}\langle B_iv_i,w_{i+1}^\alpha\rangle
\langle B_iv_i,w_{i+1}^\beta \rangle.
$$
Conversely, if $\mathbf{v}$ is a trajectory of
the system $(E^\perp,\Lambda^\perp)$,
then there exists a trajectory $\mathbf{u}$ of the system $(E,\Lambda)$,
defined $\mod\;\Gamma$ such that
$I_i(u_i,u_{i+1})=0$ and $\Pi\mathbf{u}=\mathbf{v}$.
\end{theorem}

For the proof we will need an evident

\begin{lemma}
\label{lem3.2}
Let  $\mathbf{u}$,~$\mathbf{v}$ be such that $u_i-v_i\in F_i$:
$$
u_i=v_i+\lambda_{\beta i} w_i^\beta.
$$
Then            $I_i^\alpha(u_i,u_{i+1})=c^\alpha$ for all
         $\alpha=1,\dots,k$ and all $i$ if and only if
\begin{equation}
\label{eq3.16}
\Delta\lambda_{\alpha i}=\lambda_{\alpha\,i+1}-\lambda_{\alpha i}
=g_{\alpha\beta i}\bigl(c^\beta-I_i^\beta(v_i,v_{i+1})\bigr).
\end{equation}
\end{lemma}

Equation  \eqref{eq3.16} follows from
$$
c^\alpha=I_i^\alpha(u_i,u_{i+1})=I_i^\alpha(v_i,v_{i+1})
+g_i^{\alpha\beta}\Delta\lambda_{\beta i}.
$$

\noindent{\it Proof of Theorem~\ref{th3.1}.}
Let   $\mathbf{u}=(u_i)$, \,$u_i\in E_i$,
be a trajectory of
$(E,\Lambda)$ such that $I_{i}^\alpha(u_i,u_{i+1})=0$.   Then
for any variation $\phi_i\in E_i$ such that $\phi_i=0$ except for
$i=1,\dots, n$, we have
$$
0=h(\mathbf{u},\phi)=\sum_{i=1}^n\langle A_iu_i-B_{i}^* u_{i+1}
-B_{i-1} u_{i-1},\phi_i\rangle.
$$
Let   $\mathbf{v}=\Pi\mathbf{u}$. By  \eqref{eq3.15} and~\eqref{eq3.16},
$$
I_i^\alpha(v_i,v_{i+1})=-\langle B_iv_i,w_{i+1}^\alpha\rangle,\quad
\Delta\lambda_{\alpha i}=
g_{\alpha\beta i} \langle B_iv_i,w_{i+1}^\beta\rangle.
$$

Choose the variation~$\phi$ such that $\phi_i\in E_i^\perp$. Using
\begin{equation}
\label{eq3.17}
A_iw_i^\alpha=B_{i-1}w_{i-1}^\alpha+B_i^*w_{i+1}^\alpha
\end{equation}                                                     and
  $\langle B_{i-1}w_{i-1}^\beta,\phi_{i}\rangle=0$, we obtain
\begin{align*}
h(\mathbf{u},\phi)&=\sum_{i=1}^n \langle A_iv_i-B_{i}^* v_{i+1}
-B_{i-1} v_{i-1},\phi_i\rangle-\sum_{i=1}^n \Delta\lambda_{\beta i}
\langle B_i\phi_i,w_{i+1}^\beta\rangle
\\
&=h(\mathbf{v},\phi)-\sum_{i=1}^n g_{\alpha\beta i}
\langle B_iv_i,w_{i+1}^\alpha\rangle
\langle B_i\phi_i,w_{i+1}^\beta\rangle=h^\perp(\mathbf{v},\phi),
\end{align*}
where
\begin{equation}
\label{eq3.18}
h^\perp(\mathbf{v},\mathbf{v})=h(\mathbf{v},\mathbf{v})-\sum_{i=1}^n
g_{\alpha\beta i}\langle B_iv_i,w_{i+1}^\alpha\rangle
\langle B_iv_i,w_{i+1}^\beta\rangle.
\end{equation}
This is the  quadratic action functional for the
 system $(E^\perp,\Lambda^\perp)$.
Since $\phi_i\in E_i^\perp$, $i=1,\dots,n$, are arbitrary, $\mathbf{v}$ is a
trajectory of $(E^\perp,\Lambda^\perp)$.
We skip the proof of the converse.
\qed \medskip

For $n$-periodic $\mathbf{v}$ the
bilinear form  $h^\perp$ on~$X^\perp=E_1^\perp\times\dots\times E_n^\perp$
equals
\begin{equation}
\label{eq3.19}
h^\perp(\mathbf{v},\mathbf{v})=\sum_{i=1}^n \langle A_i^\perp v_i-
{B_i^\perp}^* v_{i+1}-B_{i-1}^\perp v_{i-1},v_i \rangle,\qquad
v_{n+1}=v_1,\quad
v_0=v_n,\
\end{equation}
where
$$
A_i^\perp =R_i(A_i-C_i),\qquad
B_i^\perp=R_{i+1}B_i
$$       and
  $R_i\colon E_i^*\to (E_{i}^\perp)^*$ is the restriction map. The
reduced variational system is
$$
A_i^\perp v_i-{B_{i}^\perp}^* v_{i+1}-B_{i-1}^\perp v_{i-1}=0.
$$

Let   $\mathbf{H}^\perp,\mathbf{B}^\perp\colon X^\perp\to{X^\perp}^*$ be
the corresponding linear operators. We also put
$$
\sigma^\perp\beta^\perp=\biggl(\,\prod\det B_i^\perp\biggr)^{-1},\qquad
\beta^\perp>0, \quad
\sigma^\perp \in \{1,-1\}.
$$
Then analogously to  \eqref{eq2.12},
$$
\det\mathbf{B}^\perp=(-1)^{m-k} \prod \det B_i^\perp.
$$
Hill's theorem \ref{th2.1}, applied to the reduced system, gives

\begin{corollary}
\label{cor3.1} The following reduced Hill formula holds:
$$
\det(\widetilde P-I)=\frac{\det\mathbf{H}^\perp}{\det\mathbf{B}^\perp}
=\sigma^\perp (-1)^{m-k} \beta^\perp \det\mathbf{H}^\perp.
$$
\end{corollary}

To use this formula for stability problems, we need to know $\sigma^\perp$
and the Morse index of $h^\perp$. However, the relation between
the Morse indices of $h$ and $h^\perp$ is not evident. The reason is that a
periodic sequence $\mathbf{v}\in X^\perp$ in general corresponds to a
non-periodic sequence $\mathbf{u}\notin X$ such that $I_i(u_i,u_{i+1})=0$.
We discuss this problem  in the next two sections.

Lemma \ref{lem3.1}  implies $\operatorname{ind}  B_i =
\operatorname{ind}  B_i^\perp + \operatorname{ind}  G_i \pmod 2$.
Therefore,
\begin{equation}
\label{eq3.20}
\sigma^\perp=(-1)^{\sum\operatorname{ind} B_i^\perp}=\sigma (-1)^{\sum\operatorname{ind} G_i}.
\end{equation}

\subsection
({Degeneracy of  \$h\$})
{Degeneracy of  $h$}
\label{ssec3.6} We denote by
      $\Gamma\subset X$ the space of periodic solutions corresponding to
$V\subset\operatorname{Ker}(P-I)$. It is spanned by
$\mathbf{w}_1,\dots,\mathbf{w}_k\in X$.
Since   $\mathbf{H}\mathbf{w}=0$ for $\mathbf{w}\in\Gamma$,
the Hessian bilinear form
$h(\mathbf{u},\mathbf{u})=\langle\mathbf{H} \mathbf{u},\mathbf{u}\rangle$
is degenerate 
and defines a bilinear form  $\hat h$
on $\widehat X=X/\Gamma$. To compare  $h$ with  $h^\perp$
we need to restrict  $h$ to the level set of  $I$.
Let
\begin{align}
\label{eq3.21}
Y&=\{\mathbf{u}\in X: I_1(u_1,u_{2})=\dots=I_n(u_n,u_1)\},
\\
\label{eq3.22}
Z&=\{\mathbf{v}\in X: v_i\in F_i\}=\{\mathbf{v}\in X: v_i=
\lambda_{\alpha i} w_i^\alpha\}.
\end{align}

\begin{proposition}
\label{pr3.6}
The spaces  $Y$ and  $Z$ are $h$-orthogonal, that is,
  $h(\mathbf{u},\mathbf{v})
   =   
0$ for all $\mathbf{u}\in Y$ and
$\mathbf{v}\in Z$. Moreover,  $Y$ is the
$h$-orthogonal complement to    $Z$:
$$
Y=\{\mathbf{u}\in X:h(\mathbf{u},\mathbf{v})=0\ \text{for all} \
\mathbf{v}\in Z\}.
$$ The
restriction of $h$ to $Z$ is
$$
h(\mathbf{v},\mathbf{v})=\sum_{i=1}^n g_i^{\alpha\beta}
\Delta\lambda_{\alpha i}\Delta\lambda_{\beta i}
=\sum_{i=1}^n\langle G_i\Delta\lambda_i,\Delta\lambda_i\rangle,
\qquad
v_i=\lambda_{\alpha i} w_i^\alpha.
$$
\end{proposition}

\noindent{\it Proof.}
Let   $\mathbf{u}\in Y$ and
$\mathbf{v}=(\lambda_{\alpha i} \mathbf{w}_i^\alpha) \in Z$. Then
\begin{align*}
h(\mathbf{u},\mathbf{v})&=\sum_{i=1}^n\langle A_iu_i-B_i^*u_{i+1}
-B_{i-1}u_{i-1},\lambda_{\alpha i} w_i^\alpha\rangle
\\
&=\sum_{i=1}^n\lambda_{\alpha i}\bigl(\langle A_iu_i,w_i^\alpha\rangle
-\langle B_iw_{i}^\alpha,u_{i+1}\rangle
-\langle B_{i-1}u_{i-1},w_i^\alpha\rangle\bigr)
\\
&=\sum_{i=1}^n\lambda_{\alpha i}\bigl(\langle A_iw_i^\alpha,u_i\rangle
-\langle B_iu_{i},w_{i+1}^\alpha\rangle
-\langle B_{i-1}w_{i-1}^\alpha,u_i\rangle\bigr)
\\
&=\sum_{i=1}^n\langle A_iw_i^\alpha-B_{i}^*w_{i+1}^\alpha
-B_{i-1}w_{i-1}^\alpha,\lambda_{\alpha i}u_i\rangle=0
\end{align*}
by       \eqref{eq3.17}; we used that
$$
\langle B_iw_{i}^\alpha,u_{i+1}\rangle
=\langle B_iu_{i},w_{i+1}^\alpha\rangle+c^\alpha,\qquad
\langle B_{i-1}u_{i-1},w_i^\alpha\rangle
=\langle B_{i-1}w_{i-1}^\alpha,u_i\rangle-c^\alpha.
$$

Conversely, if  $\mathbf{v}=(\lambda_{\alpha i} w_i^\alpha)$ and
$h(\mathbf{u},\mathbf{v})=0$ for all  $\lambda_\alpha$, then
$$
\langle A_iu_i-B_i^*u_{i+1}-B_{i-1}u_{i-1},w_i^\alpha\rangle=0.
$$
Using  \eqref{eq3.17} we obtain
$$
I_i^\alpha(u_i,u_{i+1})=I_{i-1}^\alpha(u_{i-1},u_i).
$$
Thus,          $\mathbf{u}\in Y$.

Next we compute the restriction of $h$ to~$Z$. Let
$\mathbf{v}=(\lambda_{\alpha i}\mathbf{w}_i^\alpha)$.
Then by      \eqref{eq3.17},
\begin{align*}
h(\mathbf{v},\mathbf{v})&=\sum_{i=1}^n
\bigl(\lambda_{\alpha i}\lambda_{\beta i}
\langle A_iw_i^\alpha,w_i^\beta\rangle-\lambda_{\alpha i}\lambda_{\beta\,i-1}
\langle B_{i-1}w_{i-1}^\alpha,w_i^\beta\rangle
\\
&\qquad
-\lambda_{\alpha i}\lambda_{\beta\,i+1}\langle B_{i}w_{i}^\alpha,
w_{i+1}^\beta\rangle\bigr)
\\
&=\sum_{i=1}^n\bigl(\lambda_{\alpha i}(\lambda_{\beta i}-\lambda_{\beta\,i+1})
\langle B_{i}w_{i}^\alpha,w_{i+1}^\beta\rangle
\\
&\qquad+\lambda_{\beta i}(\lambda_{\alpha i}-\lambda_{\beta\,i-1})
\langle B_{i-1}w_{i-1}^\beta,w_{i}^\alpha\rangle\bigr)
\\
&=\sum_{i=1}^n\bigl(g_i^{\alpha\beta} \lambda_{\alpha i}
(\lambda_{\beta i}-\lambda_{\beta i+1})+
g_{i-1}^{\alpha\beta}\lambda_{\alpha i}(\lambda_{\beta i}-
\lambda_{\beta i-1})\bigr)
\\
&=\sum_{i=1}^n g_i^{\alpha\beta}\Delta\lambda_{\alpha i}
\Delta\lambda_{\beta i}.
\end{align*}
\qed \medskip

We obtain a quadratic form on  $Z$:
$$
h\big|_Z(\mathbf{v},\mathbf{v})
=\langle\mathbf{H}_Z\mathbf{v},\mathbf{v}\rangle
=\sum_{i=1}^n\langle G_i\Delta\lambda_{i},\Delta\lambda_{i}\rangle
=\langle\mathbf{G}\lambda,\lambda\rangle,
$$
where $\lambda_i\in\mathbb{R}^k$, \,$\mathbf{H}_Z\colon Z\to Z^*$, and
the operator $\mathbf{G}\colon \mathbb{R}^{kn}\to \mathbb{R}^{kn}$
is defined by
$$
(\mathbf{G}\lambda)_i=G_{i-1}\Delta\lambda_{i-1}-G_i\Delta\lambda_{i}.
$$
\goodbreak
\noindent
We have
$$
\operatorname{Ker}\mathbf{G}=\{\lambda\in\mathbb{R}^{kn}: G_1\Delta\lambda_{1}=
\dots=G_n\Delta\lambda_{n}\}.
$$
Thus,          $\Delta\lambda_{i}=G_i^{-1}c$.
For $\lambda\in\operatorname{Ker} \mathbf{G}$  
the equation
$\vphantom{\Big|}\smash[b]{\displaystyle\sum_{i=1}^n} \Delta\lambda_{i}=0$
implies       $\overline{G}c=0$, where $\overline{G}$ is the
matrix  \eqref{eq3.13}.

Now we impose the non-degeneracy assumption~\ref{conB}: the matrix
$\overline G$ is nondegenerate. Then
$$
\operatorname{Ker}\mathbf{G}=\{\lambda\in\mathbb{R}^{nk}:\lambda_{1}=\dots=\lambda_{n}\}
\quad\text{and}\quad
\operatorname{Ker}\mathbf{H}_Z=\Gamma.
$$
Let     $\widehat Z=Z/\Gamma$.

\begin{proposition}
\label{pr3.7} The
form  $h\big|_{\widehat Z}$ is nondegenerate and
\begin{equation}
\label{eq3.23}
\operatorname{ind} h\big|_{\widehat Z}=\sum_{i=1}^n \operatorname{ind} G_i-\operatorname{ind}\overline G.
\end{equation}
\end{proposition}

\noindent{\it Proof.}
By Proposition       \ref{pr3.6}, in the coordinates
$\mu_i=\Delta\lambda_i\in\mathbb{R}^k$
we have
$h\big|_{\widehat Z}=\chi\big|_\Theta$, where $\chi$ is the following
quadratic form on  $\mathbb{R}^{nk}$:
$$
\chi(\mu,\mu)=\sum_{i=1}^n \langle G_i \mu_i,\mu_i \rangle,
$$
and the space  $\Theta\subset\mathbb{R}^{nk}$ is defined by the condition
$$
\Theta=\biggl\{\mu\in\mathbb{R}^{nk}: \sum_{i=1}^n \mu_i=0\biggr\}.
$$
Below we use the same notation $\chi$ for the corresponding bilinear form.
Consider the $k$-dimensional space
$$
\Xi=\{\mu\in\mathbb{R}^{nk}: \mu_i=G_i^{-1}\nu, \
\nu\in\mathbb{R}^k, \ j=1,\dots,n\}.
$$
Since  $\overline G$ is nondegenerate, $\mathbb{R}^{nk}=\Xi\oplus\Theta$
and moreover,  $\chi(\mu,\xi)=0$ for any  $\mu\in\Theta$ and~$\xi\in \Xi$.
Therefore, the spaces  $\Xi$ and  $\Theta$ are $\chi$-orthogonal and
$$
\sum_{i=1}^n \operatorname{ind} G_i=\operatorname{ind}\chi=\operatorname{ind}\chi\big|_\Xi+\operatorname{ind}\chi\big|_\Theta
=\operatorname{ind}\overline G+\operatorname{ind} h\big|_{\widehat Z}.
$$
\qed \medskip

\begin{proposition}
\label{pr3.8}
$Y+Z=X$ and  $Y\cap Z=\Gamma$.
\end{proposition}

This follows from

\begin{lemma}
\label{lem3.3}
For any    $\mathbf{v}\in X$ there exists $\mathbf{u}=\Phi\mathbf{v}\in Y$,
unique $\mod \Gamma$, such that  $\mathbf{u}-\mathbf{v}\in Z$.
Explicitly,     $u_i=v_i+\lambda_{\beta i} w_i^\beta$, where
    the $\lambda_{\beta_i}$ satisfy  \eqref{eq3.16} and
\begin{equation}
\label{eq3.24}
c^\alpha=\kappa^{\alpha\beta}\sum_{i=1}^n g_{\beta\delta i}
I_i^\delta(v_i,v_{i+1}),\qquad
(\kappa^{\alpha\beta})=(\bar g_{\alpha\beta})^{-1}={\overline G}{}^{-1}.
\end{equation}
\end{lemma}

The map     $\Phi\colon X\to \widehat Y=Y/\Gamma$ satisfies
$\Phi(Z)=0$ and~$\Phi\big|_Y$ is the identity $\mod\Gamma$.
Proposition \ref{pr3.8} follows immediately.

To prove Lemma       \ref{lem3.3}, for given   $\mathbf{v}\in X$ we
find   $\mathbf{u}$ such that
$$
I_i^\alpha(u_i,u_{i+1})=c^\alpha=\mathrm{const}.
$$
By Lemma  \ref{lem3.2}, the $\lambda_{\alpha i}$
satisfy       \eqref{eq3.16}. If   $\lambda_{\alpha i}$ is
$n$-periodic in  $i$, then
   $\displaystyle\sum_{i=1}^n \Delta\lambda_{\alpha i}=0$, which gives
  \eqref{eq3.24}.
If  \eqref{eq3.24} holds, then equation \eqref{eq3.16}
determines $n$-periodic     $\lambda_{\alpha i}$ modulo a constant independent
of  $i$. Thus, $\mathbf{u}$ is defined uniquely modulo
$\Gamma$.

Let     $\widehat X=X/\Gamma$ and let
$\pi_\Gamma\colon X\to\widehat X$ be the corresponding canonical projection.
Then there exists a linear map  $\widehat\Pi\colon \widehat X\to X^\perp$
such that the
following diagram is commutative:
$$
\xymatrix{
X \ar@{->}[dr]_\Pi\ar@{->}[rr]^{\pi_\Gamma}&& \widehat X\ar@{->}[dl]^{\widehat\Pi}
\\
& X^\perp
}
$$

\begin{corollary}
\label{cor3.2}
The maps    $\widehat\Pi\big|_{\widehat Y}\colon \widehat Y\to X^\perp$ and
$\Phi\big|_{X^\perp}\colon X^\perp\to\widehat Y$ are
mutually inverse isomorphisms.
\end{corollary}

\noindent{\it Proof.}
If   $\Pi\mathbf{v}=0$ for $\mathbf{v}\in Y$, then  $\mathbf{v}\in Z$ and
hence $\mathbf{v}\in Y\cap Z=\Gamma$ by
Proposition~\ref{pr3.7}. The equation $Y+Z=X$ implies that
$\Pi\big|_Y$ is surjective. We also have $\Phi(Z)=0$.
\qed \medskip

The spaces     $\widehat Y,\widehat Z\subset \widehat X$ are orthogonal
with respect to the bilinear form  $\hat h$ on~$\widehat X$.
By Proposition       \ref{pr3.7}, $\hat h$ is nondegenerate on~$\widehat Z$
and its index is given by     \eqref{eq3.23}.

From the point of view of the Routh reduction it is natural to consider the
space
$$
Y^0=\{\mathbf{u}\in X:I_i(u_i,u_{i+1})=0, \ i=1,\dots,n\}\subset Y.
$$
Define    $\mathbf{d}_\alpha=(d_{\alpha i})\in X^*$ by
\begin{equation}
\label{eq3.25}
d_{\alpha i}=g_{\alpha\beta\,i-1}B_{i-1}w_{i-1}^\beta
-g_{\alpha\beta i}B_{i}^*w_{i+1}^\beta.
\end{equation}
Note that    $\langle d_{\alpha i},w_i^\alpha\rangle=0$. We have
$$
\langle\mathbf{d}_\alpha,\mathbf{u}\rangle
=\sum_{i=1}^n g_{\alpha\beta i}I_i^\beta(u_{i},u_{i+1}),\qquad
\mathbf{u}\in X.
$$

\begin{proposition}
\label{pr3.9}
\begin{equation}
\label{eq3.26}
Y^0+Z=\{\mathbf{v}\in X:\langle\mathbf{d}_\alpha,\mathbf{v}\rangle=0, \
\alpha=1,\dots,k\}.
\end{equation}
\end{proposition}

\noindent{\it Proof.}
Let   $\mathbf{v}\in X$. Then  $\mathbf{v}\in Y^0+Z$,
\,$\mathbf{u}=\Phi(\mathbf{v})\in Y^0$. By Lemma       \ref{lem3.3},
$$
0=\sum_{i=1}^ng_{\alpha\beta i}I_i^\beta(v_i,v_{i+1})
=\sum_{i=1}^n\langle d_{\alpha i},v_i\rangle
=\langle\mathbf{d}_\alpha,\mathbf{v}\rangle.
$$
\qed \medskip

Next we compute the restriction $h\big|_Y$.  Let $h^\perp$
be the  bilinear action form for the reduced system
$(E^\perp,\Lambda^\perp)$, and let $h^\top=h\circ\Phi\big|_{X^\perp}$, where
$\Phi$ is the operator from Corollary \ref{cor3.2},
be the bilinear form on $X^\perp$ corresponding to $h\big|_Y$,
see the commutative diagram
$$
\xymatrix
@R=12pt@C=40pt
{\widehat Y \ar@<-.8ex>
[dd]
_{\widehat\Pi|^{}_{\widehat Y}}
\ar@{->}[dr]^{\hat h|^{}_{\widehat Y}}&
\\
&\mathbb{R}
\\
\ar@{->}[ur]_{h^{\top}}\ar@<-.2ex>
[uu]_{\Phi}X^{\perp}}
$$

\begin{proposition}
\label{pr3.10}
For any    $\mathbf{v}\in X^\perp$,
\begin{equation}
\label{eq3.27}
h^\top(\mathbf{v},\mathbf{v})=h^\perp(\mathbf{v},\mathbf{v})+
\bar g_{\alpha\beta} c^\alpha c^\beta,
\end{equation}                     where
the coefficients~$c^\alpha(\mathbf{v})$ are defined by~\eqref{eq3.24}.
\end{proposition}

This follows from a more general formula which we prove next.

\begin{lemma}
\label{lem3.4}
Let   $\mathbf{v}\in X$ and~$\mathbf{u}=\Phi(\mathbf{v})$. Then
\begin{equation}
\label{eq3.28}
h(\mathbf{u},\mathbf{u})=h(\mathbf{v},\mathbf{v})-
\bar g_{\alpha\beta}c^\alpha c^\beta+\sum_{i=1}^n g_{\alpha\beta i}
I_i^\alpha(v_i,v_{i+1})I_i^\beta(v_i,v_{i+1}),
\end{equation} where
the coefficients~$c^\alpha(\mathbf{v})$ are defined by~\eqref{eq3.24}.
\end{lemma}

\noindent{\it Proof.}
By Proposition \ref{pr3.6},
$$
h(\mathbf{u},\mathbf{u})=h(\mathbf{v},\mathbf{v})+\sum_{i=1}^n
g_i^{\alpha\beta}\Delta\lambda_{\alpha i}\Delta\lambda_{\beta i}.
$$
By    \eqref{eq3.16} and  \eqref{eq3.24},
\begin{align*}
\sum_{i=1}^n g_i^{\alpha\beta}\Delta\lambda_{\alpha i}\Delta\lambda_{\beta i}
&=\sum_{i=1}^n g_{\alpha\beta i}\bigl(c^\alpha-I_i^\alpha(v_i,v_{i+1})\bigr)
\bigl(c^\beta-I_i^\beta(v_i,v_{i+1})\bigr)
\\
&=\sum_{i=1}^ng_{\alpha\beta i}I_i^\alpha(v_i,v_{i+1})I_i^\beta(v_i,v_{i+1})
-\bar g_{\alpha\beta}c^\alpha c^\beta.
\end{align*}
\qed \medskip

\noindent{\it Proof of Proposition~\ref{pr3.10}.}
We use \eqref{eq3.28} with~$\mathbf{v}\in X^\perp$. Then
$$
I_i^\alpha(v_i,v_{i+1})=-\langle B_iv_i,w_{i+1}^\alpha\rangle,
$$
and so
$$
\sum_{i=1}^n g_{\alpha\beta i}I_i^\alpha(v_i,v_{i+1})I_i^\beta(v_i,v_{i+1})
=\sum_{i=1}^n g_{\alpha\beta i}\langle B_iv_i,w_{i+1}^\alpha\rangle
\langle B_iv_i,w_{i+1}^\beta\rangle.
$$
Last we use \eqref{eq3.18}.
\qed \medskip

\subsection
{The indices of~$h$ and~$h^\perp$}
\label{ssec3.7}
In this subsection we discuss the relation between
$\operatorname{ind} h\big|_Y=\operatorname{ind} \hat h\big|_{\widehat Y}$ and~$\operatorname{ind} h^\perp$.
Using the isomorphisms $\Phi\colon X^\perp\to \widehat Y$ and
$\widehat\Pi\colon\widehat Y\to X^\perp$, we compare instead the indices of
         $h^\perp$ and~$h^\top=h\circ \Phi$ on~$X^\perp$.
As mentioned earlier,  $h^\top$ and~$h^\perp$ coincide on  $X_0^\perp=\Pi Y^0$.

We need some assumptions on the unit eigenspace of the linear Poincar\'e map
$P\colon W\to W$.
Suppose $V=\operatorname{Ker}(P-I)$ is isotropic. It is well known
(see, for example, \cite{27}) that the generalized eigenspace
 $N=\operatorname{Ker}(P-I)^{2m}$    is symplectic. Since $V\subset N$
is isotropic, $\dim N\ge 2k$.
We consider the least degenerate case $\dim N=2k$. Then
$$
N=\operatorname{Ker}(P-I)^2=\{v\in W:Pv-v\in V\}
$$
is symplectic and~$V=\operatorname{Ker}(P-I)$ is a Lagrangian subspace of
  $N$.
Consider the bilinear form
\begin{equation}
\label{eq3.29}
b(v,w)=\omega\bigl((P-I)v,w\bigr),\qquad
v,w\in N.
\end{equation}
A computation shows \cite{27} that $b$ is symmetric.

Let $\Sigma$ be the set of trajectories $\mathbf{v}$ of the variational
system corresponding to the vectors $v\in N$. Then $\Gamma\subset\Sigma$.
Trajectories in $\Sigma\setminus\Gamma$ are not periodic.
The projection
$$
\Omega^\perp=\Pi\Sigma\subset X^\perp
$$
consists of periodic sequences. We have a natural map
$\Psi=\Phi\Pi\colon\Sigma\to \widehat Y$. Set
$\widehat\Omega=\Psi\Sigma$.
We will see that
$\Omega^\perp$ is orthogonal to $X_0^\perp=\Pi Y^0$ with respect
to   $h^\perp$, and $\hat\Omega$ is orthogonal to $\widehat Y$
with respect to $h$.
Since $h^\top$ and $h^\perp$ coincide on $X_0^\perp$, the difference of their indices
is determined by
their restrictions to the complement of $X_0^\perp$. Thus if
$\Omega^\perp \oplus X_0^\perp=X^\perp$, the difference of the indices
is determined by $h^\top\big|_{\Omega^\perp}$ and $h^\perp\big|_{\Omega^\perp}$.
However, in general $X^\perp\ne \Omega^\perp \oplus X_0^\perp$.
To ensure this expansion           we need
an extra non-degeneracy  condition \ref{conC} below.

Take a basis $w^1,\dots,w^k$ in $V$ and the
conjugate basis $q_1,\dots,q_k$ in a Lagrangian complement of $V$ in $N$. Then
$w^1,\dots,w^k,q_1,\dots,q_k$ is a symplectic basis in $N$ and
\begin{equation}
\label{eq3.30}
\omega(w^\alpha,w^\beta)=0,\quad
\omega(q_\alpha,q_\beta)=0,\quad
\omega(w^\alpha,q_\beta)=\delta^\alpha_\beta,\quad
Pq_\alpha=q_\alpha+s_{\alpha\beta} w^\beta.
\end{equation}

Combining  \eqref{eq3.29} and~\eqref{eq3.30}, we obtain
$$
b(q_\alpha,q_\beta)=s_{\alpha\beta}.
$$
Hence the matrix $S=(s_{\alpha\beta})$ is symmetric. Define
symmetric matrices               $A=(a_{\alpha\beta})$ and
$A^\perp=(a^\perp_{\alpha\beta})$:
\begin{equation}
\label{eq3.31}
a_{\alpha\beta}=s_{\alpha\delta}\kappa^{\delta\varepsilon}s_{\varepsilon\beta}
-s_{\alpha\beta}, \qquad
a^\perp_{\alpha\beta}=s_{\alpha\beta}-\bar g_{\alpha\beta}.
\end{equation}
Below we need another non-degeneracy assumption.

\begin{condition}
\label{conC}                                                     The
matrix  $A^\perp=(a^\perp_{\alpha\beta})$ is nondegenerate.
\end{condition}

\begin{definition}
We say that a periodic trajectory is nondegenerate      $\mod V$
if the non-degeneracy conditions  \ref{conA},~\ref{conB}, and~\ref{conC}
hold.
\end{definition}

\begin{theorem}
\label{th3.2}
Suppose that     $\mathbf{x}$ is nondegenerate $\mod V$. Then
\begin{equation}
\label{eq3.32}
\operatorname{ind} h-\operatorname{ind} h\big|_Z-\operatorname{ind} h^\perp=\operatorname{ind} A-\operatorname{ind} A^\perp.
\end{equation}
\end{theorem}

We prove Theorem \ref{th3.2} in \S\,\ref{ssec3.8}.

\begin{corollary}
\label{cor3.3}
Suppose that  $\mathbf{x}$ is nondegenerate $\mod V$. Then
$$
\operatorname{ind} h=\operatorname{ind} h^\perp+\operatorname{ind}
h\big|_Z+\operatorname{ind} b+\operatorname{ind}\overline{G}\, \pmod 2.
$$
\end{corollary}

Indeed, $SA^\perp=A \overline G$ implies
$$
\operatorname{sign}\det A=\operatorname{sign} (\det S\det G\det A^\perp).
$$
Therefore,     $\operatorname{ind} A-\operatorname{ind}
A^\perp=\operatorname{ind} b+\operatorname{ind} \overline{G}\pmod 2$.

Proposition    \ref{pr3.7} implies 

\begin{corollary}
\label{cor3.4}
$\operatorname{ind} h=\operatorname{ind}
h^\perp+\displaystyle\sum_{i=1}^n\operatorname{ind} G_i+\operatorname{ind} b\, \pmod2$.
\end{corollary}

Equation  \eqref{eq3.20} combined with Corollary~\ref{cor3.4} imply

\begin{corollary}
\label{cor3.5}
$\sigma (-1)^{\operatorname{ind} h}=\sigma^\perp(-1)^{\operatorname{ind}
h^\perp+\operatorname{ind} b}$.
\end{corollary}

\subsection
{The spaces  $\Omega$ and  $\Omega^\perp$}
\label{ssec3.8}
Consider solutions
$$
\mathbf{w}^\alpha=(w_i^\alpha), \qquad
\mathbf{q}_\alpha=(q_{\alpha i})
$$
of the variational system corresponding to the symplectic basis
$w^\alpha$, $q_\alpha$  of the space $N$. They form a basis in $\Sigma$.
The solutions $\mathbf{w}^\alpha\in \Gamma$ are $n$-periodic and satisfy
\eqref{eq3.12}.
Equation \eqref{eq3.30} implies that $\mathbf{q}_\alpha$ satisfy
\begin{equation}
\label{eq3.33}
q_{\alpha,i+n}-q_{\alpha i}=s_{\alpha\beta} w_i^\beta.
\end{equation}
Since the basis is symplectic,
\begin{gather}
\label{eq3.34}
I_i^\alpha(q_{\beta i},q_{\beta\,i+1})=
\langle B_i w_i^\alpha,q_{\beta\,i+1}\rangle-
\langle B_i q_{\beta i},w_{i+1}^\alpha\rangle
=\omega(w^\alpha,q_\beta)=\delta_\beta^\alpha,
\\
\label{eq3.35}
\langle B_j q_{\alpha j},q_{\beta\,j+1}\rangle-
\langle B_j q_{\beta j},q_{\alpha\,j+1}\rangle
=\omega(q_\alpha,q_\beta)=0.
\end{gather}

Let     $\mathbf{q}_\alpha^\perp=\Pi\mathbf{q}_\alpha\in
\Omega^\perp\subset X^\perp$. Then
$q_{\alpha i}^\perp=q_{\alpha i}-\lambda_{\alpha\beta i}^\perp w_i^\beta$,
where
\begin{equation}
\label{eq3.36}
\lambda_{\alpha\beta\,i+n}^\perp-\lambda_{\alpha\beta i}^\perp=
s_{\alpha\beta},\qquad
\lambda_{\alpha\beta i}^\perp=g_{\gamma\beta\,i-1}\langle B_{i-1}
w_{i-1}^\gamma, q_{\alpha i}\rangle.
\end{equation}
For any    $\alpha=1,\dots,k$ define
$\widehat{\mathbf{q}}_\alpha\in \widehat Y$ by
$\widehat{\mathbf{q}}_\alpha=\Psi\mathbf{q}_\alpha$. Then
$$
\hat q_{\alpha i}=q_{\alpha i}-\nu_{\alpha\beta i} w_i^\beta,
$$ where the
coefficients  $\nu_{\alpha\beta i}$ are chosen so that the  $\hat q_\alpha$
are periodic and
\begin{equation}
\label{eq3.37}
I_j(\hat q_{\beta j},\hat q_{\beta\,j+1})=
-\langle B_j\hat q_{\beta j},w_{j+1}^\alpha\rangle
+\langle B_j w_j^\alpha,\hat q_{\beta\,j+1}\rangle=c_\beta^\alpha
\end{equation} are
independent of~$j$. Then
\begin{equation}
\label{eq3.38}
\nu_{\alpha\beta\,i+n}-\nu_{\alpha\beta i}=s_{\alpha\beta},
\qquad
\Delta\nu_{\alpha\beta i}=
s_{\alpha\gamma}\kappa^{\gamma\delta} g_{\delta\beta i},
\qquad
c_\beta^\alpha=\delta_\beta^\alpha-s_{\beta\gamma}\kappa^{\gamma\alpha}.
\end{equation}
Since   $\widehat{\mathbf{q}}_\alpha$ is defined modulo    $\Gamma$, we
have  $\widehat{\mathbf{q}}_\alpha\in\widehat
Y$ and  $\widehat{\mathbf{q}}_\alpha=\Phi\mathbf{q}_\alpha^\perp$.

We have the spaces
$$
\Omega^\perp=\Pi\Sigma=\operatorname{span}(\mathbf{q}^\perp_1,\dots,
\mathbf{q}^\perp_k)\subset X^\perp,
\quad\
\widehat\Omega=\Psi\Sigma=\Phi\Omega^\perp
=\operatorname{span}(\widehat{\mathbf{q}}_1,\dots,
\widehat{\mathbf{q}}_k)\subset \widehat Y.
$$
We define $\widehat Y^0$ as the image of  $Y^0$ under the projection
$\pi_\Gamma\colon X\to \widehat X=X/\Gamma$.

\begin{theorem}
\label{th3.3}
Suppose that the non-degenaracy conditions
                        \ref{conA},~\ref{conB}, and  \ref{conC} hold. Then

\begin{numa}

\item
\label{th3.3a}
$\widehat Y=\widehat\Omega\oplus\widehat Y^0$ and
$X^\perp=\Omega^\perp\oplus X_0^\perp$;

\item
\label{th3.3b}
the maps
$\widehat\Pi\big|_{\widehat\Omega}\colon \widehat\Omega\to\Omega^\perp$ and
$\widehat\Pi\big|_{\widehat Y^0}\colon \widehat Y^0\to X_0^\perp$
are linear isomorphisms;

\item
\label{th3.3c}
$h\big|_{Y^0}=h^\perp\big|_{X_0^\perp} \circ\Pi\big|_{Y^0}$,
\,$\hat h\big|_{\widehat Y^0}=
h^\perp\big|_{X_0^\perp}\circ\widehat\Pi\big|_{\widehat Y^0}$;

\item
\label{th3.3d}
for any    $\mathbf{u}\in Y^0$ and  $\alpha,\beta=1,\dots,k$,
$$
h(\widehat{\mathbf{q}}_\alpha,\mathbf{u})=0, \qquad
h(\widehat{\mathbf{q}}_\alpha,\widehat{\mathbf{q}}_\beta)=a_{\alpha\beta};
$$

\item
\label{th3.3e}
for any    $\mathbf{u}^\perp\in\Pi Y^0$ and  $\alpha,\beta=1,\dots,k$,
$$
h^\perp(\mathbf{q}_\alpha^\perp,\mathbf{u}^\perp)=0, \qquad
h^\perp(\mathbf{q}_\alpha^\perp,\mathbf{q}_\beta^\perp)=
a_{\alpha\beta}^\perp.
$$
\end{numa}
\end{theorem}

The proof of Theorem  \ref{th3.3} is contained in \S\,\ref{ssecA.1}.

Now we prove Theorem  \ref{th3.2}. Recall that by
Proposition          \ref{pr3.6}, $Y\perp_h Z$ and
$\widehat Y\perp_{\hat h}\widehat Z$. By Theorem \ref{th3.3} we have the
$\hat h$-orthogonal expansion
$\widehat X=\widehat Z\oplus\widehat\Omega\oplus\widehat Y^0$
and the  $h^\perp$-orthogonal expansion
  $X^\perp=\Omega^\perp\oplus X_0^\perp$ and
$\hat h\big|_{\widehat Y^0}=h^\perp\big|_{\widehat\Pi\widehat Y^0}
\circ\widehat\Pi\big|_{\widehat Y^0}$. Therefore,
$$
\operatorname{ind}\hat h-\operatorname{ind}\hat h\big|_{\widehat Z}-\operatorname{ind} h^\perp
=\operatorname{ind}\hat h\big|_{\widehat\Omega}-\operatorname{ind} h^\perp\big|_{\Omega^\perp}
=\operatorname{ind}(a_{\alpha\beta})-\operatorname{ind}(a^\perp_{\alpha\beta}).
$$
It remains to use the equations
$$
\operatorname{ind} h=\operatorname{ind}\hat h,\quad
\operatorname{ind} h\big|_Z=\operatorname{ind}\hat h\big|_{\widehat Z},\quad
\operatorname{ind} h\big|_\Omega=\operatorname{ind}\hat h\big|_{\widehat\Omega}.
$$

\subsection
{Degeneracy for
$\rho$-index form}
\label{ssec3.9}
The connection between the indices for the original and the reduced systems
   is much simpler for the 
$\rho$-index form.  We
take complex        $\rho\in S^1$ and perform the same computation
for the corresponding Hermitian form~$h_\rho$ on the complex space
$$
X_\rho=\{\mathbf{u}=(u_j)_{j\in \mathbb{Z}}:u_j\in E_j, \ u_{j+n}=\rho u_j\}
$$
of quasiperiodic sequences.  We define
$Y_\rho,Z_\rho\subset X_\rho$ by the same formulae \eqref{eq3.21}:
\begin{align}
\label{eq3.39}
Y_\rho&=\{\mathbf{u}\in X_\rho: I^1(u_1,u_2)=\dots=I^n(u_n,u_{n+1})\},
\\
\label{eq3.40}
Z_\rho&=\{\mathbf{v}\in X_\rho: v_i=\lambda_{\alpha i} w_j^\alpha, \
\lambda_{\alpha\,i+n}=\rho\lambda_{\alpha j}\}.
\end{align}

The main difference is that for          $\rho\ne 1$,
$\mathbf{w}^\alpha\notin X_\rho$. This implies, in particular, that
$$
Y_\rho=Y_\rho^0=\{\mathbf{u}\in X_\rho: I_j(u_j,u_{j+1})=0, \ j=1,\dots,n\},
$$

\begin{proposition}
\label{pr3.11}
The spaces   $Y_\rho,Z_\rho$ are      $h$-orthogonal, that is,
        $h(\mathbf{u},\overline{\mathbf{v}})=\nobreak0$ for all
$\mathbf{u}\in Y_\rho$ and  $\mathbf{v}\in Z_\rho$. If
$h(\mathbf{u},\overline{\mathbf{v}})=0$ for all  $\mathbf{v}\in
Z_\rho$, then $\mathbf{u}\in Y_\rho$. The restriction of  $h_\rho$ to~$Z_\rho$
is given by
\begin{gather*}
h(\mathbf{v},\overline{\mathbf{v}})=\sum_{i=1}^ng_i^{\alpha\beta}
\Delta\lambda_{\alpha i} \Delta\bar\lambda_{\beta i}
=\sum_{i=1}^n\langle G_i\Delta\lambda_i,\Delta\bar\lambda_i\rangle,
\\
\mathbf{v}=(\lambda_{\alpha i}\mathbf{w}_i^\alpha),
\qquad
\lambda_{\alpha\,n+i}=\rho\lambda_{\alpha i}.
\end{gather*}
\end{proposition}

The proof is the same as for    $\rho=1$ (see Proposition  \ref{pr3.6}). 

Suppose assumption  \ref{conA} holds.\footnote{We do not need
assumption  \ref{conB} in this section.} Then for  $\rho\ne 1$,
               $Y_\rho$ is the $h$-orthogonal complement to  $Z_\rho$
and     $h\big|_{Z_\rho}$ is nondegenerate. Indeed,
$$
h\big|_Z(\mathbf{v},\mathbf{v})=
\langle \mathbf{G}_\rho\lambda,\bar\lambda\rangle,
$$
where $\lambda_i\in\mathbb{C}^{k}$ and
$$
(\mathbf{G}_\rho\lambda)_i=G_{i-1} \Delta\lambda_{ i-1}
-G_i \Delta\lambda_{i},\qquad
\lambda_{i+n}=\rho\lambda_i.
$$

\begin{lemma}
\label{lem3.5}
\begin{equation}
\det\mathbf{G}_\rho=(-1)^k\rho^{-k}(\rho-1)^{2k}\prod_{i=1}^n\det G_i.
\end{equation}
\end{lemma}

This follows from Hill's formula \eqref{eq2.18} applied to~ the DLS
$(F,\Lambda_F)$ with the bilinear action form~$h\big|_Z$.  The corresponding
Poincar\'e map  $P_F$ has the matrix of the form
$$
P_F \sim \begin{pmatrix}
I & \overline{G}\,
\\
0 & I
\end{pmatrix}.
$$
Hence          $\det(P-\rho I)=(\rho-1)^{2k}$. The operators $B_i$
for the system
$(F,\Lambda_F)$ are equal to the~$G_i$.

\begin{proposition}
\label{pr3.12}
Let   $\rho\in S^1$, \,$\rho\ne 1$. Then  $X_\rho=Y_\rho\oplus Z_\rho$ and
  $Y_\rho\cap Z_\rho=\nobreak\{0\}$.
\end{proposition}

\noindent{\it Proof.}
Let   $\mathbf{u}\in X_\rho$. If we want to find  $(\lambda_{\alpha i})$,
\,$\lambda_{\alpha\,i+n}=\rho\lambda_{\alpha i}$, such that
\begin{equation}
\label{eq3.42}
(u_i+\lambda_{\alpha i} w_i^\alpha)\in Y_\rho,
\end{equation}
then  \eqref{eq3.26} gives
$$
\lambda_{\alpha\,n+1}-\lambda_{\alpha 1}
=\sum_{i=1}^n g_{\alpha\beta i} I_i^\beta(u_i,u_{i+1}).
$$
Thus
$$
\lambda_{\alpha 1}=(\rho-1)^{-1}\sum_{i=1}^n g_{\alpha\beta i}
I_i^\beta(u_i,u_{i+1}).
$$
Similarly we find  $\lambda_{\alpha 2},\dots,\lambda_{\alpha n}$.
\qed \medskip

Formula~\eqref{eq3.42} defines a projection $\Phi_\rho\colon X_\rho\to Y_\rho$,
\,$\Phi_\rho Z_\rho=0$, \,$\Phi_\rho\big|_{Y_\rho}=I$. We have
$$
\operatorname{ind} h_\rho=\operatorname{ind}
h_\rho\big|_{Z_\rho}+\operatorname{ind} h_\rho\big|_{Y_\rho}.
$$ The
projection $\Pi\colon X_\rho\to X_\rho^\perp$ gives an isomorphism of
$h_\rho\big|_{Y_\rho}$ and the  $\rho$-index form  $h_\rho^\perp$
for the reduced system     $(E^\perp,\Lambda^\perp)$. Lemma  \ref{lem3.5}
implies that the Hill  $\rho$-determinants for the original and the reduced
system are related by 
$$
\det H_\rho=\det H_\rho^\perp \det \mathbf{G}_\rho=
\det H_\rho^\perp(2-\rho-\rho^{-1})^k\prod_{i=1}^n\det G_i.
$$
Hence
$$
\operatorname{ind}
h_\rho\big|_{Z_\rho}=\sum_{i=1}^n\operatorname{ind} G_i\,\pmod 2.
$$

\begin{corollary}
\label{cor3.6}
If   $\rho\in S^1$, \,$\rho\ne 1$, then
$$
\operatorname{ind} h_\rho=\operatorname{ind}
h_\rho^\perp+\sum_{i=1}^n\operatorname{ind} G_i\,\pmod 2,
\qquad
\operatorname{null} h_\rho= \operatorname{null} h_\rho^\perp.
$$
\end{corollary}

\section{Reversible version}
\label{sec4}

\subsection{Reversible DLS}
\label{ssec4.1}
Let   $S\colon M\to M$ be a smooth involution: $S^2
  =
\mathrm{id}_M$. We say that a
DLS is         $S$-reversible if   $S$ is a time reversing symmetry for
    $L$:            for any  $x,y\in M$,
\begin{equation}
\label{eq4.1}
L(Sx,Sy)=L(y,x).
\end{equation}
Equivalently, the Lagrangian $L$
is invariant under the involution
$\widetilde S\colon M^2\to M^2$, \,$\widetilde S(x,y)=(Sy,Sx)$.

The simplest example is    $S=\mathrm{id}$, that is, $L(x,y)=L(y,x)$
(for example, a billiard system or a
standard map).               A non-trivial  $S$ appears in the
system  \eqref{eq2.5} if the potential is even.
Then  $S(x) = -x$. An analogous
possibility exists in billiards with some symmetry conditions.

\begin{proposition}
\label{pr4.1}
 Suppose that $T\colon M^2\to M^2$ is generated by an
$S$-reversible DLS and $T(x,y) = (y,z)$. Then
$$
T\bigl(S(z),S(y)\bigr)=\bigl(S(y),S(x)\bigr).
$$
Hence         $T$ is conjugate to $T^{-1}$:
\begin{equation}
\label{eq4.2}
T^{-1}\circ \widetilde S=\widetilde S\circ T.
\end{equation}
\end{proposition}

The proof follows by differentiating the identity
$$
L(x,y)+L(y,z)=L(Sz,Sy)+L(Sy,Sx).
$$

If   $\mathbf{x}=(x_i)$ is a periodic orbit of a DLS, then
   $\widetilde{\mathbf{x}}=(S x_{-i})$
is also a periodic orbit.
 A periodic orbit is called reversible if $\widetilde{\mathbf{x}}=\mathbf{x}$
modulo translations.
 The group $\mathbb{Z}_n=\mathbb{Z}/n\mathbb{Z}$ acts on  the set of $n$-periodic
sequences $(x_i)_{i\in\mathbb{Z}}$ in $M$ by translation
$(x_i)_{i\in\mathbb{Z}}\mapsto  (x_{i+k})_{i\in\mathbb{Z}}$, and we should
identify periodic orbits obtained in such a way. Any $n$-periodic
sequence is determined by $(x_1,\dots,x_n)\in M^n$ and the
translation group $\mathbb{Z}_n$ acts on $M^n$   by cyclic permutations.
Thus the set of periodic sequences is the quotient
$\mathscr{M}=M^n/\mathbb{Z}_n$. Define an involution $R\colon \mathscr{M}\to
\mathscr{M}$ by
$R(\mathbf{x})=(S x_{-i})$. Let $\mathscr{M}_+=\mathbf{x}\in
\mathscr{M}:R\mathbf{x}=\mathbf{x}\}$ be
the set of fixed points of $R$. Thus $\mathbf{x}\in \mathscr{M}_+$ if and
only if
$S(x_{j-i})=x_{i}$ for some $j\in\mathbb{Z}$ and all $i\pmod n$.

\begin{proposition}
\label{pr4.2}
$\mathbf{x}\in \mathscr{M}_+$ is a reversible periodic orbit if and only
if
$\mathbf{x}$ is a critical point of the functional
$\mathscr{A}_+
    =\nobreak 
\mathscr{A}\big|_{\mathscr{M}_+}$. 
\end{proposition}

This is a well known property of functions invariant under an involution.
Indeed,
let $X=T_\mathbf{x} \mathscr{M}$. Then
$J=dR(\mathbf{x})\colon X\to X$ is an involution. Denote
$E_\pm=\{\xi\in X: J\xi=\pm\xi\}$. Then  $X=E_+\oplus E_-$ and
$E_+=T_\mathbf{x} \mathscr{M}_+$. Since  $\mathscr{A}$ is
$R$-invariant, we have
$d\mathscr{A}(\mathbf{x})\circ J=d\mathscr{A}(\mathbf{x})$. Thus,
$d\mathscr{A}(\mathbf{x})\xi=0$ for all  $\xi\in E_-$.

For any critical point      $\mathbf{x}\in \mathscr{M}_+$ let
 $h=d^2\mathscr{A}(\mathbf{x})$ be the Hessian, that is, 
 the second differential of
  $\mathscr{A}$. This is a bilinear form on~$X$. Then
$h\big|_{E_+}=d^2\mathscr{A}_+(\mathbf{x})$ is the second differential
of      $\mathscr{A}_+$. For     $\xi=\xi_++\xi_-$, \,$\xi_\pm\in E_\pm$,
we obtain
$$
h(\xi,\xi)=h(\xi_+,\xi_+)+h(\xi_-,\xi_-).
$$
Indeed, since  $h$ is $J$-invariant,
$$
h(\xi_-,\xi_+)=h(J\xi_-,J\xi_+)=h(-\xi_-,\xi_+)=-h(\xi_-,\xi_+).
$$
If we represent  $h$ by a linear operator  $\mathbf{H}\colon X\to X^*$, then
   $J^*\mathbf{H} J=\mathbf{H}$, and so  $\mathbf{H}E_\pm=E_\pm^*$.

Let us introduce on $M$ an $S$-invariant Riemannian metric. It
defines an $R$-invariant   metric ${(\,\cdot\,{,}\,\cdot\,)}$ on $X$.  Then
$h(\xi,\eta)=( H\xi,\eta)$, where $H\colon X\to X$ is a symmetric
operator. The spaces $E_\pm$ are orthogonal with respect to the
metric and $HE_\pm\subset E_\pm$. Denote $H_\pm=H|_{E_\pm}$.
We obtain

\begin{proposition}
\label{pr4.3}
$H=H_+\oplus H_-$  and~$\det H=\det H_+\det H_-$. 
\end{proposition}

Reversible periodic trajectories   $\mathbf{x}\in \mathscr{M}_+$ are of
3 types     $\tau=0,1,2$ depending on the number of fixed points of
 $S$ they contain.
\begin{itemize}
\itemindent16pt
\item[{\it Type\/} 0:] $n=2k$ and~$\mathbf{x}=(x_1,\dots,x_k,Sx_k,\dots,Sx_1)$.

\item[{\it Type\/} 1:] $n=2k-1$ and
$\mathbf{x}=(x_1,\dots,x_k,Sx_k,\dots,Sx_2)$, where  $x_1=Sx_1$.

\item[{\it Type\/} 2:] $n=2k-2$ and
$\mathbf{x}=(x_1,\dots,x_k,Sx_{k-1},\dots,Sx_2)$, where $x_1=Sx_1$ and
$x_k=Sx_k$.
\end{itemize}

For all types
$\mathbf{x}=\mathbf{x}_\tau(\mathbf{y})\in\mathscr{M}_+$ is determined by
$\mathbf{y}=(x_1,\dots,x_k)\in M^k$. Thus the action
functional  $\mathscr{A}_+$ on~$\mathscr{M}_+$ gives a function on $M^k$:
$$
\mathscr{A}_\tau(\mathbf{y})=\mathscr{A}_+(\mathbf{x}_\tau(\mathbf{y}))
=\mathscr{A}_+(\mathbf{x}),\qquad
\tau=0,1,2.
$$
Denote
$$
\mathscr{B}(\mathbf{y})=\sum_{i=1}^{k-1} L(x_i,x_{i+1}).
$$

\begin{lemma}
\label{lem4.1}
$\mathbf{x}=\mathbf{x}_\tau(\mathbf{y})$ is a periodic orbit of type
$\tau=0,1,2$ if and only if $\mathbf{y}$ is a critical point of
\begin{alignat*}{2}
\mathscr{A}_0(\mathbf{y})&=2\mathscr{B}(\mathbf{y})+
L(Sx_1,x_1)+L(x_k,Sx_k),&\qquad
\mathbf{y}&\in M^k,
\\
\mathscr{A}_1(\mathbf{y})&=2\mathscr{B}(\mathbf{y})+L(x_k,Sx_k),&\qquad
\mathbf{y}&\in N\times M^{k-1},
\\
\mathscr{A}_2(\mathbf{y})&=2\mathscr{B}(\mathbf{y}),&\qquad
\mathbf{y}&\in N\times M^{k-2}\times N,
\end{alignat*}                         respectively.
\end{lemma}

The functional~$\mathscr{A}$ on~$\mathscr{M}$ admits a similar representation.
For example, consider
the case of periodic orbits of type  0. A
point $\mathbf{x}\in M^n$ can be written as
$$
\mathbf{x}=(y_1,\dots,y_k,Sz_k,\dots,Sz_1).
$$
Then, since      $L(Sz_{i+1},Sz_i)=L(z_i,z_{i+1})$,
$$
\mathscr{A}(\mathbf{x})=\mathscr{B}(\mathbf{y})+
\mathscr{B}(\mathbf{z})+L(y_k,Sz_k)+L(Sz_1,y_1).
$$
For $\mathbf{x}$ a periodic orbit of type 0 we have
$\mathbf{y}=\mathbf{z}=(x_1,\dots,x_k)$. We write
$\mathbf{u}\in T_\mathbf{x} M^n$ as
$$
\mathbf{u}=(v_1,\dots,v_k,J_kw_k,\dots,J_1w_1),\qquad
\mathbf{v},\mathbf{w}\in T_\mathbf{y} M^k,
$$
where $J_i=dS(x_i)\colon T_{x_i}M\to T_{Sx_i}M$. Taking the second differential
of      $\mathscr{A}$ we get
$$
h(\mathbf{u},\mathbf{u})=k(\mathbf{v},\mathbf{v})+
k(\mathbf{w},\mathbf{w})-\langle B_{k}^*J_kw_{k},v_k\rangle-
\langle v_1,B_0J_1w_1\rangle,
$$
where
$$
k(\mathbf{v},\mathbf{v})=\sum_{i=1}^k \langle A_iv_i-B_{i-1}v_{i-1}-
B_i^*v_{i+1},v_i\rangle,\qquad
v_0=0,\quad
v_{k+1}=0,
$$
and  $B_k=B(x_k,Sx_k)$, \,$B_0=B(Sx_1,x_1)$. Denote
$$
C_1=-B_0J_1\colon E_1\to E_1^*,\qquad
C_k=-B_{k}^*J_k\colon E_k\to E_k^*.
$$
Note that              $C_1=C_1^*$ and  $C_k=C_k^*$ are symmetric..
Thus
\begin{equation}
\label{eq4.3}
h(\mathbf{u},\mathbf{u})=k(\mathbf{v},\mathbf{v})+k(\mathbf{w},\mathbf{w})
+\langle C_kv_{k},w_k\rangle+\langle C_1v_1,w_1\rangle.
\end{equation}

Let us compute the corresponding bilinear forms $h\big|_{E_\pm}$.
For $\mathbf{u}\in E_\pm$ we have  $\mathbf{w}=\pm\mathbf{v}$, so
  $\mathbf{u}$ is determined by      $\mathbf{v}$:
$$
h_\pm(\mathbf{u},\mathbf{u})=h_\pm^0(\mathbf{v},\mathbf{v})=
2k(\mathbf{v},\mathbf{v})\pm\langle C_kv_{k},v_k\rangle \pm
\langle C_1v_1,v_1\rangle.
$$
Similarly, for any      $\tau=0,1,2$ a vector $\mathbf{u}\in E_\pm$
is determined by      $\mathbf{v}=(v_1,\dots,v_k)\in T_\mathbf{y} M$.

\begin{lemma}
\label{lem4.2}
For a reversible orbit of type  $\tau=0,1,2$,
$h_\pm(\mathbf{u},\mathbf{u})=h_\pm^\tau(\mathbf{v},\mathbf{v})$ has the
form
\begin{align*}
h_\pm^0(\mathbf{v},\mathbf{v})&=2k(\mathbf{v},\mathbf{v})
\pm\langle C_kv_{k},v_k\rangle\pm\langle C_1v_1,v_1\rangle,
\\
h_\pm^1(\mathbf{v},\mathbf{v})&=2k(\mathbf{v},\mathbf{v})
\pm\langle C_kv_{k},v_k\rangle,
\\
h_\pm^2(\mathbf{v},\mathbf{v})&=2k(\mathbf{v},\mathbf{v}).
\end{align*}
The domain of  $h_\pm^\tau$ is  $V_\pm^\tau$, where
\begin{align*}
V_\pm^0&=T_\mathbf{y} M^k,
\\
V_\pm^1&=\{\mathbf{v}\in T_\mathbf{y} M^k:J_1v_1=\pm v_1\},
\\
V_\pm^2&=\{\mathbf{v}\in T_\mathbf{y} M^k:J_1v_1=\pm v_1, \ J_kv_k=\pm v_k\}.
\end{align*}
\end{lemma}

Consider  a periodic orbit of type  0. Then the domains of
$h_+^0$ and $h_-^0$ coincide and
$$
h_+^0(\mathbf{v},\mathbf{v})-h_-^0(\mathbf{v},\mathbf{v})=
2\langle C_1v_1,v_1\rangle+2\langle C_kv_k,v_k\rangle.
$$

\begin{corollary}
\label{cor4.1}
Suppose that $\mathbf{x}\in \mathscr{M}_+$ is a periodic
orbit of type 0 which is a nondegenerate local minimum point of
$\mathscr{A}_+$. If the symmetric operators
$C_1$ and $C_k$ are non-positive, then $\mathbf{x}$ is a
nondegenerate local minimum for $\mathscr{A}$.  If $\mathbf{x}\in
\mathscr{M}_+$ is a
nondegenerate local maximum point of $\mathscr{A}_+$ and $C_1,C_k$ are
non-negative, then $\mathbf{x}$ is a nondegenerate local maximum for
$\mathscr{A}$. In both  cases $\det H=\det H_-\det H_+>0$.
\end{corollary}

For periodic orbits of type $\tau= 1$        or $  2$, the domains
of $h_+^\tau$  and $h_-^\tau$ are different. When $S=\mathrm{id}$, then
$V_-^\tau\subset
V_+^\tau$.

\begin{proposition}
\label{pr4.4}
Suppose that $S=\mathrm{id}_M$, and let $\mathbf{x}$ be a
reversible periodic trajectory of type $\tau$ which is a
nondegenerate minimum for $\mathscr{A}_\tau$.
Then in each of the three cases 
\begin{itemize}
\item[{\rm(a)}]      $\tau=2$,
\item[{\rm(b)}]      $\tau=1$ and  $C_k$ is non-positive,
\item[{\rm(c)}] $\tau=0$ and  $C_k$,~$C_1$ are non-positive
\end{itemize}
$\mathbf{x}$ is a nondegenerate minimum for $\mathscr{A}$.
\end{proposition}

If   $\mathbf{y}$ is a nondegenerate maximum for  $\mathscr{A}_\tau$, then
a similar statement  holds provided that
  $C_k$ and~$C_1$
are non-negative rather than non-positive.

\begin{proposition}
\label{pr4.5}
Let   $S=\mathrm{id}$.  If         $\mathscr A_+$~
is a minimal periodic orbit of type  $2$, then
      $h_\rho$ is positive definite for all    $\rho\in S^1$.
Hence              $\mathbf{x}$ is hyperbolic.  
\end{proposition}

\noindent{\it Proof.}
Take complex               $\mathbf{u}\in T_\mathbf{x}^{\mathbb{C}}M^{2k}$.
Then
$$
h_\rho(\mathbf{u},\overline{\mathbf{u}})=
k(\mathbf{u}^+,\overline{\mathbf{u}}{\kern.2ex}^+)
+k(\mathbf{u}^-,\overline{\mathbf{u}}{\kern.2ex}^-),
$$         where
    $\mathbf{u}^+=(u_1,\dots,u_k)$ and
$\mathbf{u}^-=(J_{2k}u_{2k},\dots,J_{k+1}u_{k+1})$ are complex vectors from
   $T_\mathbf{y}^{\mathbb{C}} M^k$. Since         $k$
is positive definite,         $h_\rho$ is positive definite.
\qed \medskip

\subsection{Some applications}
\label{ssec4.2}

\begin{corollary}
\label{cor4.2}
Reversible geometric degeneracy $\det H_+ = 0$ of a reversible
trajectory implies the dynamical degeneracy $\det
(P-I) = 0$.
\end{corollary}

Next we give some statements on dynamical stability (in fact,
instability) of reversible trajectories.

\begin{corollary}
\label{cor4.3}
Let   $\mathbf{x}$
be a reversible periodic trajectory
such that $\mathbf{y}$ is a nondegenerate minimum point of $\mathscr{A}_+$.
Suppose also that $\sigma(\mathbf{x})(-1)^m<0$. In the case
$\operatorname{type}(\mathbf{x})=1$ we also need the condition that $C_k$ is
non-positive, and in the case   $\operatorname{type}(\mathbf{x})=0$ that $C_1$
and $C_k$ are non-positive.
  Then $\mathbf{x}$ has a real multiplier $>1$. In particular
$\mathbf{x}$ is dynamically unstable.
\end{corollary}

\noindent{\it Proof.} By Corollary \ref{cor4.1}, both $h_\pm$ are positive
definite. Therefore $\det H_\pm>0$. Now by \eqref{eq2.8} and
Proposition \ref{pr4.3}     we have $\det(P-I)=-1$. It remains to
use Proposition \ref{pr2.2}.
\qed \medskip

Consider, for example, the  DLS generated by    billiards in a domain in
$\mathbb{R}^{m+1}$ bounded by a hypersurface $M$.
Suppose that the billiard hypersurface $M$ is
symmetric with respect to a hyperplane in $\mathbb{R}^{m+1}$,
for definiteness passing through 0.
Then the symmetry $S = S_e\colon \mathbb{R}^{m+1}\to \mathbb{R}^{m+1}$  is given by
$S_e(x)=x-2e\langle x,e\rangle$, where $e$ is the unit normal vector.

\begin{proposition}
\label{pr4.6}
Let   $x,y\in M$, \,$y=S_e x$, be a pair of symmetric points and let $B=B(y,x)$.
Then the operator  $C=-BS_e\colon T_yM\to T_y^*M$ is symmetric and
positive definite: $\langle C v,v\rangle<0$ for non-zero  $v\in T_yM$.
\end{proposition}

Indeed, by
\eqref{eq2.6},
$$
\langle C v,w\rangle=\frac{-\langle S_e v,w\rangle
+\langle S_e v,e\rangle\langle w,e\rangle}{|x-y|}=
-\frac{\langle v-e\langle v,e\rangle,w-e
\langle w,e\rangle\rangle}{|x-y|}\,,
$$
where we used  $S_e v=v-2e\langle v,e\rangle$.

\begin{corollary}
\label{cor4.4}
Let   $\mathbf{x}$ be an $S_e$-reversible  periodic billiard
trajectory of type $\tau$ such that $\mathbf{y}$ is a nondegenerate
minimum of the length functional $\mathscr{A}_+$. If
$m+ \tau$ is odd, then $\mathbf{x}$ has  a real multiplier
greater than 1.
\end{corollary}

\noindent{\it Proof.}
Consider the case $\tau=0$. By Proposition \ref{pr4.6}, the symmetric
operators $C_1$ and $C_k$ are negative definite. Therefore by Proposition
\ref{pr4.4}, $\mathbf{x}$ is a nondegenerate minimum of $\mathscr{A}$.
Note also that
$\sigma(\mathbf{x}) = (-1)^n = (-1)^{\tau} > 0$. Now it remains to use
Corollary \ref{cor2.6}. The cases $\tau = 1,2$ are analogous.
\qed \medskip

Any billiard is  $S$-reversible for $S=\mathrm{id}$.
Any reversible periodic trajectory
$\mathbf{x}$ is of type~2. By Proposition  \ref{pr2.1},
$\sigma(\mathbf{x})=(-1)^n>0$.

\begin{corollary}
\label{cor4.5}
Any\,   $\mathrm{id}$-reversible billiard trajectory which gives
a nondegenerate
minimum of the functional            $\mathscr A_+$, 
is hyperbolic.
\end{corollary}

This follows from Proposition \ref{pr4.5}.

\section{Hill's formula for a  continuous Lagrangian system}
\label{sec5}

Consider a continuous Lagrangian system $(M,\mathscr{L})$ with the
configuration space $M^m$ and smooth\,\footnote{Actually, $C^2$ is enough.}
$\tau$-periodic Lagrangian $\mathscr{L}(x,\dot x,t)$ on
$TM\times\mathbb{R}$. We assume that    $\mathscr{L}$ is  strictly convex
in velocity $\dot x\in T_xM$.
Then $\tau$-periodic trajectories  are
critical points of the action functional
$$
\mathscr{A}(\gamma)=  \int_0^\tau \mathscr{L}(\gamma(t),\dot\gamma(t), t)\,dt
$$
 on the  space $\Omega$ of  $\tau$-periodic $W^{1,2}$ curves
 $\gamma\colon\mathbb{R}\to M$.
The goal of this section is to prove an analogue of Theorem
\ref{th2.1} for continuous Lagrangian systems.

\subsection{Continuous Hill determinant}

The second variation of the functional
$\mathscr{A}$ at  $\gamma\in\Omega$ is a symmetric bilinear form $h(\xi,\eta)$ on the set
$X$ of  $\tau$-periodic $W^{1,2}$  vector fields $\xi(t)\in E_t =
T_{\gamma(t)} M$ along $\gamma$. It is defined by
$$
h(\xi,\xi)=\frac{d^2}{d\alpha^2}\bigg|_{\alpha=0}
\mathscr{A}(\gamma_\alpha),\qquad
\gamma_0=\gamma,
$$                                where
    $\gamma_\alpha\colon\mathbb{R}\to M$
is  a smooth $\tau$-periodic variation of
$\gamma$.
 Define a positive definite scalar product on $E_t$ by
$$
(v,w)=\langle B(t)v,w\rangle, \qquad
B(t)=\mathscr{L}_{\dot x\dot x}(\gamma(t),\dot\gamma(t),t).
$$

\begin{proposition}
\label{pr5.1}
 $h$ can be  uniquely represented in the form
\begin{equation}
\label{eq5.1}
h(\xi,\eta)=\int_0^\tau\Bigl(\bigl(D\xi(t),D\eta(t)\bigr)
+\bigl(U(t)\xi(t),\eta(t)\bigr)\Bigr)\,dt,
\end{equation}
where $U(t)=U^*(t)$ is a symmetric linear operator and $D$ is a covariant derivative,
that is, a linear   differential operator  such that
\begin{equation}
\label{eq5.2}
\frac{d}{dt}(\xi,\eta)=(D\xi,\eta)+(\xi,D\eta),\qquad
\frac{d}{dt}(f\xi)=\dot f\xi+f\,D\xi
\end{equation}
for smooth vector fields $\xi(t),\eta(t)\in E_t$ and a scalar function $f(t)$.
\end{proposition}

\noindent{\it Proof.}
Let   $\nabla$ be any  covariant derivative.\footnote{
  A covariant derivative is not uniquely defined: for an
  antisymmetric operator $A(t)$, $\nabla+A(t)$ is also a covariant derivative.
  We use a covariant derivative because the derivative is undefined unless $E_t$ is $t$-independent.}
  A standard computation shows that $h$ can be written in the form
$$
h(\xi,\eta)=\int_0^\tau\Bigl(\bigl((\nabla\xi(t),\nabla\eta(t)\bigr)
+\bigl(W(t)\xi(t), \nabla\eta(t)\bigr)
+\bigl(V(t)\xi(t), \eta(t)\bigr)\Bigr)\,dt,
$$
where
$V(t),W(t)$: $E_t\to E_t$ are linear operators, and  $V(t)$ is
symmetric with respect to the metric: $V(t) = V^*(t)$.

By integration by parts $h$ can be represented in the form
  \eqref{eq5.1}, where
$$
D\xi=\nabla\xi+W-W^*,\qquad
U=V-W-W^*=U^*.
$$
Hence  $D$ is also a covariant derivative.
\qed \medskip

Note that $D$ and $U$ are invariantly determined by $h$,
that is,   they are coordinate
independent and do not change by a calibration of the Lagrangian.

Equations \eqref{eq5.2} imply that $D$ is skew-symmetric
relative to the $L^2$ scalar product
$$
(\xi,\eta)_2=\int_0^\tau\bigl(\xi(t),\eta(t)\bigr)\,dt.
$$
Therefore
\begin{equation}
\label{eq5.3}
h(\xi,\eta)=\bigl((-D^2+U)\xi,\eta\bigr)_2=(\mathbf{H}\xi,\eta)_2,
\end{equation}    where
    $\mathbf{H}=-D^2+U$ is the Hessian of  $\mathscr{A}$ with respect
to the
  $L^2$-metric.

The variational system of the periodic trajectory $\gamma$ has the
form
\begin{equation}
\label{eq5.4}
D^2\xi(t)=U(t)\xi(t).
\end{equation}
This is a linear Lagrangian system.
We use the following definition.

\begin{definition}
\label{def5.1}
Let   $E=\{E_t\}_{t\in\mathbb{R}/\tau\mathbb{Z}}$ be a smooth vector bundle.
Suppose it is equipped with a metric   ${(\,\cdot\,{,}\,\cdot\,)}$ compatible
with a covariant derivative  $D$ 
and a symmetric linear operator  $U(t)\colon E_t\to E_t$.
Denote by $(E,\Lambda)$ the linear Lagrangian system
with the quadratic Lagrangian
\begin{equation}
\label{eq5.5}
\Lambda(\xi,D\xi)=\frac{1}{2}\,(D\xi,D\xi)+\frac{1}{2}\,(U\xi,\xi)
\end{equation}
and Lagrange's equations \eqref{eq5.4}.
\end{definition}

Trajectories $\xi(t)$, $0\le t\le \tau$, of the system $(E,\Lambda)$ are extremals of the quadratic
action functional
\begin{equation}
\label{eq5.6}
\frac{1}{2}\, h(\xi,\xi)=\int_0^\tau \Lambda(\xi,D\xi)\,dt
\end{equation}
for variations with fixed $\xi(0),\xi(\tau)$. Thus $h(\xi,\phi)=0$
for any smooth $\phi(t)\in E_t$ such that $\phi(0)=0$ and $\phi(\tau)=0$.

The system $(E,\Lambda)$ is the linearization of $(M,\mathscr{L})$ at
$\gamma$. In what follows we can forget about the non-linear Lagrangian
system $(M,\mathscr{L})$ and work with the  linear system $(E,\Lambda)$.

Let  $P:W\to W$ be the linear Poincar\'e map of the  trajectory
$\gamma$.
Since a solution $\xi(t)$ of the variational system  is uniquely
determined by $(\xi(0),D\xi(0))$, $W$ can be identified with $E_0\oplus E_0$.
Then $P$   is the monodromy operator of
the variational system:
$$
P\bigl(\xi(0),D\xi(0)\bigr)=\bigl(\xi(\tau),D\xi(\tau)\bigr).
$$

Define the  $W^{1,2}$-scalar product on the Hilbert space  $X$
by
$$
\llangle\xi,\eta\rrangle=(D\xi,D {\eta})_2+(\xi,{\eta})_2=
(\mathbf{B}\xi,\eta)_2,\qquad
\mathbf{B}=-D^2+I.
$$
Then  $h(\xi,\bar\eta)=\llangle H\xi,\bar\eta\rrangle$,
where the self-adjoint operator $H =\mathbf{B}^{-1}\mathbf{H}$
is the Hessian of   $\mathscr{A}$ with respect to the
  $W^{1,2}$-scalar product.

We have $H=I+K$, where  $K=(-D^2+I)^{-1}(U-I)$ is compact, with eigenvalues
                          $\lambda_k=O(k^{-2})$, \,$k=1,2,\dots$, so that
$$
\operatorname{tr} |K|=\sum_{k=1}^\infty|\lambda_k|<\infty.
$$
Thus, the           {\it Hill determinant}
$$
\det H=\prod_{k=1}^\infty(1+\lambda_k)
$$
converges absolutely.

Let   $Q\colon E_0\to E_0$ be the monodromy operator
of the equation of parallel transport:
$$
Q\eta(0)=\eta(\tau),\qquad
D\eta(t)=0.
$$

\begin{theorem}
\label{th5.1}
$\det(I-P)=\sigma(-1)^m\beta \det H$, where
\begin{equation}
\label{eq5.7}
\beta=e^{-m\tau}{\det}^{2}(e^\tau I-Q),\qquad
\sigma=\det Q.
\end{equation}
\end{theorem}

Since $Q$ is an orthogonal operator,  $\beta> 0$ and $\sigma=\pm 1$
depending on whether the bundle $E$ is orientable, that is,
 if the trajectory $\gamma$
preserves or reverses  orientation. If $M$ is orientable, then
$\sigma=1$ always.

Theorem \ref{th5.1} follows from a more general result of the
next subsection.

\begin{example}
\label{ex5.1}{\rm
Suppose the Lagrangian system has one degree of freedom and
the bundle $E$ is  trivial. Then  $Q=1$ and 
  $\mathbf{H}\xi=-\ddot\xi+a(t)\xi$.  Since $\det P=1$,
$$
\rho^{-1} \det(\rho I-P)=\rho+\rho^{-1}-2+\det(I-P).
$$
If   $\rho$ is a multiplier, Theorem      \ref{th5.1} gives
\begin{equation}
\label{eq5.8}
\det H=\frac{\rho+\rho^{-1}-2}{e^\tau+e^{-\tau}-2}\,.
\end{equation}
Set     $\tau=2\pi$ and represent the operator~$H$
in the basis  $\{e^{int}\}$. If
$$
a(t)=\sum_{n\in\mathbb{Z}}a_ne^{int},\qquad
\xi(t)=\sum_{n\in\mathbb{Z}}\xi_ne^{int},
$$
then
$$
H\xi=(-D^2+I)^{-1} \big( -D^2+U \big)\,\xi=\sum_{n\in\mathbb{Z}}(n^2+1)^{-1}
\biggl(n^2+\sum_{k\in\mathbb{Z}} a_ke^{ikt}\biggr) \xi_n\,e^{int}.
$$
Hence \eqref{eq1.3} is the matrix of     $H$, and  \eqref{eq5.8}
gives Hill's formula \eqref{eq1.2}.}
\end{example}

\subsection{Relation to Hill's formula for discrete Lagrangian systems}
\label{ssec5.2}
A continuous Lagrangian system locally, near a periodic orbit
$\gamma$,  defines  a discrete Lagrangian system. Take a partition
$0=t_0<t_1<\dots<t_n=\tau$ of $[0,\tau]$ and let   $x_i=\gamma(t_i)$.
If the  $\Delta t_i=t_{i+1}-t_i$ are small enough, the points $x_i$, $x_{i+1}$ are
non-conjugate   along $\gamma$. Then there is a neighbourhood $U_i$
of $(x_i,x_{i+1})$ in $M^2$ such that for each $(x,y)\in U_i$ there exists a
unique trajectory $u_{x,y}$: $[t_i,t_{i+1}]\to M$ close to
$\gamma\big|_{[t_i,t_{i+1}]}$ and joining $x$ and $y$.
 Define a discrete Lagrangian $L_i$ on $U_i$ by
$$
L_i(x,y)=\mathscr{A}(u_{x,y}).
$$
The discrete action functional  $A$ is
$$
A(\mathbf{y})=\sum_{i=1}^n L_i(y_i,y_{i+1}),\qquad
(y_i,y_{i+1})\in U_i, \quad
y_{n+1}=y_1.
$$
Then  $\mathbf{x}$
is a critical point of $A$, that is,   the periodic orbit of the DLS
corresponding to the periodic orbit $\gamma$
 of the continuous Lagrangian system. It is easy to see that
 $$
 B_i=-\partial_{12} L_i(x_i,x_{i+1})=\frac{1}{\Delta t_i}\bigl(B(t_i)
+O(\Delta t_i)\bigr).
 $$
 Thus $L_i$ is a  discrete Lagrangian. 

 The definition of the Hill determinant $\det H=\det(\mathbf{B}^{-1}\mathbf{H})$
for a CLS is
similar to  \eqref{eq2.8}      for a DLS.
 Discretization of the operator $\mathbf{H}=-D^2+U$
 corresponds to the operator $\mathbf{H}$ in  \eqref{eq2.8}.
 However, the operator $\mathbf{B}=-D^2+I$ does not correspond to the operator
 $\mathbf{B}$ in \eqref{eq2.11}.

The choice of $\mathbf{B}$ is natural for DLS, but not so for CLS, where
instead of $I$ we could add almost anything. This is the reason for the
strange  coefficient $\beta$ in \eqref{eq5.7}. If we use an analogue
of discrete $\mathbf{B}$, then $\mathbf{B}^{-1}\mathbf{H}$ will be unbounded.

\subsection{Generalized Hill determinant}
\label{ssec5.3}
For a given $\rho\in S^1$ let $X_\rho$ be the vector space of
complex     $\rho$-quasiperiodic locally       $W^{1,2}$ vector fields
$\xi(t)\in E_t$ such that  $\xi(t+\tau)=\rho\xi(t)$. Define a  Hermitian
$\rho$-index form        \cite{28},~\cite{27} on~$X_\rho$
by  \eqref{eq5.1}:
\begin{equation}
\label{eq5.9}
h(\xi,\bar\eta)=\int_0^\tau\Bigl(\bigl(D\xi(t),D\overline{\eta}(t)\bigr)+
\bigl(U(t)\xi(t),\overline{\eta}(t)\bigr)\Bigr)\,dt.
\end{equation}
We also denote by  $X$ the complexification   $X=X^{\mathbb{C}}$,
that is, the set of complex      $\tau$-periodic      $W^{1,2}$-vector fields
along $\gamma$. For definiteness choose    $\ln\rho$ so   that
$$
0\le\operatorname{Im}\ln\rho<2\pi,
$$
and let   $\mu=\tau^{-1}\ln\rho$.
Identifying  $X$ and~$X_\rho$ by the map
\begin{equation}
\label{eq5.10}
X\ni\xi\mapsto e^{\mu t}\xi(t)\in X_\rho,
\end{equation}
we obtain a Hermitian form~$h_\rho$ on~$X$:
\begin{align*}
h_\rho(\xi,\bar\eta)&=h(e^{\mu t}\xi,\overline{e^{\mu t}\eta}\,)
=\bigl((D+\mu I)\xi,\overline{(D+\mu I)\eta}\,\bigr)_2+(U\xi,\bar\eta)_2
\\
&=-((D+\mu I)^2\xi,\bar\eta)_2+(U\xi,\bar\eta)_2
=\bigl((-(D+\mu I)^2+U)\xi,\bar\eta\bigr)_2
\\
&=(\mathbf{H}_\rho\xi,\bar\eta)_2,\quad\text{where} \
\mathbf{H}_\rho=-(D+\mu I)^2+U.
\end{align*}
(We used that               $\bar\mu=-\mu$ and  $D$ is real and
antisymmetric).
Define the  $\rho$-Hessian operator $H_\rho\colon X\to X$ by
$h_\rho(\xi,\bar\eta)=\llangle H_\rho\xi,\bar\eta\,\rrangle$. Then
\begin{equation}
\label{eq5.11}
H_\rho=\mathbf{B}^{-1}\mathbf{H}_\rho=(-D^2+I)^{-1}\bigl(-(D+\mu I)^2+U\bigr).
\end{equation}

We henceforth assume that $\rho\ne 0$ may take any complex values.
The operator $H_\rho$ is self adjoint for $|\rho|=1$, but not in
general. Although this is not a trace class operator: 
 $\operatorname{tr}|H_\rho-I|$ diverges for $\rho\ne 1$,
we can,
following Poincar\'e \cite{2}, define the generalized Hill
determinant $\det H_\rho$ by means of the finite-dimensional
approximation
\begin{equation}
\label{eq5.12}
\det H_\rho=\lim_{N\to\infty}\det H_\rho^{(N)},\qquad
H_\rho^{(N)}=P_NH_\rho P_N^*\colon X^{(N)}\to X^{(N)},
\end{equation}
where $P_N$ is the orthogonal projection onto the finite-dimensional
eigenspace $X^{(N)}$ of the operator $D$ corresponding to the eigenvalues in
$$
\Lambda_N=\{\nu\in\Lambda=\sigma(D):|\nu|\le N\}.
$$

\begin{theorem}
\label{th5.2}                            The
determinant  \eqref{eq5.12} converges and
\begin{equation}
\label{eq5.13}
\rho^{-m}\det(P-\rho I)=\sigma(-1)^m\beta\det H_\rho.
\end{equation}
\end{theorem}

We present the proof of Theorem  \ref{th5.2} in \S\,\ref{ssec5.5}.

For $\rho=1$ we obtain Theorem \ref{th5.1}.

\subsection{Some applications}
\label{ssec5.4}
Suppose that     $\mathscr{L}(x,\dot x)=(\dot x,\dot x)/2$, where
${(\,\cdot\,{,}\,\cdot\,)}$ is a Riemannian metric on~$M$.
The periodic orbit $\gamma$ is
a closed geodesic. The quadratic Lagrangian of the variational
system  has the form \eqref{eq5.5}, where $D$ is the
Levi-Civita covariant derivative along $\gamma$ and
$U(t)\xi=R\bigl(\xi,\dot\gamma(t)\bigr)\xi$ with $R$ the curvature tensor.

The variational system  has a periodic solution
$\dot\gamma(t)$
 and a first integral $ (D\xi,\dot\gamma)= \dfrac{d}{dt}(\xi,\dot\gamma)$
periodic in time. Hence $P$ has two unit multipliers,
and  $\det H = \det(I - P) = 0$. Let us present a reduced version of
Hill's formula. More general results will be proved in the next subsection
(see Corollary \ref{cor6.1}).

Let $E_t^\perp=\{u\in E_t:(u,\dot\gamma(t))=0\}$.
If $\xi(t)\in E_t^\perp$, then
 $D\xi(t),U(t)\xi(t)\in E_t$.  Denote by
$H^\perp\colon X^\perp\to X^\perp$  the restriction of  $H$  to the
invariant subspace
$$
X^\perp=\{\xi\in X: \xi(t) \in E_t^\perp\}.
$$
Let   $P^\perp\colon W^\perp\to W^\perp$,
\,$W^\perp=E_0^\perp \times E_0^\perp$,
be the monodromy
operator corresponding to solutions $\xi(t)\in E_t^\perp$  of the variational
system. Let  $Q^\perp\colon E_0^\perp\to E_0^\perp$ be the map of parallel
transport
along 
 $\gamma$. Applying
Theorem \ref{th5.1}  to the linear Lagrangian system
$(E^\perp,\Lambda^\perp)$, where
$\Lambda^\perp=\Lambda\big|_{E^\perp}$, we obtain the following result
  \cite{5}. 

\begin{corollary}
\label{cor5.1}
Hill's formula for the reduced system has the form
$$
\det(P^\perp-I)=\sigma(-1)^{m-1}\beta^\perp \det H^\perp,\qquad
\beta^{\perp}=e^{-(m-1)\tau}{\det}^2(Q^\perp-e^\tau I).
$$
\end{corollary}

Let us formulate another corollary to Theorem \ref{th5.2}.
For complex  $\rho$
let   $X^\perp$ be the space of complex vector fields       $\xi(t)\in E_t$
and let $H_\rho^\perp=H_\rho\big|_{X^\perp}$. Then
\begin{equation}
\label{eq5.14}
\rho^{-(m-1)}\det(P^\perp-\rho I)=\sigma(-1)^{m-1}\beta^\perp
\det H_\rho^\perp.
\end{equation}

Let us present a proof of \eqref{eq5.14} from  \cite{5},
which will be generalized in \S\,\ref{ssec6.1}. Let
$Z=\{\xi\in X:\xi(t)=\lambda(t)\dot\gamma(t)\}$. We write any
$\xi\in X$ as     $\xi(t)=\eta(t)+\lambda(t)\dot\gamma(t)$, where
    $\eta\in X^\perp$. Then
\begin{equation}
\label{eq5.15}
h(\xi,\bar\xi\,)=h(\eta,\bar\eta)+\int_0^\tau |\dot\lambda|^2\,dt.
\end{equation}
Hence   $H_\rho=H_\rho\big|_{Z}\oplus H_\rho^\perp$  and
$\det H_\rho= \det H_\rho\big|_{Z}\det H_\rho^\perp$,    where
    $\det H_\rho\big|_{Z}$
is the
Hill determinant for the  system with the quadratic Lagrangian
$|\dot\lambda|^2/2$.  The characteristic polynomial of the
corresponding monodromy matrix is $(\rho-1)^2$. Thus by
\eqref{eq5.8},
$$
\det H_\rho\big|_{Z}=-\frac{e^\tau(\rho-1)^2}{\rho(e^\tau-1)^{2}}\,.
$$
But
$$
\det(\rho I-P)=(\rho-1)^2\det(\rho I-P^\perp),\qquad
\det(e^\tau I-Q)=(e^\tau-1)\det(e^\tau I-Q^\perp),
$$
which implies  \eqref{eq5.14}.

Next we discuss applications to stability of periodic trajectories,
 similar to the discrete
case. For $\rho\in S^1$ define the   $\rho$-index (\cite{28},  \cite{27})
$\operatorname{ind}_\rho\gamma$   of a periodic trajectory $\gamma$  as the index
of the Hermitian form $h_\rho$.  Then $\operatorname{ind} \gamma =
\operatorname{ind}_1 \gamma$
is the Morse index of $\gamma$. It  equals  the number of negative
eigenvalues of the operator $H$.  If $\rho$ is not an eigenvalue of
$P$,
$$
(-1)^{\operatorname{ind}_\rho\gamma}=\operatorname{sign}\det H_\rho=
\sigma(-1)^m \operatorname{sign}\bigl(\rho^{-m} \det(\rho I-P)\bigr).
$$
The argument of the  $\operatorname{sign}$ function is real for  $|\rho|=1$
since the characteristic polynomial is reciprocal.

The next result is proved in  \cite{5}.

\begin{corollary}
\label{cor5.2}
Suppose the trajectory   $\gamma$ is nondegenerate and
  $\sigma(-1)^{m+\operatorname{ind}\gamma}<0$. Then  $\gamma$ has
a real multiplier $\rho>1$.
\end{corollary}

Indeed, the characteristic polynomial $F(\rho)=\det(\rho I-P)$
satisfies $F(1)<0$ and $F(+\infty)=+\infty$. Hence $F$ has a real
root $\rho > 1$.

Corollary \ref{cor5.2} is not true if $\gamma$ is degenerate.
Suppose for example, that $\gamma$ is a closed geodesic.\footnote{The general
degenerate case is discussed
in \S\,\ref{sec6}.} Then  \eqref{eq5.15}
implies $\operatorname{ind} h_\rho=\operatorname{ind} h_\rho^\perp$, and
so
$$
(-1)^{\operatorname{ind}_\rho\gamma}=\operatorname{sign}\det H_\rho^\perp=
\sigma(-1)^{m-1} \operatorname{sign}\bigl(\rho^{1-m} \det(\rho I-P^\perp)\bigr).
$$

\begin{corollary}
\label{cor5.3}           Suppose
the closed geodesic  $\gamma$ is nondegenerate and
$\sigma (-1)^{m+\operatorname{ind}\gamma} >0$.
Then the characteristic polynomial
$F(\rho)=\det(\rho I-\nobreak P)$ has a real root $\rho > 1$.
Therefore,  $\gamma$ is exponentially unstable.
\end{corollary}

This is proved in \cite{5} and  \cite{29} using Hill's formula and
also recently in  \cite{11} using the theory of Maslov index  \cite{7}. 
In particular, nondegenerate closed geodesics of locally minimal
length on an even-dimensional orientable manifold are exponentially
unstable. Degenerate geodesics are linearly unstable, but
in general instability will not be exponential and so has no relevance
for applications to Lyapunov stability. 

Suppose $m = 1$ and let the $2\tau$-periodic trajectory $\gamma^2$ be
$\gamma$ traversed twice. {\it If $\gamma^2$ is nondegenerate (that is, $\pm1$
are not multipliers),
 then
$\gamma$ has hyperbolic (elliptic) type if and only if $\operatorname{ind}\gamma^2$ is even
(odd).}

Indeed, $\gamma$ and $\gamma^2$ are
simultaneously elliptic or hyperbolic. The multipliers  of $\gamma^2$ are
  squares of the multipliers of $\gamma$. Hence,
 $\gamma^2$ is hyperbolic if and only if its   multipliers
 are real and positive, or, equivalently,
$\operatorname{sign}\det(I-P^2)=(-1)^{1+\operatorname{ind}\gamma^2}=-1$.
Similarly, the ellipticity of
$\gamma^2$ is equivalent to   $(-1)^{1+\operatorname{ind}\gamma^2} =
1$.

For the geodesic case, we obtain the following result of Poincar\'e. Let
$\gamma$ be a closed geodesic on a 2-dimensional Riemannian
manifold. If $\gamma^2$ is nondegenerate, then $\gamma$ has
hyperbolic (elliptic) type if and only if $\operatorname{ind}\gamma^2$ is even
(odd).

Suppose now $\rho=-1$.  Then $(-1)^{\operatorname{ind}_{-1}\gamma}=
\sigma \operatorname{sign}
F(-1)$. Thus, if
$\operatorname{ind}_{-1}\gamma$ is odd, there exists a real multiplier $\rho<-1$.
Note that  the space $X_{-1}$ corresponds to antiperiodic variations
such that $\xi(\tau)=-\xi(0)$. Since $2\tau$-periodic vector fields
are sums of $\tau$-periodic and $\tau$-antiperiodic,
$$
\operatorname{ind}_{-1}\gamma=\operatorname{ind}\gamma^2-\operatorname{ind}\gamma.
$$
Thus, if      $\operatorname{ind}\gamma$ and~$\operatorname{ind}\gamma^2$
are not even or odd
simultaneously, then $\gamma$ is unstable.

\subsection{Proof of Theorem~\ref{th5.2}}
\label{ssec5.5}
We follow  \cite{5}, see also  \cite{6}. The method goes
back to Poincar\'e's proof (\cite{2},~\cite{30})
of Hill's result \cite{1}.

The real skew-Hermitian operator  $D=\,\overline{\!D}=-D^*$ has
compact resolvent       $(D+\mu I)^{-1}$. Its spectrum
$\Lambda=\sigma(D)\subset i\mathbb{R}$ coincides with the set
of characteristic exponents of the equation $D\eta(t)=0$ of parallel transport.
Thus
$$
\Lambda=\{\nu:\det(Q-e^{\tau\nu}I)=0\}.
$$
If   $\nu\in\Lambda$, then  $-\nu$ and~$\nu+\omega$ belong to  $\Lambda$,
where $\omega=2\pi i/\tau$.

Let   $\rho_1,\dots,\rho_m$ be the roots of   $\det(Q-\rho I)=0$. 
Since   $|\rho_j|=1$, we may represent them as
$\rho_j=e^{\nu_j\tau}$, where $0\le \operatorname{Im}\nu_j<2\pi/\tau$. Then
\begin{equation}
\label{eq5.16}
\Lambda=\bigcup_{j=1}^m(\nu_j+\omega\mathbb{Z}).
\end{equation}

First suppose that       $\mu\notin \Lambda$. Then  $H_\rho=ST$, where
$$
S(\mu)=-(-D^2+I)^{-1}(D+\mu I)^2,
\qquad
T(\mu)=I-(D+\mu I)^{-2}U.
$$
Since   $P_ND=DP_N$, by  \eqref{eq5.12} we have
$$
\det H_\rho=\det S\det T.
$$ The
finite-dimensional approximation  \eqref{eq5.12} of the determinant
$$
\det T(\mu)=f(\mu)=\lim_{N\to\infty} \det(P_N T P_N^*),
$$
converges absolutely for $\mu\notin\Lambda$ since
$$
\operatorname{tr}|(D+\mu I)^{-2}U|<\infty.
$$
Hence $f$ is a holomorphic function on
$\mathbb{C}\setminus\Lambda$ having at points in $\Lambda$ poles of
multiplicity not greater than double the multiplicity of the
corresponding points of the spectrum of $D$.

The function $f$ is periodic: $f(\mu + \omega) \equiv  f(\mu)$.
Indeed, if $\xi\in X^{\mathbb{C}}$, then $e^{\omega t}\xi
\in X^{\mathbb{C}}$ and
$$
\bigl(I-(D+\mu I)^{-2}U \bigr) e^{\omega t}\xi=
e^{\omega t}\bigl(I-(D+(\mu+\omega)I)^{-2}U \bigr)\xi,
$$
so                $T(\mu)$ and  $T(\mu+\omega)$ are similar. Thus
$f(\mu)=\phi(e^{\mu\tau})$,
where  $\phi(\rho)$ is a  meromorphic function
having poles at the roots $\rho_1,\dots,\rho_m$ of $\det(\rho I -
Q)$. The multiplicity of the pole is at most twice the multiplicity
of the corresponding root.

Hence  there exists a polynomial $g(\rho)$ of degree   $\le 2m - 1$
such that the functions $\phi(\rho)$ and $g(\rho)\det^{-2}(\rho
I-Q)$ have the same principal parts of the Laurent expansion at each
pole. Since $\phi(\rho) \to 1$ as $|\rho|\to +\infty$, by
Liouville's theorem,
\begin{equation}
\label{eq5.17}
\phi(\rho)=1+g(\rho){\det}^{-2}(\rho I-Q).
\end{equation}

The determinant $\det S$ converges conditionally. By~\eqref{eq5.12},
\begin{align*}
&\det\bigl(-(-D^2+I)^{-1}(D+\mu I)^2\bigr)=\lim_{N\to\infty}\,
\prod_{\nu\in\Lambda_N}\frac{(\nu+\mu)^2}{\nu^2-1}
\\
&\qquad=\lim_{N\to\infty}(-\mu^2)^k\prod_{\nu\in\Lambda_N,\,i\nu>0}
\biggl(\frac{\nu^2-\mu^2}{\nu^2-1}\biggr)^2
=(-1)^k\lim_{N\to\infty}\,\prod_{\nu\in\Lambda_N}
\frac{\nu^2-\mu^2}{\nu^2-1}\,,
\end{align*}                                        where
  $k$ is the multiplicity of zero in the spectrum of $D$. We have used
that
$\Lambda=-\Lambda$.

From  \eqref{eq5.16}
it follows that the last product converges
absolutely. Hence it is a  holomorphic function of
$\rho_1,\dots,\rho_m,\rho$ for $\rho_j\ne e^{\pm\tau}$ and $\rho\ne
0$. To compute the product, we will use Euler's formula (see, for example,
\cite{30}):
$$
\prod_{n\in\mathbb{Z}}\biggl(1-\frac{\mu^2}{(\nu+\omega n)^2}\biggr)
=\frac{\cosh\mu\tau-\cosh\nu\tau}{1-\cosh\nu\tau}\,,\qquad
\nu\notin \omega\mathbb{Z}.
$$

Suppose first that       $\rho_j\ne 1$ and~$\rho_i\ne \rho_j $
for $i\ne j$.
Equivalently, $\nu_j\notin\omega\mathbb{Z}$ and
$\nu_i-\nu_j\notin\omega\mathbb{Z}$ for $i\ne j$. Then by  \eqref{eq5.16},
\begin{align*}
\prod_{\nu\in\Lambda}\frac{\nu^2-\mu^2}{\nu^2-1}
&=\prod_{\nu\in\Lambda}\biggl(1-\frac{\mu^2}{\nu^2}\biggr)
\biggl(1-\frac{1}{\nu^2}\biggr)^{-1}
\\
&=\prod_{j=1}^m\,\prod_{n\in\mathbb{Z}}
\biggl(1-\frac{\mu^2}{(\nu_j+\omega n)^2}\biggr)
\biggl(1-\frac{1}{(\nu_j+\omega n)^2}\biggr)^{-1}
\\
&=\prod_{j=1}^m \frac{\cosh\mu\tau-\cosh\nu_j\tau}{\cosh\tau-\cosh\nu_j\tau}
=\prod_{j=1}^m \frac{\rho+\rho^{-1}-\rho_j-\rho_j^{-1}}
{e^\tau+e^{-\tau}-\rho_j-\rho_j^{-1}}
\\
&=\prod_{j=1}^m \frac{e^\tau(\rho-\rho_j)^2}{\rho(e^\tau-\rho_j)^2}
=\frac{e^{m\tau}\det^2(\rho I-Q)}{\rho^{m}\det^2(e^\tau I-Q)}\,.
\end{align*}
By continuity this holds  for any $\rho_1,\dots,\rho_m\ne
e^{\pm\tau}$ and $\rho\ne 0$. Hence
$$
\det S=(-1)^k\frac{e^{m\tau}\det^2(\rho I-Q)}{\rho^{m}\det^2(e^\tau I-Q)}\,.
$$
By \eqref{eq5.12} and  \eqref{eq5.17},
\begin{equation}
\label{eq5.18}
\rho^m\det H_\rho=(-1)^k \beta^{-1}\bigl({\det}^2(\rho I-Q)+g(\rho)\bigr).
\end{equation}
Thus $G(\rho)=\rho^m\det H_\rho$ is a polynomial of degree $2m$ in
$\rho$ with leading coefficient $(-1)^k  \beta^{-1}$.

We claim that the polynomials $G(\rho)$ and $F(\rho)=\det(P-\rho
I)$ have the same roots.  It is sufficient to prove that if
$F(\rho)=0$, then $G(\rho)=0$ and the root has at least the same multiplicity.

If $F(\rho)=0$, there exists a non-zero $\tau$-periodic vector field
$\xi\in X$ such that $(-D^2 +U)e^{\mu t}\xi(t) = 0$. Thus $H_\rho\xi=0$.
Suppose first that $\rho$ is not an eigenvalue of $Q$, that is,
 $\mu\notin\Lambda$. Then $H_\rho=ST$, where $S$ is invertible
and $\operatorname{tr}|T-I|<\infty$. Hence $T\xi=0$ implies $\det T=0$ (see,
for example, \cite{31}). Then   $\det H_\rho=\det S\det T=0$, and so
$G(\rho)=0$.

If $\mu\in\Lambda$, we can repeat the same argument replacing $S$ and $T$,
for instance,  by
$$
\widetilde S=(-D^2+I)^{-1}\bigl(-(D+\mu I)^2+I\bigr),\qquad
\widetilde T=\bigl(-(D+\mu I)^2+ I\bigr)^{-1}\bigl(-(D+\mu I)^2+U\bigr).
$$

We have proved that
$$
G(\rho)=(-1)^{k}\beta^{-1}F(\rho).
$$
It remains to show that $(-1)^k=\sigma(-1)^m$.   Indeed, $k$ is the
dimension of the subspace on which the orthogonal operator $Q$ is
the identity, while $\sigma = (-1)^n$, where $n$ is the dimension of
the subspace on which $Q$ is a reflection. Since the dimension $m -
k - n$ of the complementary subspace is even, \eqref{eq5.13} is
proved.

\section{Degeneracy in Hill's formula}
\label{sec6}

In this section we consider the case when the periodic orbit
$\gamma$ is degenerate, that is, the variational system has a non-zero
$\tau$-periodic solution $\zeta$. Equivalently,  the linear
Poincar\'e map $P$ has multiplier 1.
 Usually,
this happens  if the Lagrangian system has a time periodic first
integral $\mathscr J$  which is nondegenerate on  $\gamma$.
Then, as proved by Poincar\'e, the variational system has  a non-zero
periodic solution and a non-trivial linear time periodic first
integral which is the linearization of $\mathscr J$ at $\gamma$. Here are two
standard examples.

\textbf{1. Autonomous Lagrangian system.} Then
$(M,\mathscr{L})$ has the energy integral
$$
\mathscr{H}(x,\dot x)=\langle p,\dot x\rangle-\mathscr{L}(x,\dot x),\qquad
p=\mathscr{L}_{\dot x}(x,\dot x).
$$
The variational system of a periodic orbit $\gamma$ has a periodic
solution $\zeta(t)=\dot\gamma(t)$. A particular case is a closed
geodesic in a Riemannian metric.

\textbf{2. A Lagrangian system with symmetry.} Suppose the Lagrangian
system  $(M,\mathscr{L})$ admits a symmetry group    $\psi_s\colon
M\to M$, \,$s\in\mathbb{R}$, preserving  $\mathscr{L}$. Let
  $\mathbf{w}(x)=\dfrac{d}{ds}\bigg|_{s=0}\psi_s(x)$
be the corresponding symmetry field. Then
$$
\mathscr{J}(x,\dot x,t)=\langle p,\mathbf{w}(x)\rangle,\qquad
p=\mathscr{L}_{\dot x}(x,\dot x,t),
$$
is the Noether first integral. The variational system of a
periodic orbit $\gamma$ has a $\tau$-periodic solution
$\zeta(t)=\mathbf{w}(\gamma(t))$. Here is one concrete example.

\textbf{Planar 3-body problem.} Here
$$
\mathscr{L}(x,\dot x)=\frac{1}{2}\sum_{i=1}^3m_i|\dot x_i|^2+
\sum_{i\ne j}\frac{m_im_j}{|x_i-x_j|}\,,\qquad
x_i\in\mathbb{R}^2.
$$
Fix the centre of mass at the origin, so that
$$
M=\biggl\{x=(x_1,x_2,x_3)\in(\mathbb{R}^2)^3: \sum_{i=1}^3m_ix_i=0, \
x_i\ne x_j\biggr\}.
$$
Rotations of  $\mathbb{R}^2$ preserve  $\mathscr{L}$ and the
corresponding symmetry field is         $\mathbf{w}(x)=(Jx_1,Jx_2,Jx_3)$,
\,$J=\begin{pmatrix} 0 & -1 \\ 1 & 0 \end{pmatrix}$.
The Noether integral is the angular momentum
$$
\mathscr{J}(x,\dot x)=\sum_{i=1}^3\langle Jx_i,m_i\dot x_i\rangle.
$$
Since the system is autonomous, we have double degeneracy of any
periodic orbit on which $\mathscr{H}$  and $\mathscr{J}$ are independent:
1 is an eigenvalue of $P$ with multiplicity at least 4.  Stationary periodic
solutions (orbits of the symmetry group) have lower-order  degeneracy. 

In several recent years many periodic solutions for the 3-body problem have
been found by variational methods  \cite{13}. However, we see that the
ordinary Hill formula is degenerate for them. In this section we put forward
an approach to this problem.           

\textbf{3.
General Hamiltonian commutative symmetry.} These examples are particular
cases  of Hamiltonian symmetries.
 Let us look at the  system $(M,\mathscr{L})$   from
the Hamiltonian point of view.  Let $\mathscr{H}$ be the Hamiltonian.
Suppose that the system admits an   algebra $\mathfrak{g}$ of Hamiltonian
symmetry fields $\mathbf{v}$
generated by integrals $\mathscr{J}_\mathbf{v}$.  Let $\gamma$ be a
$\tau$-periodic solution in the phase space and $P $ the corresponding
monodromy operator.  For any $\mathbf{v}\in\mathfrak{g}$, $\zeta(t)=\mathbf{v}
(\gamma(t))$  is a periodic solution of the variational system. Therefore
\begin{equation}
\label{eq6.1}
Pw=w, \qquad
w=\mathbf{v}(\gamma(0)), \quad
\mathbf{v}\in \mathfrak{g}.
\end{equation}
If the system is autonomous, the Hamiltonian vector field of the system
$(M,\mathscr{L})$ will be in  $\mathfrak{g}$, and the corresponding eigenvector
  $w$ is  $\dot\gamma(0)$.

In the present paper we consider only the case when the $k$-dimensional algebra
$\mathfrak{g}$ is commutative. Then the corresponding eigenspace $V\subset
\operatorname{Ker}(P-I)$ is isotropic,
and the multiplicity of  eigenvalue 1 is at least $2k$. Hamiltonian reduction
makes it possible to remove this degeneracy, but then the reduced system loses
the natural Lagrangian structure.

The classical way to remove autonomous degeneracy is to pass from the
Hamilton   action functional to the Maupertuis action functional on
the energy level \cite{32}.  The classical way to remove symmetry
degeneracy in a Lagrangian system is the Routh method \cite{32}.  We
briefly describe it here.

Suppose the system $(M,\mathscr{L})$ admits  $k$
commuting independent symmetry fields     $\mathbf{w}^1,\dots,\mathbf{w}^k$
on~$M$:
$$
[\mathbf{w}^\alpha,\mathbf{w}^\beta]=0, \qquad
\alpha,\beta=1,\dots,k.
$$
The corresponding flows of symmetry        $\psi_{s}^\alpha$ commute. Let    $G$ be
the (local) commutative group acting on  $M$
by  $x\mapsto\psi^1_{s_1}\circ\dots\circ \psi^k_{s_k}(x)$,
\,$s\in\mathbb{R}^k$. Suppose that     $\widetilde M=M/G$ is
a smooth manifold and  $\pi\colon M\to\widetilde M$ a smooth fibration.

The Noether integrals  $\mathscr{J}^\alpha=\langle
p,\mathbf{w}^\alpha\rangle$ give a vector integral
$\mathscr{J}(x,\dot x,t)\in\nobreak\mathbb{R}^k$. The Routh method reduces
the Lagrangian system $(M,\mathscr{L})$ with fixed value
$\mathscr{J}=c\in\mathbb{R}^k$ of the Noether integral to a Lagrangian system
$(\widetilde M,\,\widetilde{\!\mathscr{L}})$ on the reduced configuration
space   $\widetilde M$.

For $c=0$ the reduced Lagrangian is defined by\,\footnote{Recall that
we assume summation in repeated Greek indices.}
\begin{equation}
\label{eq6.2}
\widetilde{\!\mathscr L}(x,\dot
x,t)=\min_{s\in\mathbb{R}^k}\mathscr{L} \bigl(x,\dot x+s_\alpha
\mathbf{w}^\alpha(x),t\bigr),
\end{equation}
provided the minimum exists, for example,
$\mathscr{L}$ is superlinear in velocity. Since
$\,\widetilde{\!\mathscr L}$ depends only on  $\tilde x=\pi(x)$ and
$\tilde {\dot x}=d\pi(x)\dot x$, it can be regarded as a function
on $T\widetilde M\times\mathbb{R}$.

For $c\ne 0$ take closed       $G$-invariant    1-forms $\nu_1,\dots,\nu_k$
on~$M$ such that  $\nu_\alpha(\mathbf w^\beta)\equiv
\delta_\alpha^\beta$.\footnote{Such  $\nu_\alpha$
exists globally if the fibration $\pi\colon M\to\tilde M$ is trivial. In general
the first Chern class provides an obstruction. However $\nu_\alpha$ always exists
in a neighbourhood of a periodic orbit.}
Then if we replace the Lagrangian by gauge-equivalent
$$
\,\widehat{\!\mathscr L}(x,\dot x,t)=\mathscr{L}(x,\dot x,t)-
c^\alpha \nu_\alpha(\dot x),
$$
Lagrange's equations do not change, but the Noether integrals will be replaced
by
$\,\,\widehat{\!\!\mathscr J}^\alpha\,{=}\,\mathscr{J}^\alpha-c^\alpha$.
Hence the value $c$ of the Noether integral is replaced by 0
and so the Routh function $\,\widetilde{\!\mathscr L}$
can be defined by \eqref{eq6.2}. The following theorem folds (see  \cite{32}).

\begin{theorem}[{\rm Routh}]
\label{th6.1} Let   $x(t)$ be a trajectory of the system $(M,\mathscr{L})$
with  $\mathscr{J}=c$. Then  $\tilde x(t)=\pi(x(t))$ is a trajectory of
the system $(\widetilde M,\,\widetilde{\!\mathscr L})$. Conversely, if
$\tilde x(t)$ is a  trajectory of the system $(\widetilde
M,\,\widetilde{\!\mathscr L})$, then there exists a trajectory $x(t)$
of the system
$(M,\mathscr{L})$ with  $\mathscr{J}=c$ such that  $\widetilde
x(t)=\pi(x(t))$.
\end{theorem}

If   $\gamma$~ is a periodic orbit of the system $(M,\mathscr{L})$ and
 $\widetilde\gamma$ the corresponding orbit of the system
$(\widetilde M,\,\widetilde{\!\mathscr L})$,
then their variational systems are related by a
linear version of Routh's method. In the next section we describe the
Routh reduction for a  linear Lagrangian system. It applies in a
more general case, for example, when the Lagrangian system has
non-Noether integrals. In particular, the linear Routh reduction
includes the linearized  Maupertuis reduction on an energy level.

\subsection{Routh reduction in a linear Lagrangian system}
\label{ssec6.1}
If the linear Poincar\'e map $P\colon W\to W$ of the periodic orbit $\gamma$
has  eigenvalue 1, then to any eigenvector $w=Pw$ there corresponds a
non-zero $\tau$-periodic solution $\zeta(t)$  of the variational
system $(E,\Lambda)$. As proved by Poincar\'e, the variational
system has a linear $\tau$-periodic first integral
$$
I_\zeta(\xi,D\xi)=(\zeta,D\xi)-(\xi,D\zeta).
$$
Indeed, by \eqref{eq5.4}
$$
\frac{d}{dt} I_\zeta\bigl(\xi(t),D\xi(t)\bigr)=(\zeta,D^2\xi)-(\xi,D^2\zeta)=
(\zeta,U\xi)-(\xi,U\zeta)=0.
$$
In fact, $I_\zeta(\xi,D\xi)=\omega(w,v)=J_w(v)$ is the
value of the symplectic form on the vectors $v,w\in W$
corresponding to $\xi$,~$\zeta$.

Suppose the Poincar\'e map $P$ has several eigenvectors corresponding to
unit eigenvalue. Let $V\subset \operatorname{Ker}(P-I)$ be an isotropic
subspace and  let $\Gamma\subset X$ be the corresponding vector space of periodic
solutions of the variational system $(E,\Lambda)$. Let $w^1,\dots,w^k$ be
a basis in $V$ and $\zeta^1,\dots,\zeta^k\in\Gamma$ the corresponding independent
solutions. The variational system has first integrals
$$
I^\alpha (\xi,D\xi)=(\zeta^\alpha,D\xi)-(\xi,D\zeta^\alpha),\qquad
\alpha=1,\dots, k,
$$
in involution
\begin{equation}
\label{eq6.3}
I^\alpha(\zeta^\beta,D\zeta^\beta)=(\zeta^\alpha,D\zeta^\beta)-
(\zeta^\beta,D\zeta^\alpha)=\omega(w^\alpha,w^\beta)=0.
\end{equation}
We write shortly  $I=(I^1,\dots,I^k)$.

Denote
$$
F_t=\{\zeta(t):\zeta\in\Gamma\}=\operatorname{span}\{\zeta^1(t),\dots,\zeta^k(t)\}.
$$
To simplify the presentation we use the following non-degeneracy assumption.

\begin{ncondition}
\label{nconA}
$\dim F_t=k$ {\it for all  $t$.}
\end{ncondition}

Equivalently,         $\zeta^1(t),\dots,\zeta^k(t)\in E_t$ are independent
for all $t$. Thus the Gram matrix
\begin{equation}
\label{eq6.4}
G=(g^{\alpha\beta}),\qquad
g^{\alpha\beta}(t)=\bigl(\zeta^\alpha(t),\zeta^\beta(t)\bigr),
\end{equation} is
nondegenerate for all $t$.

In Appendix \ref{ssecA.3} we will show that this assumption is unnecessary.
In fact, the set $\Sigma=\{t\in\mathbb{R}/\tau\mathbb{Z}:\dim F_t<k\}$ is finite and the
family $(F_t)_{t\notin\Sigma}$ can be  extended to a smooth $k$-dimensional vector
 bundle $(F_t)_{t\in\mathbb{R}/\mathbb{Z}}$. We will  show that   everything
in this section works without the non-degeneracy assumption~\ref{nconA}.

We describe Routh reduction for the linear system $(E,\Lambda)$.
The reduced configuration spaces  $\widetilde E_t=E_t/F_t$ can be
identified with
$$
E_t^\perp=\{u\in E_t:(u,w)=0 \ \text{for all} \ w\in F_t\}
$$
via the orthogonal projection    $\Pi=\Pi_t\colon E_t\to E_t^\perp$.
For a smooth  field      $\xi(t)\in E_t$ denote    $D^\perp\xi(t)=\Pi_t D\xi(t)$.
Explicitly,
\begin{equation}
\label{eq6.5}
\Pi\xi=\xi-g_{\alpha\beta}(\xi,\zeta^\beta)\zeta^\alpha, \qquad
D^\perp\xi=D\xi-g_{\alpha\beta}(D\xi,\zeta^\beta)\zeta^\alpha,
\end{equation} where
    $G^{-1}=(g_{\alpha\beta})$ is the inverse of the Gram matrix
$G=(g^{\alpha\beta})$. In Appendix  \ref{ssecA.3} we show that $\Pi$
and  $D^\perp$ are smooth also when the non-degeneracy assumption fails.

Define the Routh Lagrangian $\Lambda^\perp$ on
   $E^\perp=(E_t^\perp)$ by
$$
\Lambda^\perp(\eta,D^\perp\eta)=\frac{1}{2}\,(D^\perp\eta,D^\perp\eta)+
\frac{1}{2}\,(U^\perp\eta,\eta),
\qquad
\eta(t),D^\perp\eta(t)\in E_t^\perp,
$$
where the symmetric operator $U^\perp(t)\colon E_t^\perp\to E_t^\perp$
is given by
$$
(U^\perp u,u)=(Uu,u)-3 g_{\alpha\beta}(u,D^\perp\zeta^\alpha)
(u,D^\perp\zeta^\beta),\qquad
u\in E_t^\perp.
$$
Thus,  $U^\perp=\Pi U-3C$, where
$$
Cu=g_{\alpha\beta}(u,D^\perp\zeta^\alpha)D^\perp\zeta^\beta
$$
($C$ is independent on the choice of the basis).

The bilinear action form of the system $(E^\perp,\Lambda^\perp)$
is
\begin{equation}
\label{eq6.6}
\frac{1}{2}\,h^\perp(\eta,\eta)=
\int_0^\tau\Lambda^\perp(\eta,D^\perp\eta)\,dt,\qquad
\eta(t)\in E_t^\perp.
\end{equation}

We have Routh's theorem for linear Lagrangian systems.

\begin{theorem}
\label{th6.2}
Let  $\xi(t)\in E_t$ be a solution of the system $(E,\Lambda)$ such that
$I(\xi,D\xi)\equiv 0$. Then $\eta(t)=\Pi\xi(t)\in E_t^\perp$ is a
solution of  the system $(E^\perp,\Lambda^\perp)$. Conversely, if
$\eta(t)$ is a solution of the system $(E^\perp,\Lambda^\perp)$, then
there exists a solution $\xi(t)$ of  the system $(E,\Lambda)$, defined
$\mod  \Gamma$, such that $I(\xi,D\xi)=0$ and $\eta(t)=\Pi\xi(t)$.
\end{theorem}

For the proof we need the following evident result.

\begin{lemma}
\label{lem6.1}
Let
\begin{equation}
\label{eq6.7}
\eta(t)-\xi(t)\in F_t,\qquad
\xi(t)=\eta(t)+\lambda_\alpha(t)\zeta^\alpha(t).
\end{equation}
Then  $\xi$ satisfies $I_i^\alpha(\xi,D\xi)=c^\alpha$
for all  $\alpha=1,\dots,k$ if and only if
\begin{equation}
\label{eq6.8}
\dot\lambda_\alpha=g_{\alpha\beta}\bigl(c^\beta-I^\beta(\eta,D\eta)\bigr).
\end{equation}
\end{lemma}

Indeed,
$I^\alpha(\xi,D\xi)=I^\alpha(\eta,D\eta)+g^{\alpha\beta}\dot\lambda_\beta$.

\medskip\noindent{\it Proof of Theorem  \ref{th6.2}.}
A vector field $\xi(t)$, \,$0\le t\le \tau$, is a solution of
$(E,\Lambda)$ if and only if
$$
h(\xi,\phi)=\int_0^\tau \bigl((D\xi,D\phi)+(U\xi,\phi)\bigr)\,dt=0
$$
for any smooth variation           $\phi(t)\in E_t$ such that
$\phi(0)=\phi(\tau)=0$.

Suppose $I(\xi,D\xi)=0$ and let    $\eta=\Pi\xi$.
We need to show that for every smooth variation   $\phi(t)\in E_t^\perp$
such that  $\phi(0)=\phi(\tau)=0$ we have
$$
h^\perp(\eta,\phi)=\int_0^\tau \bigl((D^\perp\eta,D^\perp\phi)+
(U^\perp\eta,\phi)\bigr)\,dt=0,
$$
where  $h^\perp$ is the bilinear form \eqref{eq6.6}
corresponding to the Routh system.

By Lemma \ref{lem6.1}, $\xi=\eta+\lambda_\alpha\zeta^\alpha$, where
   $\eta(t)\in E_t^\perp$ and
\begin{equation}
\label{eq6.9}
\dot\lambda_\alpha=-g_{\alpha\beta}I^\beta(\eta,D\eta)
=2g_{\alpha\beta}(\eta,D\zeta^\beta).
\end{equation}
Since  $(\eta,\zeta^\alpha)=(\phi,\zeta^\alpha)=0$, \eqref{eq6.5}
gives
$$
D\eta=D^\perp\eta-g_{\alpha\beta}(\eta,D\zeta^\alpha)\zeta^\beta,\qquad
D\phi=D^\perp\phi-g_{\alpha\beta}(\phi,D\zeta^\alpha)\zeta^\beta.
$$
We obtain
\begin{align*}
h(\xi,\phi)&=\int_0^\tau\bigl((D\eta+\dot\lambda_\alpha \zeta^\alpha+
\lambda_\alpha D\zeta^\alpha,D\phi)+(U\eta,\phi)+
(\lambda_\alpha U\zeta^\alpha,\phi)\bigr)\,dt
\\
&=\int_0^\tau\bigl((D\eta,D\phi)+(U\eta,\phi)+
(\dot\lambda_\alpha\zeta^\alpha+\lambda_\alpha D\zeta^\alpha,D\phi)+
(\lambda_\alpha D^2\zeta^\alpha,\phi)\bigr)\,dt
\\
&=\int_0^\tau\biggl((D^\perp\eta,D^\perp\phi)+(U\eta,\phi)
+\bigl(g_{\alpha\beta}(\eta,D\zeta^\alpha)\zeta^\beta,
g_{\delta\varepsilon}(\phi,D\zeta^\delta)\zeta^\varepsilon\bigr)
\\
&\qquad+\frac{d}{dt}\bigl((\lambda_\alpha D\zeta^\alpha,\phi)
-2(\dot\lambda_\alpha \zeta^\alpha,\phi)\bigr)\biggr)\,dt
\\
&=\int_0^\tau\bigl((D^\perp\eta,D^\perp\phi)+(U\eta,\phi)
-3 g_{\alpha\beta}(\phi,D\zeta^\alpha)(\eta,D\zeta^\beta)\bigr)\,dt
\\
&=h^\perp(\eta,\phi)=0.
\end{align*}
Hence $\eta$ is a trajectory of $(E^\perp,\Lambda^\perp)$. We skip the
proof of the converse.
\qed \medskip

Now we can write Hill's formula for the reduced linear Poincar\'e map
$\widetilde P\colon \widetilde W\to\widetilde W$. Let
$$
X^\perp=\{\eta\in X:\eta(t)\in E_t^\perp\}
$$
and let
$$
H^\perp=(- D^{\perp 2}+I)^{-1}(-D^{\perp 2}+U^\perp)\colon X^\perp\to X^\perp
$$
be the Hessian operator for the reduced system
$(E^\perp,\Lambda^\perp)$. Let $Q^\perp\colon E_0^\perp\to E_0^\perp$ be
the operator of parallel transport corresponding to $D^\perp\eta=0$
and let $\sigma^\perp=\det Q^\perp=\pm 1$. If  assumption  \ref{nconA}
holds, then $\sigma=\sigma^\perp$ because the bundle $F$ is oriented.
In general                             $F$ can  be non-oriented.

\begin{corollary}
\label{cor6.1}
\begin{equation}
\label{eq6.10}
\det(\widetilde P-I)=\sigma^\perp(-1)^{m-k}\beta^\perp \det H^\perp,\qquad
\beta^\perp=e^{(m-k)\tau}{\det}^{-2}(e^\tau I-Q^\perp)>0.
\end{equation}
\end{corollary}

For the geodesic problem    $Q^\perp=Q\big|_{X^\perp}$,
\,$H^\perp=H\big|_{X^\perp}$, \,$\widetilde P=P^\perp$,
\,$\sigma=\sigma^\perp$,   and
we obtain Corollary  \ref{cor5.1}.

Note that in general $H^\perp\ne H\big|_{X^\perp}$, except when
$D^\perp\zeta=0$.
 The reason is that if $\eta\in X^\perp$ is $\tau$-periodic, $\lambda$ in
\eqref{eq6.9} is not periodic in general, and so $\xi\notin X$. Hence the space
$X^\perp$ of periodic
$\eta(t)\in E_t^\perp$ does not correspond to the space of periodic $\xi(t)\in E_t$
such that $I(\xi,D\xi)=0$. Thus $h^\perp$ is not the restriction of $h$ to
$X^\perp$ as in the geodesic case. Hence we need to discuss the relation between
$h^\perp$ and $h$.

\subsection{Elimination of degeneracy in the action functional}
\label{ssec6.2}
As in \eqref{eq3.21}, define     subspaces $Y,Z\subset X$ by
\begin{align}
\label{eq6.11}
Y&=\{\xi\in X: I(\xi,D\xi)\equiv \mathrm{const}\},
\\
\label{eq6.12}
Z&=\{\xi\in X:\xi(t)\in F_t \ \text{for all~$t$}\}.
\end{align}
If   $\zeta^1,\dots,\zeta^k \in \Gamma$ are basis periodic solutions, then
$$
Z=\biggl\{\xi(t)=\lambda_\alpha(t)\zeta^\alpha(t):
\lambda_\alpha(t+\tau)=\lambda_\alpha(t), \
\int_0^\tau g^{\alpha\beta}\dot\lambda_\alpha \dot\lambda_\beta\,dt
<\infty\biggr\}.
$$

\begin{lemma}
\label{lem6.2}
For any $\eta\in X$ there exists $\xi=\Phi(\eta)\in Y$, unique $\mod  \Gamma$,
  such that $\xi-\eta\in Z$. Explicitly, $\eta=\xi-\lambda_\alpha\zeta^\alpha$,
  where $\lambda_\alpha$ satisfies  \eqref{eq6.8} with
$c^\alpha=c^\alpha(\eta)$ given by
\begin{equation}
\label{eq6.13}
c^\alpha=\kappa^{\alpha\beta}\int_0^\tau g_{\beta\delta}
I^\delta(\eta,D\eta)\,dt,\qquad
(\kappa^{\alpha\beta})=(\bar g_{\alpha\beta})^{-1},\quad
\bar g_{\alpha\beta}=\int_0^\tau g_{\alpha\beta}\,dt.
\end{equation}
\end{lemma}

This follows from  Lemma \ref{lem6.1} for periodic $\xi$ and $\eta$.
We have defined  a projection $\Phi\colon X\to\widehat Y=Y/\Gamma$ which is
identical on $Y$
and $\Phi =0$ on $Z$. We obtain

\begin{proposition}
\label{pr6.1}
$Y\cap Z=\Gamma$  and $Y+Z=X$.
\end{proposition}

\begin{proposition}
\label{pr6.2}
The spaces  $Y$ and  $Z$ are $h$-orthogonal. that is,  $h(\xi,\eta)=0$
for all  $\xi\in Z$ and $\eta\in Y$. The restriction of $h$ to $Z$ has the
form
\begin{equation}
\label{eq6.14}
h(\lambda_\alpha\zeta^\alpha,\lambda_\beta\zeta^\beta)
=\int_0^\tau g^{\alpha\beta}\dot\lambda_\alpha \dot\lambda_\beta\,dt.
\end{equation}
\end{proposition}

\noindent{\it Proof.}
Take    $\xi\in Z$, \,$\xi(t)=\lambda_\alpha(t)\zeta^\alpha(t)$. Then
\begin{align}
h(\lambda_\alpha\zeta^\alpha,\eta)&=\int_0^\tau
\Bigl(\bigl(D(\lambda_\alpha\zeta^\alpha),D\eta\bigr)
+(U\lambda_\alpha\zeta^\alpha,\eta)\Bigr)\,dt
\nonumber
\\
&=\int_0^\tau\bigl((\dot\lambda_\alpha\zeta^\alpha,D\eta)
+(\lambda_\alpha D\zeta^\alpha,D\eta)
+(\lambda_\alpha D^2\zeta^\alpha,\eta)\bigr)\,dt
\nonumber
\\
&=\int_0^\tau\biggl(\frac{d}{dt}(\eta,\lambda_\alpha D\zeta^\alpha)
+\dot\lambda_\alpha I^\alpha(\eta,D\eta)\biggr)\,dt
\nonumber
\\
&=(\eta,\lambda_\alpha D\zeta^\alpha)\big|_0^\tau
+\int_0^\tau \dot\lambda_\alpha I^\alpha(\eta,D\eta)\,dt
\label{eq6.15}
\end{align}
(we have used that  $\zeta$ satisfies the variational system).
If $\eta=\lambda_\alpha\zeta^\alpha$, then
$I^\alpha(\eta,D\eta)=g^{\alpha\beta}\dot\lambda_\beta$. Hence
\begin{equation}
h(\lambda_\alpha\zeta^\alpha,\lambda_\beta\zeta^\beta)
=\int_0^\tau g^{\alpha\beta}\dot\lambda_\alpha\dot\lambda_\beta\,dt
+(\lambda_\alpha\zeta^\alpha,\lambda_\beta D\zeta^\beta)\big|_0^\tau.
\label{eq6.16}
\end{equation}
If   $\eta\in Y$ and the~$\lambda_\alpha$ are periodic, \eqref{eq6.15} gives
0 and  \eqref{eq6.16} gives \eqref{eq6.14}, which proves Proposition~\ref{pr6.2}.
\qed \medskip

Let
$$
\widehat X=X/\Gamma,\qquad
\widehat Y=Y/\Gamma,\qquad
\widehat Z=Z/\Gamma.
$$
Then  $\widehat Y\oplus\widehat Z=\widehat X$. The bilinear form~$h$
defined a form  $\hat h$ on~$\widehat X$ and
$\widehat Y\perp_{\hat h} \widehat Z$, while  $\hat h\big|_{\widehat Z}$
is positive definite.

\begin{corollary}
\label{cor6.2}
$$
\operatorname{ind} h=\operatorname{ind} \hat h\big|_{\widehat Y},\qquad
\operatorname{null} h=\operatorname{null}\hat h\big|_{\widehat Y}+k.
$$
\end{corollary}

\begin{corollary}
\label{cor6.3} The projection
$\Pi\colon X\to X^\perp$ defines an isomorphism
\begin{equation}
\label{eq6.17}
\widehat\Pi\big|_{\widehat Y}\colon\widehat Y\to X^\perp,\qquad
\bigl(\Pi\big|_{\widehat Y}\bigr)^{-1}=\Phi\big|_{X^\perp}.
\end{equation}
\end{corollary}

Indeed, if $\Pi\xi=0$ for $\xi\in Y$, then  $\xi\in Z$ and
hence   $\xi\in Y\cap Z=\Gamma$ by Proposition~\ref{pr6.1}. Similarly,
$Y+Z=X$ implies that $\Pi(Y)=X^\perp$.

Next we compute the restriction $h\big|_Y$. Let   $h^\perp(\eta,\eta) $
be the bilinear form for the reduced system  $(E^\perp,\Lambda^\perp)$.

\begin{proposition}
\label{pr6.3}
Let   $h^\top=h\big|_Y\circ \widehat\Pi^{-1}$ be the bilinear form
on $X^\perp$ corresponding to  $h\big|_Y$.
Then for any      $\eta\in X^\perp$,
\begin{equation}
\label{eq6.18}
h^\top(\eta,\eta)=h^\perp(\eta,\eta)+
\bar g_{\alpha\beta} c^\alpha c^\beta,
\end{equation}
where the $c^\alpha(\eta)=I^\alpha(\xi,D\xi)$, \,$\xi=\Phi\eta$,
are defined by  \eqref{eq6.13}.
\end{proposition}

This follows from a more general formula.

\begin{lemma}
\label{lem6.3}
Let $\xi(t)=\eta(t)+\lambda_\alpha(t)\zeta^\alpha(t)$, \,$0\le t\le \tau$,
be as in Lemma \ref{lem6.1}. Then
\begin{align}
\nonumber
h(\xi,\bar\xi\,)&=h(\eta,\eta)-
\int_0^\tau g_{\alpha\beta}(I^\alpha(\eta,D\eta)-c^\alpha)
\bigl(I^\beta(\eta,D\eta)-c^\beta\bigr)\,dt
\\
&\qquad+\bigl((\lambda_\alpha\zeta^\alpha,\lambda_\beta D\zeta^\beta)
+2(\lambda_\alpha\eta,D\zeta^\alpha)+
2c^\alpha\lambda_\alpha\bigr)\big|_0^\tau.
\label{eq6.19}
\end{align}
\end{lemma}

\noindent{\it Proof.}
We have
$$
h(\xi,\xi)=h(\eta,\eta)+2h(\eta,\lambda_\alpha\zeta^\alpha)
+h(\lambda_\alpha\zeta^\alpha,\lambda_\beta\zeta^\beta).
$$
Now    \eqref{eq6.19} follows from \eqref{eq6.15},~\eqref{eq6.16}, and
  \eqref{eq6.8}.
\qed \medskip

\noindent{\it Proof of Proposition~\ref{pr6.3}.}
If  $\eta$ and  $\lambda$ are periodic, then the boundary terms in \eqref{eq6.19}
vanish. By  \eqref{eq6.13},
\begin{align*}
&\int_0^\tau g_{\alpha\beta}\bigl(c^\alpha-I^\alpha(\eta,D\eta)\bigr)
\bigl(c^\beta-I^\beta(\eta,D\eta)\bigr)\,dt
\\
&\qquad=\int_0^\tau g_{\alpha\beta}I^\alpha(\eta,D\eta)I^\beta(\eta,D\eta)\,dt
-2\int_0^\tau g_{\alpha\beta} c^\alpha I^\beta(\eta,D\eta)\,dt
+\bar g_{\alpha\beta} c^\alpha c^\beta
\\
&\qquad=\int_0^\tau g_{\alpha\beta}I^\alpha(\eta,D\eta)I^\beta(\eta,D\eta)\,dt
-\bar g_{\alpha\beta}c^\alpha c^\beta.
\end{align*}
Next we use $\eta=\Pi\xi\in X^\perp$. Then
$I^\alpha(\eta,D\eta)=-2(\eta,D^\perp\zeta^\alpha)$, and so
$$
\int_0^\tau g_{\alpha\beta}I^\alpha(\eta,D\eta)I^\beta(\eta,D\eta)\,dt
=4\int_0^\tau g_{\alpha\beta}(\eta,D^\perp\zeta^\alpha)
(\eta,D^\perp\zeta^\beta)\,dt.
$$
Finally,
\begin{align}
h(\eta,\eta)&=\int_0^\tau\bigl((D^\perp\eta,D^\perp\bar\eta)
+g_{\alpha\beta}(\eta,D^\perp\zeta^\alpha)(\eta,D^\perp\zeta^\beta)
+(U\eta,\bar\eta)\bigr)\,dt
\nonumber
\\
&=h^\perp(\eta,\eta)+4\int_0^\tau g_{\alpha\beta}(\eta,D^\perp\zeta^\alpha)
(\eta,D^\perp\zeta^\beta)\,dt.
\label{eq6.20}
\end{align}
It remains to substitute \eqref{eq6.20} in \eqref{eq6.19}.
\qed \medskip

From the point of view of Routh reduction it is natural to consider the space
\begin{equation}
\label{eq6.21}
Y^0=\{\xi\in X:I(\xi,D\xi)\equiv 0\}\subset Y.
\end{equation}
Indeed, by \eqref{eq6.18},
$h\big|_{Y^0}=h^\perp\circ \Pi\big|_{Y^0}$.

\begin{proposition}
\label{pr6.4}
\begin{equation}
\label{eq6.22}
Y^0+Z=\biggl\{\eta\in X:\int_0^\tau g_{\alpha\beta}
(D^\perp\zeta^\beta,\eta)\,dt=0\biggr\}.
\end{equation}
\end{proposition}

\noindent{\it Proof.}
We take  $\eta\in X$. Then  $\xi=\Phi\eta\in Y^0$ provided  that
$$
0=\int_0^\tau g_{\alpha\beta}I^\beta(\eta,D\eta)\,dt
=-2 \int_0^\tau g_{\alpha\beta}(D^\perp\zeta^\beta,\eta)\,dt.
$$
Here we have used that
\begin{equation}
\label{eq6.23}
g_{\alpha\beta}I^\beta(\eta,D\eta)
=\frac{d}{dt}\bigl(g_{\alpha\beta}(\eta,\zeta^\beta)\bigr)
-2 g_{\alpha\beta}(\eta,D^\perp\zeta^\beta).
\end{equation}
\qed \medskip

We see that    $Y^0+Z=X$ if and only if $D^\perp\zeta^\alpha=0$.
Equivalently, $DZ\subset Z$. Let
\begin{equation}
\label{eq6.24}
X_0^\perp=\Pi Y^0=\biggl\{\eta\in X^\perp:\int_0^\tau g_{\alpha\beta}
(D^\perp\zeta^\beta,\eta)\,dt=0\biggr\}.
\end{equation}
Then $X_0^\perp$ has codimension  $\le k$ in~$X^\perp$.
We have $h^\top\ge h^\perp$ and~$h^\top=h^\perp$ on~$X_0^\perp$. Since
$\operatorname{ind} h^\top=\operatorname{ind} h$,
$$
\operatorname{ind} h^\perp\le \operatorname{ind} h\big|_X\le \operatorname{ind} h^\perp+k.
$$

Let   $\Omega=Y/Y^0$. Since   $\widehat\Pi\colon Y\to X^\perp$ is an isomorphism,
   $\dim\Omega\le k$. The integral $I\colon Y\to\mathbb{R}^k$ gives
a map       $\Omega\to \mathbb{R}^k$. To compare the indices of $h^\top$
and
  $h^\perp$, in \S\,\ref{ssec6.4} we construct a basis in $\Omega$,
on which  $I$ is nondegenerate.

\subsection
{Indices of~$h$ and~$h^\perp$}
\label{ssec6.3}
In this section we discuss the relation between $\operatorname{ind}
h=\operatorname{ind} h\big|_Y$ and   $\operatorname{ind} h^\perp$. Let
$P\colon W\to W$ be the Poincar\'e map. As in \S\,\ref{ssec3.7},
we assume that
$$
N=\operatorname{Ker}(P-I)^2=\{v\in W:Pv-v\in V\}
$$
is symplectic and                 $V=\operatorname{Ker}(P-I)$ is a Lagrangian
subspace in      $N$. Let   $w^1,\dots,w^k$ be a basis in  $V$ and
  $q_1,\dots,q_k$ a basis in a Lagrangian complement to~$V$ in  $N$.
We define the
matrix  $s_{\alpha\beta}$ by formula \eqref{eq3.30} and the matrices
$A=(a_{\alpha\beta})$ and~$A^\perp=(a_{\alpha\beta}^\perp)$
by \eqref{eq3.31}, where       $\varkappa^{\alpha\beta}$ and  $\bar
g_{\alpha\beta}$ are the matrices in \eqref{eq6.13}.   

\begin{definition}
We say that $\gamma$ is nondegenerate $\mod V$ if Assumption
\ref{nconA} (p.~\pageref{nconA}) holds and  $\det A^\perp\ne 0$.
\end{definition}

\begin{theorem}
\label{th6.3}
Suppose that $\gamma$ is nondegenerate  $\mod V$. Then
\begin{equation}
\label{eq6.25}
\operatorname{ind} h-\operatorname{ind} h^\perp=\operatorname{ind} A-\operatorname{ind} A^\perp.
\end{equation}
\end{theorem}

The formulation coincides with   Theorem \ref{th3.2},
but the proof is different.
We prove Theorem \ref{th6.3} in \S\,\ref{ssec6.4}.
Since $h\big|_{\widehat Z}$ is positive definite, as in the proof of
Corollary \ref{cor3.3},
 we obtain

\begin{corollary}
\label{cor6.4}
Suppose $\det A\ne 0$. Then
$$
(-1)^{\operatorname{ind} h}=(-1)^{\operatorname{ind} h^\perp+\operatorname{ind} b}.
$$
\end{corollary}

\subsection
{The spaces $\Omega$ and~$\Omega^\perp$}
\label{ssec6.4}
Consider periodic solutions $\zeta^\alpha(t)$, $\eta_\alpha(t)$ of the
system  $(E,\Lambda)$ which correspond to  $w^\alpha$,~$q_\alpha$.
Then the  $\zeta^\alpha(t)$ are periodic and satisfy \eqref{eq6.3}.
Equations \eqref{eq3.30} imply
\begin{equation}
\label{eq6.26}
\eta_\alpha(t+\tau)-\eta_\alpha(t)=s_{\alpha\beta}\zeta^\beta(t),\qquad
(\eta_\alpha,D\eta_\beta)-(\eta_\beta,D\eta_\alpha)=0
\end{equation}
and  $(\zeta^\alpha,D\eta_\beta)-(\eta_\beta,D\zeta^\alpha)=\delta_\beta^\alpha$.

For any $\alpha=1,\dots,k$ we put
$$
\hat\eta_\alpha=\eta_\alpha-\lambda_{\alpha\beta}\zeta^\beta,
$$ where
the coefficients $\lambda_{\alpha\beta}$ are chosen so that
the $\hat\eta_\alpha$ are $\tau$-periodic and
$$
(\zeta^\alpha,D\hat\eta_\beta)-(\hat\eta_\beta,D\zeta^\alpha)
=c^\alpha_\beta=\mathrm{const}.
$$
Then the~$c_\beta^\alpha$ satisfy  \eqref{eq3.38} and
$$
\lambda_{\alpha\beta}(t+\tau)-\lambda_{\alpha\beta}(t)=s_{\alpha\beta},\qquad
\dot\lambda_{\alpha\beta}=s_{\alpha\delta}\kappa^{\delta\varepsilon}
g_{\varepsilon\beta}.
$$

We define  $\eta_\alpha^\perp=\Pi\eta_\alpha$. Then
$$
\eta_\alpha^\perp=\eta_\alpha-\lambda_{\alpha\beta}^\perp \zeta^\beta,\qquad
\lambda_{\alpha\beta}^\perp(t+\tau)-\lambda_{\alpha\beta}^\perp(t)
=s_{\alpha\beta},\qquad
\lambda_{\alpha\beta}^\perp=(\eta_\alpha,\zeta^\delta) g_{\delta\beta}.
$$
Consider the spaces
$$
\Omega=\operatorname{span}(\hat\eta_1,\ldots,\hat\eta_k)=
\Phi\Omega^\perp\subset\widehat Y,
\qquad
\Omega^\perp=\operatorname{span}(\eta^\perp_1,\dots,\eta^\perp_k)
=\Pi\Sigma\subset X^\perp.
$$
We also define $\widehat\Omega$ and  $\widehat Y^0$
as the images of  $\Omega$ and $Y^0$ under the canonical projection
$\Pi_\Gamma\colon X\to \widehat X=X/\Gamma$.

\begin{theorem}
\label{th6.4}
Suppose $\det A^\perp \ne 0$. Then

\begin{numa}

\item
\label{th6.4a}
$\widehat Y=\widehat\Omega\oplus\widehat Y^0$,
\,$X^\perp=\Omega^\perp\oplus \Pi Y^0$;

\item
\label{th6.4b}
the maps    $\widehat\Pi\big|_{\widehat\Omega}\colon \widehat\Omega
\to \Omega^\perp$ and  $\widehat\Pi\big|_{\widehat Y^0}\colon \widehat Y^0\to
X_0^\perp$ are linear isomorphisms;

\item
\label{th6.4c}
$h\big|_{Y^0}=h^\perp\big|_{X_0^\perp} \circ\Pi\big|_{Y^0}$ and
$\hat h\big|_{\widehat Y^0}=h^\perp\big|_{X_0^\perp} \circ
\widehat\Pi\big|_{\widehat Y^0}$;

\item
\label{th6.4d}
for any  $\xi\in Y^0$ and~$\alpha,\beta=1,\dots,k$,
$$
h(\hat\eta_\alpha,\xi)=0, \qquad
h(\hat\eta_\alpha,\hat\eta_\beta)=a_{\alpha\beta};
$$

\item
\label{th6.4e}
for any  $\xi^\perp\in\Pi Y^0$ and~$\alpha,\beta=1,\dots,k$
$$
h^\perp(\eta_\alpha^\perp,\xi^\perp)=0, \qquad
h^\perp(\eta_\alpha^\perp,\eta_\beta^\perp)=a_{\alpha\beta}^\perp.
$$
\end{numa}
\end{theorem}

The proof of Theorem  \ref{th6.4} is contained in \S\,\ref{ssecA.2}.

\begin{corollary}
\label{cor6.5}
In the basis  $\hat\eta_1,\dots,\hat\eta_k$
$$
\hat h\big|_{\widehat\Omega}-h^\perp\circ\Pi\big|_{\widehat\Omega}
=(s_{\alpha\delta}\kappa^{\delta\varepsilon} s_{\varepsilon\beta}
-2s_{\alpha\beta}+\bar g_{\alpha\beta})=SKS-2S+\overline{G}.
$$
This quadratic form is positive definite.
\end{corollary}

Indeed, let $Q$ be the square root of $K$, that is,
the positive definite symmetric matrix such that $Q^2=K$.
Then   $Q^{-2} = \overline{G}$ and
$$
SKS-2S+\overline{G}=RR^*,\qquad
R=(S-\overline{G}\,)Q.
$$
This matrix is positive definite because $R$ is nondegenerate.

Now we prove Theorem \ref{th6.3}. Recall that by Proposition \ref{pr6.2}
and Theorem \ref{th6.4} we have the $h$-orthogonal expansion $X = Z\oplus\Omega\oplus Y^0$
and the $h^\perp$-orthogonal expansion
$X^\perp = \Omega^\perp\oplus\Pi Y^0$. The form  $h|_Z$ is positive definite
 and
$h\big|_{Y^0}=h^\perp\big|_{{X_0^\perp}} \circ\Pi\big|_{Y^0}$. Therefore,
$$
\operatorname{ind} h-\operatorname{ind} h^\perp=\operatorname{ind} h\big|_\Omega-\operatorname{ind} h^\perp\big|_{\Omega^\perp}
=\operatorname{ind} A-\operatorname{ind} A^\perp.
$$

\subsection{Example: autonomous systems}
\label{ssec6.5}
Suppose the Lagrangian system is autonomous, so the variational system of a
periodic trajectory $\gamma$ has a
periodic solution $\zeta(t)=\dot\gamma(t)$. If $\gamma$ is nondegenerate in the
autonomous sense (only two unit multipliers)
then, as proved by Poincar\'e, there exists a
 family\footnote{It seems more natural to parametrize the family by the period $\tau$.
However, this is not always possible because it may happen that
$\tau'(\alpha)=0$.} $\gamma_\alpha$ of $\tau(\alpha)$-periodic orbits such
that $\gamma_0=\gamma$
and $\tau(0)=\tau$ (\cite{33},  \cite{34}).
Let
$E(\alpha)=\mathscr{H}\big|_{\gamma_\alpha}$ and
$A(\alpha)=\displaystyle\int_{\gamma_\alpha}\langle p,dx\rangle$~
be the energy
and Maupertuis action of $\gamma_\alpha$.

\begin{lemma}
\label{lem6.4}
Suppose that $dE / d\tau \ne 0$. Then
$$
(-1)^{\operatorname{ind} b}=-\operatorname{sign}\bigl(\tau'(\alpha)
E'(\alpha)\bigr)=-\operatorname{sign}\frac{dE}{d\tau}\,.
$$
\end{lemma}

\noindent{\it Proof.}
The union of trajectories of $\gamma_\alpha$ in the phase space $TM\cong T^*M$ is
a symplectic cylinder $\Sigma$.
Restricting the Hamiltonian system to $\Sigma$ we obtain an integrable Hamiltonian system
with one  degree of freedom and Hamiltonian
$H(\vartheta,I)=E(I)$, where
$\vartheta\in\mathbb{R}/\mathbb{Z}$, \,$I\in\mathbb{R}$.
Then we can assume      $\alpha=I$, \,$\gamma_\alpha(t)=(E'(\alpha)t,\alpha)$.
Then  $\tau(\alpha)=1/\nu(\alpha)$, where the frequency  $\nu(\alpha)$ is
$E'(\alpha)$. We have $w=\begin{pmatrix} 1\\ 0\end{pmatrix}$ and
$v=\begin{pmatrix} 0\\ 1\end{pmatrix}$. The monodromy matrix
of $\gamma_\alpha$ is
$P=P_\alpha=\begin{pmatrix} 1 & \nu'/\nu \\ 0 & 1 \end{pmatrix}$.
Thus    $Pv=\begin{pmatrix} \nu'/\nu \\ 1\end{pmatrix}$, and so
$s=\omega(Pv,v)=\nu'/\nu=-\tau^{-2}\,d\tau/dE$.
\qed \medskip

As usual, we denote   $\operatorname{ind}\gamma=\operatorname{ind} h(\gamma)$.

\begin{proposition}
\label{pr6.5}
Let a periodic trajectory $\gamma$ have exactly $2$ unit multipliers.
Suppose that
$\sigma (-1)^{m+\operatorname{ind}\gamma}dE/d\tau<0$. Then  $\gamma$ has
a real multiplier $\rho > 1$.
\end{proposition}

\noindent{\it Proof.}
Since   $\sigma=\sigma^\perp$, by Corollary \ref{cor6.4} and
Lemma \ref{lem6.4},
$$
(-1)^{\operatorname{ind} h^\perp}=-\operatorname{sign}\biggl(\frac{dE}{d\tau}\biggr)(-1)^{\operatorname{ind}\gamma}.
$$
The dimension of the reduced system is   $m^\perp=m-1$. Hence
$$
\sigma (-1)^{m^\perp+\operatorname{ind} h^\perp}=
\sigma(-1)^{m-1+\operatorname{ind}\gamma}\biggl(-\operatorname{sign}\frac{dE}{d\tau}\biggr)=-1,
$$ and by
Corollary  \ref{cor5.2} applied to the reduced Hill formula
      \eqref{eq6.10}, there exists a multiplier $\rho > 1$.
\qed \medskip

\begin{example}
\label{ex6.1}                                   {\rm Suppose
a particle in $\mathbb{R}^m$ moves under the potential field
with homogeneous potential energy
$$
V(\lambda x)=\lambda^kV(x), \qquad
\lambda>0, \quad
k(k-2) \ne 0.
$$
Suppose $\gamma$ is a $\tau$-periodic solution with energy $E$.
Then  $\gamma_\lambda(t)=\lambda\gamma(\lambda^{k/2-1}t)$ is a periodic
solution with period $\tau(\lambda)=\lambda^{1-k/2}\tau$ and
energy   $E(\lambda)=\lambda^k E$. Hence
\begin{equation}
\label{eq6.27}
\frac{dE(\lambda)}{d\tau(\lambda)}=
\frac{2k}{k-2}\biggl(\frac{\tau(\lambda)}{\tau}\biggr)^{(k+2)/(k-2)}E.
\end{equation}
Thus by Lemma   \ref{lem6.4}
$$
(-1)^{\operatorname{ind} b}=\operatorname{sign}\frac{2-k}{k}\,.
$$
}
\end{example}

Consider the problem of the motion of a particle in $\mathbb{R}^m$ in the
force field generated by a homogeneous potential of degree $k$, where
$k(k-2)\ne 0$. Equations $\sigma=1$,  \eqref{eq6.27}, and Proposition
\ref{pr6.5}  immediately imply

\begin{proposition}
\label{pr6.6}
Let a periodic trajectory $\gamma$ have exactly $2$ unit multipliers.
Suppose that $(-1)^{m + \operatorname{ind}\gamma} (k-2)/k < 0$.
Then $\gamma$ has a real multiplier $\rho > 1$.
\end{proposition}

\subsection
{Degeneracy in the $\rho$-index form}
\label{ssec6.6}
We have seen that the relation between $\operatorname{ind}  h$ and
$\operatorname{ind}  h^\perp$ is not evident.
This simplifies drastically for  $\rho\ne 1$. Let  $X_\rho$, $\rho\in S^1$, be the set of complex
$\rho$-quasiperiodic vector fields. Similarly to \eqref{eq6.11}, \eqref{eq6.12}
 define subspaces $Y_\rho,Z_\rho\subset X_\rho$:
$$
Y_\rho=\{\xi\in X_\rho: I(\xi,D\xi)\equiv \mathrm{const}\},\qquad
Z_\rho=\{\xi\in X_\rho: \xi(t)\in F_t\}.
$$
It is easy to see that for $\rho\ne 1$ and~$\xi\in X_\rho$,
$I(\xi,D\xi)\equiv c$ implies $c=0$. Thus,
$$
Y_\rho=Y_\rho^0=\{\xi\in X_\rho: I(\xi,D\xi)\equiv 0\}.
$$

\begin{proposition}
\label{pr6.7}
For $\rho\ne 1$ we have $Y_\rho\cap Z_\rho=\{0\}$ and
$X_\rho=Y_\rho\oplus Z_\rho$.
\end{proposition}

\noindent{\it Proof.}                                         We will define
a projection       $\Phi_\rho\colon X_\rho\to Y_\rho$ along $Z_\rho$.
Take    $\eta\in X_\rho$ and look for           $\lambda_\alpha(t)$ such
that $\xi=\eta+\lambda_\alpha\zeta^\alpha\in Y_\rho$.
Then by  \eqref{eq6.8}, $\dot\lambda_\alpha=f_\alpha$, where
$f_\alpha(t+\tau)=\rho f_\alpha(t)$.
Hence $f_\alpha(t)=e^{\mu t} b_\alpha(t)$, where $\mu=\tau^{-1}\ln\rho$
and       $b_\alpha(t)$ is a     $\tau$-periodic function:
$$
b_\alpha(t)=\sum_{k\in\mathbb{Z}} b_{\alpha k} e^{k\omega t},\qquad
\omega=\frac{2\pi i}{\tau}\,.
$$
We obtain a unique solution  $\lambda_\alpha(t)$ such that
$\lambda_\alpha(t+\tau)=\rho\lambda_\alpha(t)$:
\begin{equation}
\label{eq6.28}
\lambda_\alpha(t)=e^{\mu t}\sum_{k\in\mathbb{Z}}
\frac{b_{\alpha k} e^{k\omega t}}{\mu+k\omega}\,.
\end{equation}
The denominator is non-zero if $\rho\ne 1$.
\qed \medskip

\begin{proposition}
\label{pr6.8}
The spaces   $Y_\rho,Z_\rho$ are $h$-orthogonal:
$$
h(\xi,\bar\eta)=0\quad\text{for all} \
\xi\in Z_\rho,\
\eta\in Y_\rho.
$$ The
restriction of $h$ to~$Z_\rho$ is positive definite for    $\rho\ne 1$:
\begin{equation}
\label{eq6.29}
h(\lambda_\alpha\zeta^\alpha,\bar\lambda_\beta\zeta^\beta)=
\int_0^\tau g^{\alpha\beta}\dot\lambda_\alpha\dot{\bar\lambda}_\beta\,dt.
\end{equation}
\end{proposition}

\noindent{\it Proof.}
Take    $\xi(t)=\lambda_\alpha(t)\zeta^\alpha(t)\in Z_\rho$ and
$\eta\in Y_\rho$. Then by \eqref{eq6.15},
\begin{equation}
\label{eq6.30}
h(\lambda_\alpha\zeta^\alpha,\bar\eta)
=(\lambda_\alpha D\zeta^\alpha,\bar\eta)\big|_0^\tau
=(|\rho|^2-1)\bigl(\lambda_\alpha(0)D\zeta^\alpha(0),D\eta(0)\bigr)=0.
\end{equation}
The proof of  \eqref{eq6.29} is similar.
\qed \medskip

Next we compute the restriction of  $h$ to~$Y_\rho$. Let
$$
X_\rho^\perp=\{\eta\in X_\rho:\eta(t)\in E_t^\perp\},
$$
and let $\Pi\colon X_\rho\to X_\rho^\perp$ be the projection~\eqref{eq6.5}.
Since   $X_\rho=Y_\rho\oplus Z_\rho=X_\rho^\perp\oplus Z_\rho$,
   $\Pi\big|_{Y_\rho}\colon Y_\rho\to X_\rho^\perp$ is an isomorphism and
its inverse is           $\Phi_\rho\big|_{X_\rho^\perp}$.

\begin{proposition}
\label{pr6.9}
For $\rho \ne 1$ the bilinear form $h_\rho^\top=h\circ\bigl(\Pi\big|_{Y_\rho}\bigr)^{-1}$
on~$X_\rho^\perp$ is equal to the Routh form  $h_\rho^\perp$. 
\end{proposition}

This follows from Lemma~\ref{lem6.1} (for complex vector fields) since
$c=0$ and  $\lambda(\tau)=\rho\lambda(0)$, \,$\eta(\tau)=\rho\eta(0)$,
\,$|\rho|=1$.

\begin{corollary}
\label{cor6.6}
For  $\rho\ne 1$ the $\rho$-index of the system $(E,\lambda)$ equals the
$\rho$-index of the Routh system $(E^\perp,\Lambda^\perp)$.
\end{corollary}

Proposition \ref{pr6.9} is not true for $\rho=1$.
Then the relation between indices is more complicated, as we saw before.

\section{Reversible case}
\label{sec7}

Suppose the Lagrangian system $(M,\mathscr{L})$ is reversible: there is an
involution $S\colon M\to M$ which is a time reversing symmetry for
    $\mathscr{L}$:
$$
\mathscr{L}(S(x),dS(x)\dot x,t)=\mathscr{L}(x,-\dot x,-t).
$$
Let     $\tau=2T$. Then for any $\tau$-periodic curve $\gamma\in\Omega$,
$$
\mathscr{A}(\gamma)=\int_{-T}^{T}\mathscr{L}(\gamma(t),\dot\gamma(t),t)\,dt
=\mathscr{A}(\widetilde\gamma),
$$                             where
    $\widetilde\gamma(t)=S\gamma(-t)$. Thus, the involution
$R\colon\Omega\to\Omega$,
$$
R(\gamma)(t)=\widetilde\gamma(t)=S(\gamma(-t)),
$$
preserves $\mathscr{A}$. A
$\tau$-periodic orbit $\gamma$ is called reversible if
$R(\gamma)=\gamma$. Then
$$
S\gamma(-t)=\gamma(t),\qquad
S\gamma(T-t)=\gamma(T+t).
$$
Hence  $\gamma(0)$ and  $\gamma(T)$ belong to the set  $N$
of fixed points of $S$. It is easy to see that $\gamma$
is a reversible periodic orbit if and only if
$\gamma_+=\gamma\big|_{[0,T]}$ is a critical point of the action functional
$$
\mathscr{A}_+(\nu)=\int_0^T \mathscr{L}\bigl(\nu(t),\dot\nu(t),t\bigr)\,dt
$$
on the set   $\Omega_+$ of curves $\nu\colon[0,T]\to M$ with end-points in
  $N$.

\goodbreak
Let  $X=T_\gamma\Omega$ be the set of vector fields along  $\gamma$
and $J=dR(\gamma)\colon X\to X$. Then
$$
(J\xi)(t)=J_{-t}\xi(-t),\qquad
J_t=dS(\gamma(t))\colon E_t\to E_{-t}.
$$
Since  $R$ preserves $\mathscr{A}$, the involution $J$
preserves the Hessian bilinear form:
$$
h(J\xi,J\eta)=h(\xi,\eta).
$$
Since the operators~$D$ and~$U$ are intrinsically associated with  $h$,
$$
J^*=J,\qquad
DJ=-JD,\qquad
UJ=JU.
$$

Let     $X_\pm=\{\xi\in X:J\xi=\pm \xi\}$. Then  $X=X_+\oplus X_-$ and
any  $\xi\in X$ is represented as $\xi=\xi_++\xi_-$, where
$\xi_\pm\in X_\pm$. Then
$$
h(\xi,\xi)=h(\xi_+,\xi_+)+h(\xi_-,\xi_-).
$$

Since   $DJ=-JD$, we have $D^2\colon X_\pm \to X_\pm$.
Hence the Hessian operator
$H=(-D^2+I)^{-1}(-D^2+U)$ commutes with~$J$, and so  $H\colon X_\pm\to X_\pm$.
Denote    $H_\pm=H\big|_{X_\pm}$.

\begin{proposition}
\label{pr7.1}
$H=H_+\oplus H_-$ and  $\det H=\det H_+\det H_-$.
\end{proposition}

Next we give more explicit formulae for $h_\pm$. Any $\xi\in X_\pm$
is determined by the restriction
$$
\xi\big|_{[0,T]}\in Y_\pm=\{\eta\in Y: \eta(0)\in E_0^\pm, \
\eta(T)\in E_T^\pm\},\qquad
Y=W^{1,2}([0,T],E),
$$
where
$$
E_0^\pm=\{v\in E_0:J_0v=\pm v\}, \qquad
E_T^\pm=\{v\in E_T: J_T v=\pm v\}.
$$

Thus we have the orthogonal decompositions
$$
E_0=E_0^+\oplus E_0^-,\qquad
E_T=E_T^+\oplus E_T^-.
$$

For $\eta\in Y_\pm$ the corresponding    $\xi\in X_\pm$ is given by
$$
\xi\big|_{[0,T]}=\eta,\qquad
\xi\big|_{[-T,0]}=\pm J\eta .
$$
Thus,
$$
h(\xi,\xi)=h(\eta,\eta)+h(J\eta,J\eta)=2K(\eta,\eta),
$$
where
$$
K(\eta,\eta)=\int_0^T \bigl((D\eta,D\eta)+(U\eta,\eta)\bigr)\,dt
$$
is the same form $h$, but considered on     $Y$.
Let   $K_\pm=K\big|_{Y_\pm}$. Then  $K_+=d^2\mathscr{A}_+(\gamma_+)$
is the second variation of the functional~$\mathscr{A}_+$.

Let us consider the case $S=\mathrm{id}$. Then  $E_0^+=E_0$, \,$E_T^+=E_T$, \,$E_0^-=0$,
\,$E_T^-=0$. Thus,  $Y_+=Y$ and
$$
Y_-=Y_0=\{\eta\in Y:\eta(0)=0, \ \eta(T)=0\}.
$$

\begin{corollary}
\label{cor7.1}
Let $m$ be odd and $S = \mathrm{id}$.
If $\gamma_+$ is a nondegenerate minimum of $\mathcal{A}_+$,
then the corresponding reversible periodic orbit $\gamma$ has a real multiplier $>1$.
\end{corollary}

Indeed, $\gamma$ preserves orientation, so $\sigma>0$. The Hessian
$K_+$ is positive definite, and hence the same is true for
$K_-=K_+\big|_{Y_0}$.

An analogue of Corollary~\ref{cor7.1}
is true also for $m$ even.

\begin{proposition}
\label{pr7.2}
If   $S=\operatorname{id}$ and  $\gamma_+$ is
 a nondegenerate minimum of $\mathscr{A}_+$, then $h_\rho$ is positive definite
for $|\rho|=1$. Hence there are no multipliers on the unit circle.
\end{proposition}

\noindent{\it Proof.} For complex $\xi\in X$, set
$$
\eta(t)=e^{\mu t}\xi(t)=u(t)+iv(t),\qquad
u(t),v(t)\in E_t.
$$
Then
$$
h_\rho(\xi,\bar\xi\,)=h(\eta,\bar\eta)=h(u,u)+h(v,v).
$$
Let us show that $h(u,u)$ is positive definite on
$W^{1,2}([-T,T],E)$.   Indeed,  $u_+=u\big|_{[0,T]}\in Y$ and
$u_-=J(u\big|_{[-T,0]})\in Y$. Thus
$$
h(u,u)=K(u_+,u_+)+K(u_-,u_-)>0,\qquad
u\ne 0.
$$
By Hill's formula       $\det(P-\rho I)\ne 0$ for $\rho\in S^1$, and so  $P$
has no multipliers on  $S^1$. For another proof see \cite{35}.
\qed \medskip

Consider again the case of a general involution $S$. Then  $Y_+\cap
Y_-=Y_0$.  Let $Y_\pm^\perp$ be  the $K_\pm$-complement  of $Y_0$
in $Y_\pm$, that is,  the set of
$\eta\in Y_\pm$ such that $K_\pm(\eta,\zeta)=0$ for all $\zeta\in Y_0$.
By integration by parts,
$$
Y_\pm^\perp=Y_\pm\cap Y^\perp,\qquad
Y^\perp=\{\eta\in Y:(-D^2+U)\eta=0\}.
$$                   The
restriction of  $K_\pm$ to $Y_\pm^\perp$ equals
\begin{equation}
K_\pm^\perp(\eta,\eta)=(D\eta,\eta)\big|_0^T,\qquad
\eta\in Y_\pm^\perp.
\label{eq7.1}
\end{equation}
Let   $K_0=K\big|_{Y_0}$ and~$K^\perp=K\big|_{Y^\perp}$. Then we have
$$
H_\pm\cong K_0\oplus K_\pm^\perp,\qquad
H\cong K_0\oplus K_0\oplus K_+^\perp\oplus K_-^\perp.
$$
It follows that if        $\det K_0\ne 0$, then
$\operatorname{sign}\det H=(-1)^{\operatorname{ind} K_+^\perp+\operatorname{ind} K_-^\perp}=(-1)^{\operatorname{ind} K^\perp}$.

\begin{proposition}
\label{pr7.3}
If the time moments  0 and~$T$ are non-conjugate, then
$$
(-1)^{\operatorname{ind}\gamma}=\operatorname{sign}\det H=(-1)^{\operatorname{ind} K^\perp}.
$$
\end{proposition}

If the time moments  0 and~$T$ are non-conjugate,
$\dim Y_+^\perp=2n$ and
  $\dim Y_-^\perp=2(m-n)$, where  $n=\dim N$.

The  quadratic form $K_+^\perp$
has a simple meaning,   the Hessian of the discrete Lagrangian
(Hamilton action function) defined locally as
$$
L(x,y)=\mathscr{A}_+(\nu),
$$
where  $\nu\colon [0,T]\to M$ is a trajectory joining  $x$ and  $y$.

\addtocontents{toc}{\string\aaaaa}

\appendix

\renewcommand{\thesection}{\Alph{section}}
\renewcommand{\thesubsection}{A.\arabic{subsection}}
\renewcommand{\theequation}{A.\arabic{equation}}
\renewcommand{\thelemma}{A.\arabic{lemma}}

\section{Appendix}
\label{secA}
\addtocontents{toc}{\string\bbbbb}

\subsection{Proof of Theorem \ref{th3.3}}
\label{ssecA.1}
\ref{th3.3a} By  \eqref{eq3.37}, for any constant vector
$\lambda^\beta$ we have:
$$
I^\alpha(\lambda^\beta\hat q_{\beta i},
\lambda^\beta\hat q_{\beta\,i+1})=c_\beta^\alpha\lambda^\beta.
$$
By condition  \ref{conC}, the matrix  $c_\alpha^\beta$ is nondegenerate.
Therefore the equation
$$
I^\alpha(\lambda^\beta\hat q_{\beta i},
\lambda^\beta\hat q_{\beta\,i+1})=r^\alpha
$$
with respect
to $\lambda^\beta$ is solvable for any constant vector  $r^\alpha$.
This implies the first statement in~\ref{th3.3a}.

To prove the second statement in  \ref{th3.3a}, we show that
$\dim\Omega^\perp=k$ and $\Omega^\perp\cap X_0^\perp= 0$. In view
of equation~\eqref{eq3.26}, the last two conditions are equivalent to the
non-degeneracy of the matrix
$$
e_{\alpha\beta}=\langle\mathbf{d}_\alpha,\mathbf{q}_\beta^\perp\rangle
=\sum_i g_{\alpha\delta i}I_i^\delta\bigl(q^\perp_{\beta i},q^\perp_{\beta\,i+1}\bigr).
$$
By \eqref{eq3.36} we have
\begin{align*}
e_{\alpha\beta}&=\sum_i g_{\alpha\delta i}
\bigl(-\langle B_i q^\perp_{\beta i},w_{i+1}^\delta\rangle
+\langle B_i w_i^\delta,q^\perp_{\beta\,i+1}\rangle\bigr)
\\
&=\sum_i g_{\alpha\delta i}\bigl(\delta_\beta^\delta
+\lambda^\perp_{\beta\varepsilon i}g_i^{\varepsilon\delta}
-\lambda^\perp_{\beta\varepsilon\,i+1} g_i^{\delta\varepsilon}\bigr)
\\
&=\bar g_{\alpha\beta}-\sum_i\bigl(\lambda^\perp_{\beta\varepsilon\,i+1}
-\lambda^\perp_{\beta\varepsilon i}\bigr)=\bar g_{\alpha\beta}-s_{\alpha\beta}
=-a_{\alpha\beta}^\perp.
\end{align*}

\ref{th3.3b} Since the map $\widehat\Pi\big|_{\widehat{Y}}\colon
\widehat Y\to X^\perp$ is an isomorphism,  $\widehat\Pi\big|_{\widehat\Omega}
\colon \widehat\Omega\to \Omega^\perp$ and $\widehat\Pi\big|_{\widehat Y^0}
\colon \widehat Y^0\to X_0^\perp$ are isomorphisms.

\ref{th3.3c} By Proposition  \ref{pr6.3}, for any
$\mathbf{v}\in X^\perp$
$$
h^\perp(\mathbf{v},\mathbf{v})=h\big|_Y \circ \Pi^{-1}(\mathbf{v},\mathbf{v})
+\bar g_{\alpha\beta} c^\alpha c^\beta .
$$
If   $\mathbf{v}\in X_0^\perp$, we have $\mathbf{v}=\Pi\mathbf{u}=
(u_i-\lambda_{\alpha i} w_i^\alpha)$, \,$\mathbf{u}\in Y^0$. Then
by Lemma  \ref{lem6.2},
$c^\alpha=I_i^\alpha(\mathbf{u}_i,\mathbf{u}_{i+1})=0$.  This implies the
first equation in  \ref{th3.3c}. To prove the second it is sufficient
to note that $h=\hat h\circ \pi_\Gamma$ and~ $\Pi=\hat\pi
\circ \Pi_\Gamma$.

\ref{th3.3d} By \eqref{eq2.10},
\begin{align*}
h(\widehat{\mathbf{q}}_\alpha,\mathbf{u})&= \sum_i \langle
A_i\hat q_{\alpha i}-B^*_i\hat q_{\alpha\,i+1}
-B_{i-1}\hat q_{\alpha\,i-1},u_i \rangle
\\
&=\sum_i \langle-A_i \nu_{\alpha\beta i} w_i^\beta
+B^*_i \nu_{\alpha\beta i+1}w_{i+1}^\beta
+B_{i-1} \nu_{\alpha\beta\,i-1} w_{i-1}^\beta,u_i \rangle
\\
&=\sum_i \langle B^*_i \Delta\nu_{\alpha\beta i} w_{i+1}^\beta
-B_{i-1} \Delta\nu_{\alpha\beta\,i-1} w_{i-1}^\beta,u_i \rangle
\\
&=s_{\alpha\gamma}\kappa^{\gamma\delta}\sum_i
\langle B^*_i g_{\delta\beta i}w_{i+1}^\beta
-B_{i-1} g_{\delta\beta\,i-1} w_{i-1}^\beta,u_i \rangle=
-s_{\alpha\gamma}\kappa^{\gamma\delta}\sum_i
\langle d_{\delta i},u_i\rangle,
\end{align*}
where the~$d_{\delta i}$ are defined by  \eqref{eq3.25}.
Proposition~\ref{pr3.9} implies the first assertion in \ref{th3.3d}.
Analogously,
\begin{align*}
h(\widehat{\mathbf{q}}_\alpha,\widehat{\mathbf{q}}_\beta)&=
s_{\alpha\gamma} \kappa^{\gamma\delta}\sum_i \langle B^*_i
g_{\delta\varepsilon i} w_{i+1}^\varepsilon
-B_{i-1}g_{\delta\varepsilon\,i-1} w_{i-1}^\varepsilon, \widehat
q_{\beta i}\rangle
\\
&=s_{\alpha\gamma} \kappa^{\gamma\delta}\sum_i g_{\delta\varepsilon i}
\bigl(\langle B_i\hat q_{\beta i},w_{i+1}^\varepsilon\rangle
-\langle B_i w_i^\varepsilon,\hat q_{\beta\,i+1}\rangle\bigr)
\\
&=-s_{\alpha\gamma}\kappa^{\gamma\delta}\bar g_{\delta\varepsilon}
c_\beta^\varepsilon=-s_{\alpha\gamma}c_\beta^\gamma=a_{\alpha\beta}.
\end{align*}
This implies the second assertion in \ref{th3.3d}.

\ref{th3.3e} We get
\begin{align*}
h^\perp(\mathbf{q}^\perp_\alpha,\mathbf{u}^\perp)&=\sum_i
\langle(A_i-C_i) q^\perp_{\alpha i}-B_i^* q^\perp_{\alpha\,i+1}
-B_{i-1} q^\perp_{\alpha\,i-1},u^\perp_i \rangle
\\
&=-\sum_i \langle C_i q^\perp_{\alpha i}+\lambda^\perp_{\alpha\beta i}
A_i w_i^\beta-\lambda^\perp_{\alpha\beta\,i+1} B_i^* w_{i+1}^\beta
-\lambda^\perp_{\alpha\beta\,i-1} B_{i-1} w_{i-1}^\beta,u^\perp_i \rangle
\\
&=\sum_i g_{\gamma\delta i}\bigl(\langle B_i w_i^\gamma,
q_{\alpha\,i+1}\rangle-\langle B_i q_{\alpha\,i},w_{i+1}^\gamma\rangle\bigr)
\langle B_i u^\perp_i,w_{i+1}^\delta \rangle
\\
&=\sum_i g_{\alpha\delta i}\langle B_i u^\perp_i,w_{i+1}^\delta\rangle.
\end{align*}
Since   $u_i^\perp=u_i-g_{\varepsilon\beta\,i-1}
\langle B_{i-1} w_{i-1}^\varepsilon,u_i \rangle v_i^\beta$, we obtain
$$
g_{\gamma\delta\,i}\langle B_i u^\perp_i,w_{i+1}^\delta \rangle
=-\langle d_{\gamma i},u_i \rangle=0.
$$
This implies the first assertion in  \ref{th3.3e}. Analogously,
\begin{align*}
h^\perp(\mathbf{q}^\perp_\alpha,\mathbf{q}^\perp_\beta)
&=\sum_i\bigl(g_{\alpha\delta i}\langle B_i q_{\beta i},w_{i+1}^\delta\rangle
- g_{\alpha\delta\,i-1}\langle B_{i-1}w_{i-1}^\delta,q_{\beta i}\rangle\bigr)
\\
&=\sum_i\bigl(-g_{\alpha\beta i}+g_{\alpha\delta i}
\langle B_i w_i^\delta,q_{\beta\,i+1}\rangle-g_{\alpha\delta\,i-1}
\langle B_{i-1} w_{i-1}^\delta,q_{\beta i}\rangle\bigr)
\\
&=-\bar g_{\alpha\beta}+s_{\alpha\beta}.
\end{align*}
This implies the second assertion in~\ref{th3.3e}.

\subsection{Proof of Theorem \ref{th6.4}}
\label{ssecA.2}
\ref{th6.4a} For any constant vector   $\lambda^\beta$ we have
$I^\alpha(\lambda^\beta\hat\eta_\beta,D\lambda^\beta\hat\eta_\beta)
=\hat c_\beta^{\,\alpha}\lambda^\beta$. therefore for any constant vector
$r^\alpha$ the coefficients~$\lambda^\beta$ can be chosen so that
$I^\alpha(\lambda^\beta\hat\eta_\beta,
D\lambda^\beta\hat\eta_\beta)=r^\alpha$. This implies the
first equation in  \ref{th6.4a}.

To prove the second equation in  \ref{th6.4a}, we show that $\dim\Omega^\perp=k$
and $\Omega^\perp\cap X_0^\perp=0$. In view of equation  \eqref{eq6.24},
it is sufficient to check that the matrix
$$
e_{\alpha\beta}=\int_0^\tau g_{\alpha\delta}
(D^\perp\zeta^\delta,\eta_\beta^\perp)\,dt
$$
is nondegenerate. We have:
\begin{align*}
e_{\alpha\beta}&=\int_0^\tau g_{\alpha\delta}
(D\zeta^\delta,\eta_\beta^\perp)\,dt=
\int_0^\tau g_{\alpha\delta}(D\zeta^\delta,\eta_\beta-
\lambda_{\beta\varepsilon}^\perp\zeta^\varepsilon)\,dt
\\
&=\int_0^\tau\bigl(g_{\alpha\delta}(D\zeta^\delta,\eta_\beta)
-g_{\alpha\delta}(\eta_\beta,\zeta^\nu)g_{\nu\varepsilon}
(D\zeta^\delta,\zeta^\varepsilon)\bigr)\,dt.
\end{align*}
Using the equation
\begin{equation}
\label{eqA.1}
0=\frac {d}{dt}\bigl(g_{\varepsilon\nu} g^{\nu\delta}\bigr)=
\bigl(\dot g_{\varepsilon\nu}\zeta^\nu+
2g_{\varepsilon\nu}D\zeta^\nu,\zeta^\delta\bigr),
\end{equation}
we continue:
\begin{align*}
e_{\alpha\beta}&=\frac{1}{2}\int_0^\tau\biggl(\frac{d}{dt}
\bigl(g_{\alpha\delta}(\zeta^\delta,\eta_\beta)\bigr)
-g_{\alpha\beta}\biggr)\,dt
=\frac{1}{2}\,(s_{\alpha\beta}-\bar g_{\alpha\beta})
=-\frac{1}{2}\,a_{\alpha\beta}^\perp.
\end{align*}

\ref{th6.4b} The maps    $\widehat\Pi\big|_{\widehat\Omega}\colon
\widehat\Omega\to\Omega^\perp$ and  $\widehat\Pi\big|_{\widehat Y^0}\colon
\widehat Y^0\to X_0^\perp$ are isomorphisms because
$\widehat\Pi\colon\widehat Y\to X^\perp$ is an isomorphism.

\ref{th6.4c} By Proposition \ref{pr6.3}, for any     $\eta\in X^\perp$
$$
h^\perp(\eta,\eta)=h\big|_Y \circ \Pi^{-1}(\eta,\eta)
+\bar g_{\alpha\beta} c^\alpha c^\beta.
$$
If   $\eta\in X_0^\perp$, we have
$\eta=\Pi\xi=\xi-\lambda_\alpha\zeta^\alpha$, \,$\xi\in Y^0$. Then by
Lemma  \ref{lem6.1}, $c^\alpha=I^\alpha(\xi,D\xi)=0$.  This implies the
first equation in  \ref{th6.4c}. To prove the second it is sufficient
to note that $h=\hat h\circ \pi_\Gamma$ and  $\Pi=\widehat\Pi
\circ \pi_\Gamma$.

\ref{th6.4d} Integrating by parts we get
\begin{align*}
h(\hat\eta_\alpha,\xi)&=\int_0^\tau(-D^2\hat\eta_\alpha+
U\hat\eta_\alpha,\xi)\,dt
=s_{\alpha\delta}\kappa^{\delta\varepsilon}\int_0^\tau
(\dot g_{\varepsilon\nu}\zeta^\nu+2g_{\varepsilon\nu}D\zeta^\nu,\xi)\,dt
\\
&=s_{\alpha\delta}\kappa^{\delta\varepsilon}\int_0^\tau \frac d{dt}
\bigl(g_{\varepsilon\nu}(\zeta^\nu,\xi)\bigr)\,dt=0.
\end{align*}
Analogously,
$$
h(\hat\eta_\alpha,\hat\eta_\beta)=
s_{\alpha\delta}\kappa^{\delta\varepsilon}\int_0^\tau(\dot g_{\varepsilon\nu}
\zeta^\nu+2g_{\varepsilon\nu}D\zeta^\nu,\hat\eta_\beta)\,dt.
$$
We define~$\mu_{\alpha\beta}$ so that
$\dot\mu_{\alpha\beta}=c_\alpha^\delta g_{\delta\beta}$. Then
$$
\bigl(\zeta^\alpha,D(\hat\eta_\beta-\mu_{\beta\delta}\zeta^\delta)\bigr)
-(\hat\eta_\beta-\mu_{\beta\delta}\zeta^\delta,D\zeta^\alpha)
=c^\alpha_\beta-\dot\mu_{\beta\delta}g^{\delta\alpha}=0.
$$
Therefore,
\begin{align*}
\int_0^\tau(\dot g_{\varepsilon\nu}\zeta^\nu+
2g_{\varepsilon\nu}D\zeta^\nu,\hat\eta_\beta)\,dt
&=\int_0^\tau\frac{d}{dt}\bigl(g_{\varepsilon\nu}
(\zeta^\nu,\hat\eta_\beta-\mu_{\beta\delta}\zeta^\delta)\bigr)\,dt
\\
&=-\int_0^\tau \mu_{\varepsilon\beta}\,dt
=-\bar g_{\varepsilon\delta}c^\delta_\beta,
\end{align*}
where we have used \eqref{eqA.1}. Finally,
$\hat h(\hat\eta_\alpha,\hat\eta_\beta)=-s_{\alpha\delta}
\kappa^{\delta\varepsilon}\bar g_{\varepsilon\nu} c^\nu_\beta
=-s_{\alpha\delta} c^\delta_\beta=a_{\alpha\beta}$.

\ref{th6.4e} Integrating by parts we get
$$
h^\perp(\eta^\perp_\alpha,\xi^\perp)=\int_0^\tau
\Bigl(\bigl((-D\Pi D+U)\eta^\perp_\alpha,\xi^\perp\bigr)
-3g_{\delta\varepsilon}(\eta^\perp_\alpha,D\zeta_\delta)
(\xi^\perp,D\zeta_\varepsilon)\Bigr)\,dt.
$$
Direct, but lengthy computation gives
$$
h^\perp(\eta^\perp_\alpha,\xi^\perp)=\int_0^\tau 2g_{\alpha\varepsilon}
(\xi^\perp,D\zeta^\varepsilon)\,dt=\int_0^\tau\frac{d}{dt}
\bigl(g_{\alpha\varepsilon}(\xi,\zeta^\varepsilon)\bigr)
=0.
$$
Here we have used  \eqref{eqA.1}. Analogously,
\begin{align*}
h^\perp(\eta^\perp_\alpha,\eta^\perp_\beta)
&=\int_0^\tau 2g_{\alpha\varepsilon}(D\zeta^\varepsilon,\eta^\perp_\beta)\,dt
=\int_0^\tau 2g_{\alpha\varepsilon}\bigl((D\zeta^\varepsilon,\eta_\beta)
-(D\zeta^\varepsilon,\zeta^\nu)
g_{\nu\vartheta}(\zeta^\vartheta,\eta_\beta)\bigr)\,dt
\\
&=\int_0^\tau\biggl(-g_{\alpha\beta}+\frac{d}{dt}
\bigl(g_{\alpha\vartheta}(\zeta^\vartheta,\eta_\beta)\bigr)\biggr)\,dt
=-\bar g_{\alpha\beta}+s_{\alpha\beta}.
\end{align*}

\subsection{Degenerate case}
\label{ssecA.3}
In this subsection we consider the case when the nonde\-generacy
assumption  \ref{nconA} on p.~\pageref{nconA} fails, that is,
$\operatorname{rank}(\zeta^1(t),\dots,\zeta^k(t))$ drops on
$\Sigma\subset\mathbb{R}/\tau\mathbb{Z}$.
We will see that the Routh reduction of the
system $(E,\Lambda)$ to $(E^\perp,\Lambda^\perp)$
and other results on elimination of degeneracy
hold with minor modifications of the proofs. Note that for DLS
condition \ref{conA}, which is similar to
\ref{nconA}, is probably necessary.

\begin{lemma}
\label{lemA.1}
The family $(F_t)_{t\not\in\Sigma}$ can be   extended to a smooth
vector bundle $F=(F_t)_{t\in\mathbb{R}/\tau\mathbb{Z}}$.
Thus the orthogonal complement
$E^\perp=(E_t^\perp)$ is a smooth vector bundle. The operator $D^\perp$ on
$E^\perp$ and the reduced
Lagrangian $\Lambda^\perp$ defined for $t\notin\Sigma$ can be
smoothly extended to $t\in\Sigma$.
\end{lemma}

\noindent{\it Proof.}
Suppose that     $0\in\Sigma$, and let
$$
\operatorname{rank}\bigl(\zeta^1(0),\dots,\zeta^k(0)\bigr)=k-l.
$$
Without loss of generality we may assume that
$$
\zeta^1(0)=\dots=\zeta^l(0)=0,\qquad
\operatorname{rank}\big(\zeta^{l+1}(0),\dots,\zeta^k(0)\big)=k-l.
$$
Then  $\operatorname{rank}\big(D\zeta^1(0),\dots,D\zeta^k(0)\big)=l$,
or else the solutions $\zeta^1,\dots,\zeta^l$ of the variational system are
dependent. Since
$$
I^\alpha(\zeta^\beta,D\zeta^\beta)=(\zeta^\alpha,D\zeta^\beta)-
(\zeta^\beta,D\zeta^\alpha)=0,
$$
we have
$$
\bigl(D\zeta^\alpha(0),\zeta^\beta(0)\bigr)=0, \qquad
\alpha=1,\dots,l,\quad
\beta=l+1,\dots,k.
$$
Thus
$$
\operatorname{rank}\bigl(D\zeta^1(0),\dots,D\zeta^l(0),
\zeta^{l+1}(0),\dots,\zeta^k(0)\bigr)=k.
$$
Since
$$
D(\zeta^\alpha-tD\zeta^\alpha)=-tD^2\zeta^\alpha=-tU\zeta^\alpha=O(t^2),
$$
we have
$$
\zeta^\alpha(t)=t D\zeta^\alpha(t)+O(t^3),\qquad
\alpha=1,\dots,l.
$$
Thus, the space
$$
F_t=\operatorname{span}\biggl(\frac{1}{t}\,D\zeta^1(t),\dots,
\frac{1}{t}\,D\zeta^l(t),
\zeta^{l+1}(t),\dots,\zeta^k(t)\biggr),\qquad t\ne 0,
$$
has a limit
$$
F_0=\operatorname{span}\bigl(D\zeta^1(0),\dots,D\zeta^l(0),
\zeta^{l+1}(0),\dots,\zeta^k(0)\bigr)
$$
as  $t\to 0$, and         $F_t$ is smooth at $t=0$.
The first statement is proved.

Since $E_t^\perp$ is smooth at $t=0$,
also        $\Pi_t\colon E_t\to E_t^\perp$ is smooth,
and hence the operator~$D^\perp$ is smooth. Finally we need to check that
the  term
$$
(C\eta,\eta)=g_{\alpha\beta}(u,D^\perp\zeta^\alpha)(u,D^\perp\zeta^\beta)
$$
in~$\Lambda^\perp$ is smooth at  $t=0$.

Denote    $j^{\alpha\beta}=\bigl(D\zeta^\alpha(0),D\zeta^\beta(0)\bigr)$,
\,$\alpha,\beta=1,\dots,l$. Then  $(j^{\alpha\beta})$
is a nondegenerate matrix and
$$
(D\zeta^\alpha,D\zeta^\beta)=j^{\alpha\beta}+O(t^2).
$$
Thus
$$
g^{\alpha\beta}=(\zeta^\alpha,\zeta^\beta)=t^2(D\zeta^\alpha,D\zeta^\beta)+
O(t^4)=t^2j^{\alpha\beta}+O(t^4),\qquad
\alpha,\beta=1,\dots,l.
$$                        The
matrix  $(g^{\alpha\beta}(0))$, \,$\alpha,\beta=l+1,\dots,k$, is nondegenerate,
while one can show that
$$
g^{\alpha\beta}=O(t^2),\qquad
\alpha=1,\dots,l,\quad
\beta=l+1,\dots,k.
$$
Thus for the inverse matrix $(g_{\alpha\beta})$ we obtain
$$
g_{\alpha\beta}=t^{-2}(j_{\alpha\beta}+O(t^2)),\qquad
\alpha,\beta=1,\dots,l.
$$
The block $(g_{\alpha\beta})$, \,$\alpha,\beta=l+1,\dots,k$,
is smooth and nondegenerate and the
block  $(g_{\alpha\beta})$, \,$\alpha=1,\dots,l$,
\,$\beta=l+1,\dots, k$, is smooth.

Since   $D^\perp\zeta^\alpha(t)=O(t^2)$ for $\alpha=1,\dots,l$,
we obtain that $C$ is smooth at $t=0$.
Thus the reduced Lagrangian $\Lambda^\perp$ is smooth on $E^\perp$.
\qed \medskip

In fact, everything we have done in \S\,\ref{sec6}
holds in the singular case. For example, let us check that the projection
 $\Phi\colon X\to \widehat Y$
along $Z$ is well defined and smooth. As in the nondegenerate case, we
 have
$$
Z=\{\xi\in X:\xi(t)\in F_t\}=\biggl\{\xi(t)=\lambda_\alpha(t)\zeta^\alpha(t):
\int_0^\tau g^{\alpha\beta}\dot\lambda_\alpha\dot\lambda_\beta\,dt
<\infty\biggr\},
$$
but now  $\lambda_\alpha$ may be singular for           $t\in\Sigma$.

Take   $\eta\in X$ and look for  $\lambda_\alpha(t)$ such that
$\xi=\eta+\lambda_\alpha\zeta^\alpha \in Y$. Then by \eqref{eq6.8} and
  \eqref{eq6.23},
$$
\dot\lambda_\alpha=g_{\alpha\beta} c^\beta-
g_{\alpha\beta}I^\beta(\eta,D\eta)=g_{\alpha\beta}c^\beta+
\frac{d}{dt}\bigl(g_{\alpha\beta}(\eta,\zeta^\beta)\bigr)
-2g_{\alpha\beta}(\eta,D^\perp\zeta^\beta).
$$
The last term is smooth for $t\in\Sigma$. Suppose again that $0\in\Sigma$.
Then we obtain
$$
\lambda_\alpha(t)=-\frac{1}{t}\sum_{\beta=1}^l j_{\alpha \beta} c^\beta
+\frac{1}{t}\sum_{\beta=1}^l j_{\alpha \beta}(\eta,\zeta^\beta)
+\text{smooth terms}, \qquad
\alpha=1,\dots,l.
$$
It follows that $\lambda_\alpha(t)\zeta^\alpha(t)$ is smooth,
so $\xi(t)$ is smooth.


\end{document}